\documentclass[11pt]{amsart}

\usepackage{amsmath}
\usepackage{amssymb}
\usepackage{amsfonts}
\usepackage{amsthm}
\usepackage{verbatim}
\usepackage{amscd}
\usepackage{cite}
\usepackage{leftidx}
\usepackage{enumerate}
\usepackage{txfonts}
\usepackage{manfnt}
\usepackage{amscd}
\usepackage{bookman}
\usepackage[mathscr]{eucal}
\usepackage{hyperref}
 \DeclareMathOperator{\ts}{TS}

 \DeclareMathOperator{\dsc}{Dsc}
 
 \DeclareMathOperator{\TT}{TT}
 \DeclareMathOperator{\odsc}{ODsc}
\DeclareMathOperator{\Bl}{\mathcal {B}\l}
 
\def\inv{^{-1}}

\DeclareMathOperator{\mgbar}{\overline\M_g}
\DeclareMathOperator{\del}{\partial}
 \DeclareMathOperator{\Hom}{Hom}

\def\refp #1.{(\ref{#1})}

\newcommand{\A}{\mathcal{A}}

\newcommand{\M}{\mathcal{M}}

\newcommand{\bE}{\mathbb{E}}

\newcommand{\bR}{\mathbb{R}}

\newcommand{\TS}{\mathrm{TS}}

\newcommand{\Cal}[1]{\mathcal #1}

\newcommand{\un}{{\underline{n}}}
\newcommand{\ul}[1]{\underline {#1}}

\def\sbr #1.{^{[#1]}}
\def\sfl #1.{^{\lfloor #1\rfloor}}

\newcommand{\subpar}[1]{_{(#1)}}

\newcommand\strad{{\mathrm {strad}}}
\newcommand\subp [1] {_{(#1)}}

\newcommand\half{\frac{1}{2}}
\newcommand\red{{\mathrm red}}

\def\bt{\boxtimes}
\def\what{\widehat}
\def\inv{^{-1}}
\def\?{{\bf{??}}}

\def\Hilb{\text{Hilb}}

\def\M{\Cal M}
\def\A{\Bbb A}

\def\C{\mathbb C}
\def\P{\mathbb P}
\def\N{\mathbb N}
\def\R{\mathbb R}
\def\Z{\mathbb Z}

\def\sgn{\text{\rm sgn} }
\def\sym{\text{\rm Sym} }

\def\Spec{\text{\rm Spec} }

\def\ls{\vskip.25in}

\def\Q{\mathbb Q}

\def\O{\mathcal O}

\def\bt{\boxtimes}

\def\Sym{\textrm{Sym}}

\def\g{\mathfrak g}

\def\1/2{\frac{1}{2}}

\def\I{\mathcal{ I}}

\def\simto{\stackrel{\sim}{\rightarrow}}
\def\2{{[2]}}
\def\l{\ell}
\def\nl{\newline}

\def\he{\mathcal{HE}}

\def\<{\langle}
\def\>{\rangle}

\def\2{{[2]}}
\def\l{\ell}
\def\ti{\mathrm{TI}}

\def\scl #1.{^{\lceil#1\rceil}}
\def\spr #1.{^{(#1)}}
\def\sbc #1.{^{\{#1\}}}

\def\subpr#1.{_{(#1)}}

\def\beq{\begin{equation*}}
\def\eeq{\end{equation*}}
\newcommand{\td}{\tilde }

\def\g3{{\Gamma\spr 3.}}
\def\gg{{\Gamma\spr 2.}}
\def\ggg{{\Gamma\spr 3.}}
\def\lom{L.\omega}
\def\omsq{\omega^2}

\newcommand{\beql}[2]{\begin{equation}\label{#1}#2\end{equation}}
\newcommand{\beqa}[2]{\begin{eqnarray}\label{#1}#2\end{eqnarray}}
\newcommand{\eqspl}[2]{
\begin{equation}\label{#1}
\begin{split}
#2\end{split}\end{equation}}
\newcommand{\eqsp}[1]{\begin{equation*}
\begin{split}#1\end{split}\end{equation*}}

\newcommand{\exseq}[3]{
0\to #1\to #2\to #3\to 0
}
\newcommand{\beginalphaenum}{
\begin{enumerate}\renewcommand{\labelenumi}{ }
\item \begin{enumerate}
}

\def\eex{\end{rm}\end{example}}
\newcommand\newsection[1]{\section{#1}\setcounter{equation}{0}
}
\newcommand\newsubsection[1]{\subsection{#1}\setcounter{equation}{0}}

\newcommand{\be}{\mathbb E}

\newtheorem{thm}{Theorem}[section]

\newtheorem*{thm*}{Theorem}
\newtheorem{cor}[thm]{Corollary}
\newtheorem*{cor*}{Corollary}

\newtheorem{lem}[thm]{Lemma}

\newtheorem*{claim*}{Claim}
\newtheorem{prop}[thm]{Proposition}

\newtheorem{defn}[thm]{Definition}
\theoremstyle{remark}

\newtheorem{rem}[thm]{Remark}

\newtheorem{example}[thm]{Example}
\newtheorem*{example*}{Example}

\newcommand{\mm}{m,m-1}
\begin{document}
\title{ Tautological module
and intersection theory\\
on Hilbert schemes of nodal curves
}
\normalsize
\author
{Ziv Ran}
\thanks{\raggedright{
arxiv.org 0905.2229 This article is a revision of portions of arXiv:0803.4512}}
\date {\today}
\address {Math Dept.  UC Riverside\nl
Surge Facility,  Big Springs Road,\nl
Riverside CA 92521}
\email {ziv.ran @ucr.edu} 

\begin{abstract}
This paper presents the rudiments of Hilbert-Mumford Intersection (HMI) theory: intersection theory on the relative Hilbert scheme of
 a family of nodal (or smooth)
curves, over a base of arbitrary dimension. We introduce an additive
group of geometric cycles, called 'tautological module', generated by diagonal loci, node
scrolls, and twists thereof. We determine recursively the intersection action on
 this group by the discriminant (
big diagonal) divisor and all its powers.
We show that this
suffices to determine arbitrary polynomials in Chern classes, in particular
Chern numbers, for the  tautological vector bundles on the Hilbert
schemes, which are closely related to enumerative geometry of families of nodal curves.

\end{abstract}
\maketitle
\tableofcontents

\setcounter{section}{-1} \newsection{Overview}
This paper is a contribution to the study of Enumerative Geometry
of nodal curves via their subschemes.
To illustrate informally, in part, what it is about, we recall a formula
from  19th century Algebraic Geometry (see
for example \cite{semple-roth}, p. 377): a
nonsingular curve  $X$
in complex projective 3-space admits an expected
 finite number $n_4$ of  4-secant lines,
and that a formula for $n_4$ in terms of   the
degree $d$ and genus $g$ of $X$ can be given: specifically,
\[24n_4=144-204d+106d^2-12d^2g-24d^3+2d^4+84dg-156g+12g^2.\]
Standard dimension-counting suggests that given a nice enough
$b$-dimensional family of space curves, it will admit a finite
number $n_{4+b}$ of $(4+b)$-secant lines, and one can ask for
a method to compute $n_{4+b}$ in terms of basic projective characters
of the family. The results of this paper provide, inter alia,
such a method, as we now proceed to describe.
\newsubsection{Setting} To fix ideas, consider a family of curves given by a flat projective
morphism \beq \pi:X\to B \eeq over an irreducible base, with fibres \beq
X_b=\pi\inv(b), b\in B \eeq which are irreducible nonsingular for the
generic $b$ and   at worst nodal for every $b$. For example, $X$ could be
the universal family of automorphism-free curves over the appropriate
open subset of $\overline{\mathcal M}_g$, the moduli space of
Deligne-Mumford stable curves. Many questions in the classical projective
and enumerative geometry of this family can be naturally phrased, and in a
formal sense solved (see for instance \cite{R}), in the context of the
{\it{relative Hilbert scheme}} \beq X\sbr m._B=\Hilb_m(X/B). \eeq This is a universal
parameter space for length-$m$ subschemes of $X$ contained in fibres of $\pi$,
and carries natural {\it {tautological vector bundle}} $\Lambda_m(E)$,  associated to any vector bundle $E$ on $X$ (e.g., the relative dualizing sheaf $\omega_{X/B}$). One specific example of the enumerative questions
which may be considered from this viewpoint is  the fundamental class in $\mgbar$ of the closure of the hyperelliptic locus.\par Typically, the
geometric questions one wants to consider can be formulated in terms of relative multiple points and multisecants
in the family, which can be described in terms of degeneracy
loci of bundle maps involving tautological bundles,  and the formal solutions involve
Chern numbers of those
tautological bundles. Thus, turning these formal solutions into meaningful
ones requires computing the Chern numbers in question. This problem
was stated but, aside from some low-degree cases,  left open in \cite{R}.
Our main purpose here is to solve this problem in general. More than that,
we shall in fact provide a calculus to compute certain images of arbitrary polynomials in the
Chern classes of the tautological bundles. In the 'absolute' case
$E=\omega_{X/B}$, the computation ultimately reduces these polynomials
to polynomials in Mumford's tautological classes \cite{Mu} on various
boundary strata of $B$. The latter are computed via a conjecture of Witten,
proved by Kontsevich \cite{Kon}. It should be mentioned
that in the case of
the symmetric product (= Hilbert scheme) of a single
smooth curve, a complete intersection theory was worked out by
Macdonald \cite{macd}. On the other hand,
the intersection theory of Hilbert schemes of smooth
\emph{surfaces} was investigated deeply by Nakajima, Lehn and
others, see \cite{nak}, \cite{L}, \cite{lehn-montreal} and references therein.
\newsubsection{Tautological module: motivation}
Now the framework for our solution is a little different to what
is commonly done in similar problems (e.g. Macdonald's set-up). Rather than compute a suitable
intersection ring, we will focus primarily on the (intersection)
action of the \emph{discriminant} or big diagonal $\Gamma\spr m.$
and its powers.
The motivation for this approach
comes from a result in \cite{R} called the 'Splitting principle'.
This says that the total Chern class of pullback  of a tautological
bundle $\Lambda_m(E)$ to the full-flag Hilbert scheme $W^m=W^m(X/B)$, which maps to the degree-$i$ Hilbert schemes
$X\sbr i._B, i\leq m$,  can
be expressed as a simple decomposable polynomial in the
(pullbacks of)  $\Gamma\spr
i., i\leq m$. The recursive analogue of this result, Cor. \ref{tautbun}
below, says
that the pullback of $c(\Lambda_m(E))$ on the
  Hilbert scheme $X\sbr\mm._B$, parametrizing
flags of schemes of lengths $m,m-1$ (which we will call the 'flaglet' Hilbert scheme), is a product of
$c(\Lambda_{m-1}(E))$ and a polynomial in in discriminants
$\Gamma\spr i., i\leq m$.  It
follows that if we assume recursively that we have some
 reasonable way to  express polynomials in
$c(\Lambda_{m-1}(E))$, say as elements of a 'Tautological module'
 $T^{m-1}_R(X/B)$ and want to do
the same for $m$, then we need to determine 2 things:
\par 1. Tautological {module} in degree $m$, $T^m_R(X/B)$, i.e. a group
 together with an action of $\Gamma\spr m.$.
 \par 2. Transfer
calculus, going from $T^{m-1}_R(X/B)$ to $T^m_R(X/B)$ via the
flaglet correspondence $X\sbr m,m-1._B$. \par Given these,
$T^m_R(X/B)$ would recursively contain all polynomials
in $\Gamma\spr i., i\leq m$, hence all polynomials in the Chern classes
of $\Lambda_m(E)$.
\newsubsection{Tautological module: elements} Given a family $X/B$
of nodal (possibly pointed) curves, the associated Tautological Module
$T^m_R(X/B)$ (Definition \ref{taut-mod-def}) is constructed recursively
in $m$,
\emph{grosso modo}, as follows (see the body of the paper for details). \par
--For $m=1$, it equals $R$, a
$\Q$-subalgebra of $H^\bullet(X)$ containing the relative canonical class
$\omega$ as well as any distinguished sections.
Here $H^\bullet$ denotes any cohomology ring coarser than (i.e. admitting a map from) the Chow ring over $\Q$.
\par -- The recursive
 step. First, decompose the tautological module according to partitions
or 'distributions':
\[ T^m_R(X/B)=\bigoplus\limits_\mu T^\mu_R(X/B)\] the sum being over all
partitions $\mu$ of weight $m$; thus, it suffices to describe
each $\mu$ summand. Then, we parametrize the boundary by a
union of families $T(\theta)$ associated to the relative nodes $\theta$ of
$X/B$, and for each of those let $X^\theta/T(\theta)$ be the corresponding
family blown up in $\theta$, which is endowed with a pair of distinguished
sections denoted $\theta_x, \theta_y$, set $R^\theta=R[\theta_x,
\theta_y]$, and define firstly the \emph{boundary tautological module}
of type $\mu$ as
\[ \del T^\mu_{R}=\bigoplus\limits_\theta T^\mu_{R^\theta}(X^\theta/T(\theta)) \]
(using recursion, we may assume this defined for $\mu$ of weight $<m$). Then define for $\mu$ of weight $m$,
 \eqspl{sectors}{
T^\mu_R(X/B)=\left( \ts_\mu(R)\right )\oplus\left(
\bigoplus\limits_{\substack{\nu\amalg\{n\}=\mu\\ 0<j<n}}
\left (\Q F^n_j\oplus\Q\Gamma\spr
m.F^n_j\right)\otimes  \del T^\nu_R(X/B)\right)} in which
\begin{itemize}\item

 $\ts_\mu(R)$, the interior part of the module,
 is of a purely topological character and can be identified
  with a formal algebraic construct,
an appropriate summand of the 'tensor-symmetric' algebra
$T(\sym(R))$,
\item $F^n_j$ is a formal symbol (for now), called a 'node scroll',
\item
$\Gamma\spr m.$ is the discriminant or big diagonal on $X\sbr m._B$, for
the purpose of the formula just a formal symbol as well,\item
$-\Gamma\spr m.F^n_j$ is called a \emph{node section}.\item  We call the two
main summands of \eqref{sectors}
 the \emph{diagonal} and \emph{node scroll} sectors of the
tautological module $T^\mu_R$ and denote them $DT^\mu_R,
NT^\mu_R$ and similarly $DT^m_R, NT^m_R$. $NT$ itself splits as
$NFT\oplus N\Gamma T$, node scrolls plus node sections.\end{itemize}
 The above definition is doubly recursive in
the sense that modulo the relatively elementary part $DT^m_R$,
the remaining part $NT^m_R$ involves
tautological modules of lower weight for (boundary) families of lower
genus (albeit with more markings). The recursive definition may
be replaced by a non-recursive one by working with
node \emph{polyscrolls}, associated to a boundary stratum defined by a
collection of nodes rather than a single one.\par
The tautological module
maps to the homology (Chow or ordinary) of the Hilbert scheme, where the
diagonal sector maps to cycles living on various diagonal loci (lifted from
analogous loci on the symmetric product), and the node scroll sector maps
to cycles on certain $\P^1$-bundles which live over the boundary and are
exceptional for the cycle map. In particular, a zero-dimensional or 'top
degree' element $\alpha\in T^m_R(X/B)$ has a well-defined \emph{cycle
degree} or 'integral' $\int\alpha\in\Q$.

\newsubsection{Tautological module: Discriminant action} Now our first main result,
the \emph{Tautological module theorem} \ref{taut-module}, describes the
action of $\Gamma\spr m.$, i.e. the $\Q[\Gamma\spr m.]$-module structure, on the $\Q$-vector space $T^m_R=NT^m_R\oplus DT^m_R $.  This
structure is an extension\eqspl{taut-ext}{
\exseq{NT^m_R}{T^m_R}{DT^m_R}}
where the module structure on the quotient $DT^m_R$, unrelated to the singularities,  is via
standard action of the big diagonal in the cohomology of a symmetric
product (which can be modelled by a second-order differential operator);
the structure on the submodule $NT^m_R$ is by the standard action (via Grothendieck's formula) of a section
$\Gamma\spr m.$ on the cohomology of  a suitable $\P^1$-bundle (and it
therefore anti-triangular with respect to the $NFT\oplus N\Gamma T$
decomposition). It can be described in terms of discriminant actions of
lower weight and lower genus. Also, the 'mixing' part of the action takes
$DT^m_R$ only into the $NFT$ summand of $NT^m_R$.

\newsubsection{Tautological module: transfer}
As indicated above, the story is completed by the
\emph{Transfer Theorem} \ref{taut-tfr}, which computes the transfer
(pull-push) operation on $T^{m-1}_R(X/B)$ via $X\sbr\mm._B$, viewed as
a correspondence between $X\sbr m._B$ and $X\sbr m-1._B$, showing
in particular that it lands in $T^m_R(X/B)$ .\par The conjunction of the
Splitting Principle, Module Theorem and Transfer Theorem computes all
polynomials in the Chern classes, in particular the Chern numbers, of
$\Lambda_m(E)$ as $\Q$-linear combinations of tautological classes on
$X\sbr m._B$.
\newsubsection{Computation}
The calculus of of the discriminant action and Chern polynomials
 has been implemented (for arbitrary base dimension) on the
computer by Gwoho Liu, in the form of a  Java program named Macnodal (in honor of MacDonald \cite{macd}). See \S \ref{macnodal} and \cite{macnodal}  for details.
The results are consistent with Cotteril's \cite{cotteril2} results for
pencils.

\newsubsection{Punctual transfer}
Finally, we discuss an analogue of the tautological module and the
transfer for \emph{punctual schemes}, i.e. those supported at
a single point, which are parametrized by the small diagonal
$\Gamma\subp{m}$, which itself is a (singular)
blowup of $X$. This case is
somewhat simpler in its formal aspects but
still goes to the heart of the complexities of the Hilbert scheme.
It has applications to enumeration of various ramification
loci.
 \newsubsection{Applications }
 A number of applications, examples and computations are
  scattered throughout the paper, especially in \S\S \ref{unordered-sec}, \ref{node-scrolls-sec}, \ref{punctual-transfer-sec} and \ref{examples}.
  In particular, multisecants in nodal families, as
  mentioned at the beginning, are fully enumerated.
In
a less elementary vein,
 the machinery of this paper is
 projected to be the first step of a project
 to compute the fundamental class in $\overline\M_g$ of
the locus of curves admitting a $g^r_d$ for given $r$ and $d$, e.g. a
$g^1_2$ (the hyperelliptics). A baby case (genus 3) can be worked
out here, thanks to the exceptional luxury that the excess
degeneracy is not excessive in dimension. To go further,
the idea is to construct an appropriate boundary modification
of the Hodge bundle, together with its
natural evaluation map to the tautological bundle associated to
 the canonical
bundle, such that the degeneracy locus (in the Hilbert scheme)
of this map would consist of the desired $g^r_d$
locus plus a 'good' excess locus, whose contribution could be
computed by Fulton-MacPherson theory. The required modification
is nontrivial, especially on the $\Delta_0$ boundary component,
and is at present known in detail  only for
$d=2$
(see \cite{grd}). Higher-degree cases are work in progress.\par
I am grateful to Gwoho Liu for many helpful discussions and for creating Macnodal. I also thank Ethan Cotteril for helpful communications about his work, especially \cite{cotteril2},
which provides an alternative method for deriving some of the same enumerative applications in the case of pencils.

\newsection{Preliminaries}
This paper is a continuation of our earlier paper \cite{structure}, whose results, terminology and notations will be used throughout. Some additional terminology and remarks
will be given in this section.
\subsection{Graph enumeration, generating functions}\label{graph} See the textbooks \cite{aigner},
\cite{goulden-jackson} or \cite{stanley} for standard techniques and results.
We will present some variants of known formulas, which will prove useful
in deriving some explicit closed formulas in our Intersection Theory
(see especially \S\ref{node-scrolls-sec}).
\subsubsection{Simple graphs}
We consider connected labelled graphs without loops on a fixed vertex-set $[n]=\{1,...,n\}$. Let $\nu_{n,m}$ be the number of these graphs with $m $ edges, none multiple. We also consider connected 'edge-weighted' graphs on $[n]$, where
each edge $e$ is assigned a positive multiplicity $m(e)$.
 Let $w_{n,m}$ be the weighted number of connected  graphs on $[n]$ where
 the edge multiplicities add up to $m$,
i.e. $\sum\limits_{\mathrm{edges}}m(e)=m$, and
where the
weight of the graph is defined as
$\frac{1}{\prod m(e)!}$. Consider the generating functions
\[T(z,y)=\sum\limits_{n=1}^\infty \sum\limits_{m=0}^\infty\frac{\nu_{n,m}}{n!}z^ny^{m}, T_b(u)=\sum\limits_{n=1}^\infty\frac{\nu_{n,n-1+b}}{n!}u^n.\]
Here $b$ represents the 1st Betti number of the graph.
Because $\nu_{n,m}=0$ for $m<n-1$, we can write
\[T(z,y)=\sum\limits_{b=0}^\infty y^{b-1}T_b(yz)\]

The classical (and elementary) Riddell-Uhlenbeck formula states that
\eqspl{uhlenbeck}{\exp(T(z,y))=1+\sum\limits_{n=1}^\infty \frac{1}{n!}z^n(1+y)^{\binom{n}{2}}
=1+\sum\limits_{n=1}^\infty\frac{1}{n!} \sum\limits_{m=0}^{\binom{n}{2}}\binom{\binom{n}{2}}{m}z^ny^m.}
This follows from the fact that the number of $m$-edge, not necessarily connected graphs on $[n]$ is $\binom{\binom{n}{2}}{m}$.
In particular for the tree case ($b=0$)
we have (Cayley's result)
\eqspl{cayley}{T_0(z)=\sum\limits_{n=1}^\infty \frac{n^{n-2}}{n!}z^n.}
\subsubsection{Edge weighting}
We define similarly
\[W(z,y)=\sum\limits_{b=0}^\infty y^{b-1}W_b(yz)=\sum\limits_{n=1}^\infty
\sum\limits_{m=0}^\infty \frac{w_{n,m}}{n!}z^ny^{m}, W_b(u)=\sum\limits_{n=1}^\infty
\frac{w_{n,n-1+b}}{n!}u^n.\]
Now the total weight of all, not necessarily connected, edge-weighted graphs of total multiplicity $m$ on $[n]$
is
\[\sum\limits_{m_1+...m_{\binom{n}{2}}=m}\frac{1}{m_1!...m_{\binom{n}{2}}!}
=\frac{1}{m!}\binom{n}{2}^m.\]
It follows similarly that
\eqspl{W-function}{
\exp(W(z,y))=1+\sum\limits_{n=1}^\infty\sum\limits_{m=0}^\infty
\frac{\binom{n}{2}^m}{n!m!}z^ny^m.
}
\subsubsection{Vertex weighting}
This formula admits a useful generalization to vertex-weighted
graphs. Suppose vertex $i$ is assigned a fixed weight $p_i$ for
$i=1,...$, where we may set $p_i=0$ for $i>n$). The $p_i$ are
regarded as indeterminates or elements of some $\Q$-algebra
$P$. The entire graph is then weighted
$\prod\limits_{i<j}\frac{(p_ip_j)^{m_{i,j}}}{m_{i,j}!}\in P$,
where $m_{i,j}$ is the multiplicity of the $(i,j)$- edge. Let
$w_{S, p.,m}$ be the total weight of all such graphs which are
connected and have vertex-set $S$ and total edge multiplicity
$m=\sum\limits_{i<j}m_{i,j}$. For a vertex-set $S$ with weights
$p.$, set
\[W_{S, p.}(z,y)=\sum\limits_{m=1}^\infty \frac{w_{S, p.,m}}{|S|!}z^{S}y^m.\]
where $z^{S}=\prod\limits_{i\in S} z_i$. The $z_i$ are formal
variables whose squares are set to zero: $z_i^2=0$; thus
$z^Sz^{S'}=z^{S\cup S'}$ whenever $S\cap S'= \emptyset$ and
otherwise $z^Sz^{S'}=0$. This generating function can be
evaluated as follows. Set $f(S,p.)=\sum\limits_{\substack {
i<j\\ i,j\in S}}p_ip_j$ and note that the total weight of all
such graphs with vertex-set $S$, without the connectedness
hypothesis is
\[\sum\limits_{\sum m_{i,j}=m}\prod\limits_{\substack
{ i<j\\ i,j\in S}}\frac{ (p_ip_j)^{m_{i,j}}}{\prod m_{i,j}!}
=\frac{1}{m!}f(S,p.)^m.\] Then we conclude similarly
\eqspl{w-n-p.-m}{ \exp(\sum\limits_{S} W_{S, p.})
=1+\sum\limits_{S,m} \frac{f(S,p.)^m}{m!|S|!}z^{S}y^m \mod
(z_i^2) } Therefore, \eqspl{}{ \sum\limits_{S} W_{S, p.}=
\sum\limits_{n=1}^\infty \frac{(-1)^{n+1}}{n}(\sum\limits_{S,m}
\frac{f(S,p.)^m}{m!|S|!}z^{S}y^m)^n \mod (z_i^2)
 }
 Note that in this formula, each set $S$ on the RHS will be a
disjoint union of sets $S$ on the LHS.\par

\subsubsection{Directed case}
 We now consider a directed
analogue of the above. Consider forward-directed edge-weighted
graphs on $[m]$. Such a graph is specified by nonnegative
integers $e_{j,i}, \forall j<i$ (= number of edges from $j$ to
$i$), and is assigned a total weight $\frac{1}{\prod
e_{j,i}!}$. We will fix the number of edges into $i$ at $k_i,
i=2,...,m$, i.e. $k_i=\sum\limits_{j<i}e_{j,i}$.
 Let
$\vec w_{m,k.}$ denote the weighted number of these graphs
that are connected. On the other hand, the weighted number
of all these graphs, possibly disconnected, is, as above
\[
\prod\limits_{i=2}^m(\sum\limits_{\sum e_{j,i}=k_i}
\frac{1}{\prod e_{j,i}!})=
\prod\limits_{i=2}^m(\sum\limits_{\sum e_{j,i}=k_i}
\frac{\binom{k_i}{e_{1,i},...,e_{i-1,i}}}{k_i!})=
\prod\limits_{i=2}^m\frac{ (i-1)^{k_i}}{k_i!}. \]
Therefore is we define a generating function in $z, y_2,...$
\[\vec W(z,y_2,...)=\sum\limits_{m=2}^\infty\sum\limits_{k_1,...k_m}\frac{\vec w_{m,k.}}{m!}z^my_2^{k_2}...y_m^{k_m}\]
then it follows that
\eqspl{}{
\exp(\vec W(z, y_2,...))=1+\sum \limits_{m,k.}\frac{1}{m!}\prod\limits_{i=2}^m\frac{ (i-1)^{k_i}}{k_i!}z^my_2^{k_2}...y_m^{k_m}.
} Thus we may consider the $\vec w_{m,l.}$ as known.
\subsection{Products, diagonals, partitions}
\label{partitions}
The intersection calculus we aim to develop is couched
in terms certain diagonal-like loci on products,
defined in the general case in terms of partitions. To facilitate working with
 these loci systematically, we now establish some conventions, notations and
 simple remarks
related to partitions. Our viewpoint on partitions is influenced by the fact that we
will mainly use them to define 'diagonal' conditions, so in particular
singleton blocks are essentially insignificant.\par
\subsubsection{b-partitions} By a \emph{block
partition or b-partition} (aka labelled partition) $(I.)$ of weight $m$ and length $r$ we mean an
expression
\[\{1,...,m\}=I_1\coprod...\coprod I_r, \forall I_j\neq\emptyset.\]
If a b-partition $(I.)$ is such that all its blocks   except for one of them, say $I$, are singletons, we will denote $(I.)$ as  $(I)$ or $(I,m]$. Given a set $X$ (or an object in a category
with products-- the modifications for this case are left to the reader), a
b-partition $(I.)$ of weight $m$ defines an ordered 'polydiagonal' subset
 of the (Cartesian) product $X^m$, which will be denoted by $X^{(I.)}$ or
 $OD\subpar{I.}$ or, if the dependence on $X$ must be explicated, $OD_{(I.), X}$:
 in the case where $X$ is a set,
we identify $X^m$ with the set of functions $\{1,...,m\}\to X$, and then
$X^{(I.)}=OD\subpar{I.}\subset X^m$ consists of the functions constant on each block. It
is the image of an injection $X^r\to X^m$ and will
sometimes also be identified with that injection.\par
\subsubsection{Partitions} A b-partition $(I.)$ determines an ordinary partition
of the same weight, viz. $(|I.|)$, which we prefer to view
via the corresponding 'length distribution'. Thus
the \emph{length distribution} associated to a b-partition $(I.)$ is the
function $\mu:\N\to \Z_{\geq 0}$ defined by
\[n\mapsto |\{j:|I_j|=n\}|\]
We call a function $\mu:\N\to \Z_{\geq 0}$ either a
\emph{distribution} or \emph{partition}. This is the same thing as 'partition' in the usual sense: in the usual partition notation, the partition corresponding to $\mu$ is $(..., n^{\mu(n)}, ...,2^{\mu(2)}, 1^{\mu(1)}).$
A distribution has \emph{weight}
$|\mu|=w(\mu)=\sum n\mu(n)=m$ \emph{degree} $d(\mu)=\sum\limits_{n>0}(n-1)\mu(n)$ and \emph{length}
$\l(\mu)=|\{n:\mu(n)>0\}|$.
 The length distribution of a b-partition of weight $m$
 has weight $m$ and conversely, any distibution
 of weight $m$ is the length distribution of some b-partition of weight $m$.
Two b-partitions are said to be \emph{equivalent} if their distributions are
the same or equivalently, if they differ by a permutation of $1,...,m$. A
distribution is viewed essentially as a collection of block sizes, and will
often be specified by specifying the non-singleton block sizes: e.g. $(n)$
for $n>1$ refers to a distribution (of some weight $m\geq n$) with unique
nonsingleton block of  size $n$. A distribution $\mu$ defines a
polydiagonal or polyblock diagonal
 \eqspl{polyblock-diag-eq}{D_\mu=D_{\mu,X}=\prod\limits_nX\spr \mu(n).\hookrightarrow X\spr
w(\mu).}
where $X\spr k.$ is the $k$th symmetric product. The
embedding is defined by repeating an element in the $n$th factor, i.e.
$X\spr\mu(n).$, $n$ times, i.e.
\[ (\sum\limits_{i=1}^{\mu(n)}x_{i,n}:n=1,2,...)\mapsto
\sum\limits_n \sum\limits_{i=1}^{\mu(n)} nx_{i,n}.\]

When $X$ has a well-defined dimension
$\dim(X)$, the codimension of $D_\mu$ in $X^{(w(\mu))}$ is $d(\mu)\dim(X)$.  As
above, $D_\mu$ may be viewed either as a locus or a map. We will write
$D_{(n),m}$ for $D_\mu$ where $\mu$ is the unique distribution of weight
$m$ with unique nonsingleton block of size $n$. Also, we will denote by
$1^m$ the unique distribution of weight $m$ supported on $\{1\}$, whose
associated polyblock diagonal if $X\spr m.$ itself.\par
 The following is an easy remark.
\begin{lem}
For a b-partition (I.) with corresponding distribution $\mu$, the degree of
the map $OD\subpar{I.}\to D_\mu$  is \[ a(\mu):=\prod\limits_n\mu(n)!\]
\end{lem}
\subsubsection{Union Operation}
Now, we will need some operations on b-partitions and associated
distributions. Let $ u_{r,s}(I.) $ be the b-partition obtained from $(I.)$ by deleting the $r$th and $s$th blocks, $r\neq s$
and inserting their union.
We let $u_{a,b}(\mu)$ be the
corresponding operation on distributions, which
corresponds to
deleting blocks of size $a,b$ and inserting a block of size $a+b$; by
definition, $u_{a,b}(I.)=\emptyset$ unless $I.$ contains blocks of sizes
$a,b$ (two blocks of size $a$, if $a=b$); in other words,
\eqspl{union-op-distributions}{
u_{a,b}(\mu)=\begin{cases} \mu-1_a-1_b+1_{a+b}, \mu\geq 1_a+1_b\\
\emptyset, {\mathrm{otherwise}} \end{cases}} Here $1_a$ is the indicator
(characteristic) function of $a$.\par
In the geometric setting, $u_{a,b}$ correspond to intersecting with a suitable diagonal, i.e.
\[D_{u_{r,s}(I.)}=D_{I.}\cap D_{i,j}\]
where $i\in I_r, j\in I_s$ are arbitrary and $D_{i,j}$ is the
pullback of the diagonal from the $i,j$ factors.

\newsubsection{Diagonal operators on tensors}\label{diagonal operators on tensors}
Given a topological space $X$, the polydiagonals of its
symmetric products are reflected algebraically in the
(co)homology of these symmetric products. The algebraic structures that result can be defined purely algebraically,
which is the purpose of this section.
\subsubsection{Tensymmetric algebra}
Let $R$ be a commutative unitary graded $\Q$-algebra.
An example to keep in mid thoughout is a subalgebra
of the cohomology of a topological space (e.g.
manifold).
Consider the so-called 'tensymmetric' algebra \eqspl{}{
\TS(R)=\bigotimes\limits_{n=\infty}^1 \sym^{\bullet}(R).
 } Here and elsewhere, unless otherwise specified $\sym$ means $\sym_\Q$.
 Let  $\alpha.$ be a simple (decomposable) element
  in this algebra.
Then $\alpha$ can be written as

\[\alpha.=\bigotimes\limits_n\prod\limits_{i=1}^{\mu(n)}\alpha_{n,i}.\]
Here the second product is the formal one in $\sym^{\bullet}(R)$ (rather than the
one in $R$, which will be denoted $._R$ or $\prod _R$);
when in doubt, the product in $\sym$ will be denoted $._S$.
In the geometric situation, the products in $\sym$ and $\bigotimes$
both correspond to external cup products, and then will often be denoted
by $\star$.
The function $\mu$ is a length distribution,
i.e. a finitely-supported function from the positive integers to the
nonnegative integers and
we call $\mu$
the distribution associated to $(\alpha.)$.
  Thus $\alpha.$ is a tensor product of symmetric tensors,
  with the one in position $n$ having (tensor) degree $\mu(n)$.
  This yields a 'grading by
distribution':
\eqsp{ \ts(R)=\bigoplus\limits_\mu\ts_\mu(R),\\
\ts_\mu(R)=\bigotimes\limits_n\sym^{\mu(n)}(R).
}
We define the \emph{weight} of an element $\alpha.\ts_\mu(R)$ as that of the associated
distribution, i.e. \mbox{$w(\alpha.)=\sum n\mu(n)$.}
 Of course, in any simple $\alpha.$,
all but finitely many tensor factors (or '$n$-block factors',
we we shall call them) equal 1. Via the natural inclusion
$\sym^{\mu(n)}(R)\to R^{\otimes\mu(n)}$, $\alpha$ may be viewed
as an element of
\[\ti_{(I_1,...,I_k)}(R):=R^{\otimes k}\]
 for any b-partition $(I.)$ so that $(|I.|)=\mu$.
 Thus, we may define the 'inflated tensor algebra' as \eqspl{}{
\ti(R)=\bigoplus\limits_{(I.)}\ti_{(I.)}(R),
\ti_{(I_1,...,I_k)}(R)=R^{\otimes k}. } Then we have a natural
inclusion
\[\ts(R)\to \ti(R)\]
which takes $\ts_\mu(R)$ to
$\bigoplus\limits_{|I.|=\mu}\ti_{(I.)}(R)$. In the other
direction, there is also a natural symmetrization map
\[\ti(R)\to\ts(R),\]
which takes $\ti_{(I.)}(R)$ to $\ts_{|I.|}(R)$.

[[[[[************cut
************]]]]]]]]]]]]
\subsubsection{Norm operator}\label{norm-operator-sec}
For an element $\theta\in
R$, we denote by \[[m]_*(\theta)\in\sym^m(R)=\ts_{1^m}(R)\] the element (symmetric tensor)
\[\theta.1^{m-1},\] and more generally by \[[m]^s_*(\theta)\ts_{1^m}(R)\] the
element \[\theta^s.1^{m-s}.\] This called the $m$-th topological norm of
$\theta$.
See \S \ref{norm section} for a geometric interpretation.
\par
\subsubsection{Diagonal operators}
We now define a weight-preserving 'projection'
\eqspl{D-dag}{D^\dag_\bullet:\sym^{\bullet}(R)\to\ts(R).}
 This is a vector with components $D^\dag_\mu$, $\mu$ ranging over all distributions. Each
$D^\dag_\mu$ is defined as follows. Let
$m=w(\mu)=\sum n\mu(n)$ for a distribution $\mu$ and consider a
decomposable element
$\beta\in\Sym^m(R)$. Then $\beta$ comes from many elements of the form
\[\beta...=\bigotimes\limits_n\prod\limits_{i=1}^{\mu(n)}\prod\limits
_{k=1}^n \beta_{n,i,k}\] (the two
 internal products are in $\sym^{\bullet}$). We call such
  $\beta...$ a lift of $\beta$. Then let $D_\mu^\dag(\beta)$ as the following sum over all possible lifts $\beta...$ of the given $\beta$:
\[D_\mu^\dag(\beta)=\sum\limits_{\beta...\mapsto \beta} (\bigotimes\limits_n\prod\limits_{i=1}^{\mu(n)}(\prod\limits
_{k=1}^n{_R}\ \beta_{n,i,k}))\] where $\prod _R$ means product in $R$ and the middle product is again product in $\sym^{\bullet}$. Each $D^\dag_\mu$ is a
projection in the sense that it admits a right inverse. This
 right inverse is the natural 'forgetful'  map
\eqspl{}{D_\mu:\ts_\mu(R)\to\Sym^m(R), m=w(\mu),\\
\bigotimes\limits_n\prod\limits_{i=1}^{\mu(n)}\alpha_{n,i}
\mapsto \prod\limits_n\prod\limits_{i=1}^{\mu(n)}1^{n-1}\alpha_{n,i}
} (internal product is the product in $\Sym$)
Assembling these together, we get a map
\[D_\bullet[]=\bigoplus\limits_\mu D_\mu[]:\ts(R)\to\sym^{\bullet}(R)\]
Often $R$ will be a graded ring, which naturally induces a gradation on
$\ts(R)$, said to be by \emph{degree} (not to be confused with weight). If
$R$ has top piece $R^d$ endowed with a linear map $\int:R^d\to\Q$,
extended by zero to $R$, then $\int$ extends to $\ts(R)$ by multiplicativity,
i.e. for $\alpha$ decomposable,
\[\int (\alpha.)=\prod\limits_{n,i}\left(\int\alpha_{n,i}\right)\]
which of course depends only on the degree-$d$ component of each
$\alpha_{n,i}$ and vanishes if one of these components does.\par
\subsubsection{Duality}
 If
$R$ is $\Q$-self-dual, $D_\mu^\dag$ also admits a more
useful (weight-
preserving) 'Gysin adjoint' \[D_{\mu\dag}:\ts_\mu(R)\to \sym^{\bullet}(R),\] defined as follows. Let
$\sym^n(R)\to R$ be the multiplication map, which by
duality corresponds to a map
\[J_n: R\to \sym^n(R).\]
Then define
\[D_{\mu\dag}(\bigotimes\limits_n\prod\limits_{i=1}^{\mu(n)}\alpha_{n,i}
)=\prod\limits_n\prod\limits_{i=1}^{\mu(n)}J_n(\alpha_{n,i}).\]
Elements in the image of $D_{\mu\dag}$ are called \emph{polyblock
diagonal classes of type} $\mu$.
\subsubsection{Ordered analogue} All of the above admit an ordered analogue,
where both $\sym(R)$ and $\ts(R)$ are replaced by the tensor algebra $\bigotimes R$,
and partitions $\mu$ are replaced by b-partitions $(I.)$
The analogue of the map $D_\mu$ is the map
\eqsp{D_{(I_1,...,I_k)}:\bigotimes\limits^m R\to \bigotimes\limits^k R,
m=|I.|,\\
\alpha_1\cdots \alpha_m\mapsto \prod\limits_{j=1}^k(\prod\limits_{i\in I_j}{_R}
\ \alpha_i)
}
 that replaces
each tensor product inside a block by the corresponding $R$-product.
\subsubsection{Interpretation:}
Given a space $X$, a partition $\mu$ of weight $m$ corresponds to a polyblock diagonal subspace $d_\mu$, a cartesian product of symmetric products of the symmetric product $X\spr m.$.
These assemble together to a finite-to-one map
\[d_{\bullet m}=\coprod\limits_{w(\mu)=m}d_\mu\to X\spr m. .\]
 If $R$ represents some kind of
cohomology ring, e.g. the Chow ring on a variety, then
the Gysin map associated to  $d_\mu$ is $D_{\mu\dag}$.
In particular, $*1$ is the
class of a point,
$J_n$ is the Gysin map for the embedding of the small diagonal in a symmetric product. Each $\alpha_{n,i}$ is considered as living on a small
 diagonal $X\subset\Sym^n(X)$ and
  $\prod\limits_{i=1}^{\mu(n)}\alpha_{n,i}$
  lives on $X\spr\mu(n).\subset X\spr n\mu(n).$.
  The map $D_\mu^\dag$ is the
pullback map induced by the inclusions $d_\mu\to \sym^{\bullet}(X)$, while
$D_\mu$ is a natural right inverse for it.
\subsection{Discriminant operator}\label{dsc-operator-sec}
Our aim now to define a 'discriminant' operation on $\TS(R)$
that corresponds to intersecting with the big diagonal for $X$ smooth. As
part of our intersection calculus for Hilbert schemes, we will later derive a
formula for intersecting with the discriminant polarization of which this
operation will form the 'classical' part.
\par To this end we first define an operation $u_{n_1, n_2}$
on $\ts(R)$ (not preserving individual $\ts_\mu(R)$)
that corresponds to uniting two blocks of sizes $n_1, n_2$, similar to the
corresponding definition for polyblocks. The definition is:
\eqspl{def-u-n1-n2}{
&\forall\alpha=\bigotimes\limits_n\prod\limits_{i=1}^{\mu(n)}
\alpha_{n,i}\in \ts_\mu(R),
\forall n_1\neq n_2:\\ &u_{n_1, n_2:\mu}(\alpha.)=\\&
\sum\limits_{i=1}^{\mu(n_1)}
\sum\limits_{j=1}^{\mu(n_2)}
...\otimes\alpha_{n_1+n_1,1}...\alpha_{n_1+n_2,\mu(n_1+n_2)}
(\alpha_{n_1,i}._R\alpha_{n_2,j})\otimes
...\hat\alpha_{n_1,i}...\hat\alpha_{n_2,j}...\\
&u_{n,n:\mu}(\alpha.)=
\sum\limits_{1\leq i<j\leq \mu(n)}
...\otimes\alpha_{2n,1}...\alpha_{2n,\mu(2n)}
(\alpha_{n,i}._R\alpha_{n,j})\otimes ...\hat\alpha_{n,i}...\hat\alpha_{n,j}... }
(here $\hat{.}$ means 'omit', as usual). In
other words,  for the case $n_1\neq n_2$: omit in all possible ways one alpha
factor from each of the $n_1$ and $n_2$ block subproducts and insert
their $R$-product in the $n_1+n_2$-block subproduct. This defines a map
\eqsp{
u_{n_1, n_2: \mu}:\ts_\mu(R)\to \ts_{u_{n_1. n_2}( \mu)}(R)
} where $u_{n_1. n_2}( \mu)$ is as in \eqref{union-op-distributions}.\par
In particular, for the 'trivial' partition $\mu=(1^m)$, we get
\[u_{1,1,(1^m)}(\alpha_1...\alpha_m)=\sum\limits_{i<j}
(\alpha_i._R\alpha_j)\otimes(\alpha_1...\hat{\alpha_i}...\hat{\alpha_j}
...\alpha_m)\in R\otimes\sym^{m-2}(R).\]
\begin{example*}
For $\mu=(2, 1^2)=(2\mapsto 1, 1\mapsto 2)$, we have
\[ u_{2,1,\mu}(\alpha_{2}\alpha_{1,1}\alpha_{1,2})=(\alpha_2._R\alpha
_{1,1})\alpha_{1,2}+(\alpha_2._R\alpha_{1,2})\alpha_{1,1}\]
The unspecified product is that of $\sym$, of course.
In particular,
\[u_{2,1,\mu}(\alpha_{2}1_R^2)=2\alpha_21_R\]
\[u_{1,1,\mu}(\alpha_{2}\alpha_{1,1}\alpha_{1,2})=
\alpha_{2}(\alpha_{1,1}._R\alpha_{1,2})\qed\]
\end{example*}
In a rather more general vein, we can associate a similar operation  to  partition
$\nu\leq\mu$:
\eqspl{}{
u_{\nu:\mu}:\ts_\mu \to\ts_{\mu-\nu'}, \nu'(n):=\min(\nu(n)-1, 0),\\
\bigotimes\limits_n\prod\alpha_{n,j}\mapsto\bigotimes\limits_n(\sum\limits
_{|A_n|=\nu(n)}(\prod\limits_{j\in A_n} {_R}\alpha_{n,j})(\prod\limits_{j\notin A_n}\alpha_{n,j}))
}This corresponds
to uniting a set of blocks corresponding to $\nu_i$,
for $i=1,...,r$, to
a single block of size $|\nu_i|$, and $R$-multiplying the corresponding $\alpha$ factors.\par
Similarly, if $\nu_1,...,\nu_r$ are partitions with $\sum\nu_i\leq \mu$, we can define inductively
\eqspl{u-op-multi-partition}
{u_{\nu_1;...;\nu_r:\mu}:\ts_\mu(R)\to\ts_{\mu\spr r.}(R),\\
u_{\nu_1;...;\nu_r:\mu}=u_{\nu_1;...;\nu_{r-1}:\mu'}\circ
u_{\nu_r:\mu}, \mu'=\mu-\nu_r'.
}
There is a (simpler) ordered analogue of this, which takes the form
\eqsp{ou_{i,j}:\bigotimes\limits^kR\to\bigotimes\limits^{k-1}R, i<j\\
\alpha_1...\alpha_i...\alpha_j...\alpha_k\mapsto \alpha_1...
(\alpha_i._R\alpha_j)...\hat{\alpha_j}...\alpha_k
}
 Next, define an
operation corresponding to multiplication by a fixed
 ring element within a block of size $n$: for an element
$\omega\in R$, define \eqspl{def-u-n-omega}{
u_{n,\omega,\mu}(\alpha.)=\sum\limits_{i=1}^{\mu(n)}...
\otimes...\hat\alpha_{n,i}\cdot_S(\alpha_{n,i}._R\omega )... } In other words, replace in all possible ways
 an element in the $n$-block by its $R$-product with $\omega$ (this might be called $R$-multiplication by $\omega$, extended to $\ts(R)$ as a 'derivation inside the $n$-block').
We can similarly define for any
$\Q$-linear map $g:R\to R$, \eqspl{u-n-g}{
u_{n,g,\mu}(\alpha.)=\sum\limits_{i=1}^{\mu(n)}...
\otimes...\hat\alpha_{n,i}\cdot_S g(\alpha_{n,i})... ,}
(i.e. $g$ extended as a derivation or 'interior multiplication by $g$ in the $n$- block).
Again, there is a simple ordered analogue, given by
\eqspl{ou-i-g}{
ou_{i,g}(\alpha_1...\alpha_i...\alpha_k)=\alpha_1...g(\alpha_i)...\alpha_k.
} When there is no confusion, we will denote $ou_{i,g}$ by $ou_{I_i, g}$.

 Finally, in terms of these,
define the 'discriminant' operator on $\TS_\mu(R)$ by \eqspl{dsc-on-sym}
{\dsc_\mu:\ts_\mu(R)\to\ts(R),\\
\dsc_\mu=\sum\limits_{n_1\geq n_2}n_1n_2u_{n_1,n_2,\mu}}
In particular, set
\[\dsc\spr m.=\dsc_{(1^m)}.\]
Also set
\eqspl{}{
U_{\omega,\mu}(\alpha.)= \sum\limits_n\binom{n}{2}
u_{n,\omega,\mu}(\alpha.), \alpha.\in\ts_\mu(R)
}
These assemble together to maps
\[ \dsc, U_\omega :\ts(R)\to\ts(R)\]
and similarly, an ordered version on the inflated tensor
product $\ti(R)$:
\eqsp{O\dsc=\sum\limits_{i<j}u_{i,j}, OU_\omega=\sum\limits_iou_{i,\omega}.}
For future reference, it is important to note that we can write
\eqspl{dsc-as-sum-Dij}{ O\dsc-OU_\omega=\sum\limits_{i<j}
D_{i,j} } where $D_{i,j}$ acts  on $\ti_{(I.)}(R)$ as $u_{i,j}$
if $i,j$ are in different blocks of $(I.)$, and as
$ou_{k,\omega}$ if $i,j$ are both in the $k$-th block
$I_k$.\par

The motivation for this definition is the following elementary
result, which could be deduced from Macdonald's work
\cite{macd}. Recall first (see \S \ref{sec-discrim} below) that
if $X$ is a smooth curve over $\C$, there is a 'half
discriminant' class, which we abusively call discriminant,
$\dsc\spr m.$ on the symmetric product $X\spr m.$, which is
half the class of the big diagonal $D\spr m.$ (locus of
nonreduced cycles), and whose pullback to the cartesian product
coincides with the big diagonal there (with multiplicity 1).
Similarly, there is a discriminant class on a cartesian product
of symmetric products $\prod X\spr m_i.$.
\begin{lem}\label{classical-dsc-lem}
Let $X$ be a smooth curve with canonical class $\omega, R=H^.(X)$.
Then\begin{enumerate}\item
The cup-product action of the discriminant on $H^.(\prod X\spr\mu(n).)=\ts_\mu(R)$ is given by $\dsc_\mu-U_{\omega,\mu}$.
\item
The  cup product action of the discriminant on polyblock diagonal
classes is given by
\eqspl{dsc=formal dsc}{[\dsc\spr m.]\cup
[D_{\mu\dag}(\alpha.)]= [D_{\bullet\dag}(\dsc_\mu(\alpha.)-
U_{\omega,\mu}(\alpha.))], \alpha.\in\ts_\mu(R) }

\end{enumerate}
\end{lem}


\begin{proof}
The second part is just an elaboration of the first, so
it suffices to prove \eqref{dsc=formal dsc}.
Since we are working over $\Q$,
it suffices to prove both sides are equal after pullback to
the Cartesian product where the pullback of
$\dsc\spr m.$ is the big diagonal (multiplicity 1 (!)), and
splits as a sum of the diagonals pulled back from
$X\times X$, namely $\sum\limits_{a<b}D_{a,b}$, and $\alpha.$ is replaced by a class
$\alpha_{(I.)}$ on an ordered polyblock diagonal. Then clearly those $a,b$ in
different blocks of sizes $n_1, n_2$ (the sizes may be different or not) give rise to $u_{n_1,
n_2}$, while those in the same block of size $n$ give rise to $u_{n,\omega}$.
\end{proof}

 This result remains true, in fact, when $X$ is nodal (as follows, e.g. from the discussion
 in \S\ref{ordered-diags}, or by an elementary dimension-counting argument).
 However, it is of little interest in that case because of the lack of geometric
 meaning of the symmetric products. On the other hand,
one of the main ingredients of our intersection calculus, to be developed
starting in the next section, is an analogue of the Lemma for Hilbert
schemes of families of nodal curves (see
 Proposition \ref{disc.diag}), where the two
sides of \eqref{dsc=formal dsc} are not equal but differ by an 'exceptional'
class called a node scroll class. The device of pulling
back to an ordered version will be used there too.\par
\subsection{ (Half-) discriminant}\label{sec-discrim} Let $X/B$
be a family of smooth curves and $D\spr m.$ the big diagonal (or discriminant)
in the relative symmetric product
$X\spr m._B$,
i.e. \mbox{$D_\mu\cap X\spr m._B$} for $\mu=(2\mapsto 1, 1\mapsto m-2)=(2, 1^{m-2})$ (also written simply as $(2)$).
This is a reduced Cartier divisor, defined locally by the discriminant function which is
a polynomial in the elementary symmetric functions of a local parameter of $X/B$. The associated line bundle $\O(D\spr m.)$ is always canonically  divisible by 2 as line bundle.
Its half is denoted $h=\dsc\spr m.$. One way to see this
is to note that $D\spr m.$, which is the branch locus of
 $\varpi:X^m_B\to X\spr m._B$,
 is also the branch locus of a flat (albeit singular) double cover
\eqspl{}{\epsilon:X^{\odot m}_B\to X\spr m._B
}
where $X^{\odot m}_B=X^m_B/\frak A_m$ is the 'orientation product',  quotient of the cartesian product by the alternating group,
which generically parametrizes an $m$-
cycle together with
an orientation. Then $h$ is defined by
\eqspl{}{\epsilon_*\O_{X^{\odot m}_B}=\O_{X\spr m._B}\oplus h\inv.
}
(i.e. $h\inv$ is canonically the kernel
of the trace map $\epsilon_*\O_{X^{\odot m}_B}\to\O_{X\spr m._B}$).
Indeed $\epsilon^*h$ is precisely the (reduced) ramification
divisor of $\epsilon$, which is half of $\epsilon^*D\spr m.$. In particular, note that $\epsilon^*h$ is effective.\par
An explicit formula for the discriminant- which extends to Hilbert schemes as well- is the following (see \cite{lehn-montreal}, p.8): let $Z\subset X\spr m.\times X$ be the tautological
subscheme and $A=p_{1*}\O_Z$. The analogous object
on the Hilbert scheme is what we call the tautological sheaf associated to the trivial bundle
 and denote by $\Lambda_m(\O_X)$ (see \S\ref{transfer-chern-sec}). Then $A$ is endowed with
a trace pairing, whence a map $A\to A^*$ which drops
rank precisely on $D\spr m.$, hence $[D\spr m.]=-2c_1(A)$.
Therefore a half-discriminant can be defined by
\eqspl{disc-via-tauto-eq}{\dsc\spr m.=-c_1(A)=-c_1(\Lambda_m(\O_X)).}
The same formula applies to define the discriminant $\Gamma\spr m.$ on the Hilbert scheme $X\sbr m._B$ (at least for any family $X/B$ of nodal curves).\par
Though not important for our purposes it amusing to note that
the two definitions of discriminant
agree. This follows from the fact that on the orientation product $X^{\odot m}_B$, the pullback
$\epsilon^*(\det(\Lambda_m(\O_X)))$ is an ideal sheaf (hence, it is the ideal sheaf of $\epsilon^*(D\spr m.)_{\red}$):
the map $\epsilon^*(\det(\Lambda_m(\O_X)))\to \O_{X^{\odot m}}$
is given by first mapping over $X^m_B$:
\[f_1\wedge... \wedge f_m\mapsto \sum\limits_{\sigma\in\frak S_m} \sgn(\sigma)\sigma^*f_1\cdots \sigma^*f_m\]
then noting that this is $\frak A_m$-invariant on $X^m_B$, hence descends to
$X^{\odot m}_B$.

\subsection{Norm}\label{norm section}
For a line bundle $L$ on a family of smooth curves $X/B$, we denote by $[m]_*(L)$
its norm on the symmetric product $X\spr m.$, defined by
\eqspl{norm-eq}{[m]_*(L)=c_1(p_{1*}(p_2^*L\otimes\O_Z))+\dsc\spr m.=c_1(\Lambda_m(L))+\dsc\spr m.}
(notations as above).
For an effective Cartier divisor $D$ on $X$, the norm of $\O(D)$ is
\eqspl{norm-effective-eq}{ [m]_*(D)=p_{1*}(p_2^*(D).Z)}
(direct image as cycle).
To see this fact (just the Riemann-Roch for the finite
map $Z\to X\sbr m._B$), use the exact
sequence
\[\exseq{\Lambda_m(\O(-D))}{\O_X}{\O_D}.\]
Thus, $[m]_*(D)$ is a divisor supported on the locus of cycles
meeting $D$. Note again that the same formula \eqref{norm-eq}
defines the norm in the Hilbert scheme setting.\par
 In terms of cohomology, the class $[m]_*(D)$ for
 $D$ effective is
just the class corresponding to $[D]1^{m-1}$ under the identification of
$H^.(\Sym^m(X))$ with $\Sym^m(H^.(X))$. \par Similarly, we set, for $s\leq
m$,
\eqspl{norm's}{[m]^s_*(D)=\varpi_*(p_1^*(D)...p_s^*(D)) .}
This corresponds to $[D]^s1^{m-s}$.\par These constructions are compatible
with that of \S\ref{norm-operator-sec}, in the sense that
the cohomology class of $[m]^s_*(D)$ is $[m]^s_*([D])$. This is clear from the above description.
\begin{rem}
Another formula for $\dsc\spr m.$ (see \cite{R2}, Corollary 2.1) is
\eqspl{}{\dsc\spr m.=[m]_*(\omega_{X/B})\otimes \omega_{X\spr m._B/B}\inv
.}
To prove this up to numerical equivalence it suffices
 to show that the pullback both sides on the relative Cartesian
 product $X^m_B$ are isomorphic. This is proved by applying the Riemann-Hurwitz formula for smooth varieties to $\varpi:X^m_B\to X\spr m._B$.
\end{rem}
\begin{rem}\label{Theta}
In Macdonald's development \cite{macd} of intersection theory on symmetric products of a single smooth curve $X/\C$, as expounded in \cite{acgh}, Ch. 8,  is based on
the classes $[m]_*[pt]=:\theta_m$ and $[\Theta_m]$ which is the pullback of the theta-divisor $\Theta(m)$ on the Jacobian $J(m)$ via the Abel-Jacobi map (his notation is different). $\Theta(g-1)$ may be realized as the locus of effective line-bundles of degree $g-1$ in $J(g-1)$.  In fact,
$[\Theta_m]$ and $\Gamma\spr m.:=\dsc\spr m.$ are related by:
\eqspl{Gamma-Theta}{[\Theta_m]=(g+m-1)\theta_m-\Gamma\spr m..} This may be seen as follows.
For $m\geq g$, let $L$ be any line bundle of degree $m+g-1$, hence with
$m$ sections. Then the map
 \[X\spr m.\to J(g-1),\]
 \[z\mapsto L(-z)\]
is surjective and pulls back $\Theta(g-1)$ to $\Theta_m$. Therefore, $[\Theta_m]$ is the degeneracy class of the natural evaluation map
\[H^0(L)\otimes\O\to\Lambda_m(L).\]
This yields \eqref{Gamma-Theta}.

For $m<g$, $\Theta_m$ is induced by $\Theta_g$ via  $X\spr m.\to X\spr m.+(g-m)[pt]\subset X\spr g.$, and $\Gamma\spr g..X\spr m.=\Gamma\spr m.+(g-m)\theta_m$, hence
\eqref{Gamma-Theta} follows again.

In any event, passing between $\Theta$ and $\Gamma$-based theories is a
matter of simple change of variable.
\end{rem}
\newsubsection{Boundary data}\label{boundary-data-sec}
Let $\pi:X\to B$ now denote an arbitrary flat
 family of
nodal curves of arithmetic genus $g$ over an irreducible base, with smooth
generic fibre. In order to specify the additional information
required to define a node scroll, we make
the following definition.
\begin{defn}
A \underline{boundary datum} for $X/B$ consists of\begin{enumerate}
\item an irreducible variety $T$ with a map  $\delta:T\to B$
unramified
to its image, where the image is a component of the boundary,
i.e. the locus in $B$ parametrizing singular curves;
\item a 'relative node' over $T$, i.e.
 a lifting $\theta:T\to X$ of $\delta$ such that each
$\theta(t)$ is a node of $X_{\delta(t)}$;
\item a labelling, continuous in $t$, of the two branches
of $X_{\delta(t)}$ along $\theta(t)$ as $x$-axis and $y$-axis.

\end{enumerate}
Given such a datum, the \underline{associated boundary family
} $X^\theta_T$ is the
normalization (= blowup) of the base-changed family $X\times_BT$
along the section $\theta$, i.e.
\[ X^\theta_T=\Bl_\theta(X\times_BT),
\] viewed as a family of curves of genus
$g-1$ with two, everywhere distinct, individually defined
 marked points $\theta_x, \theta_y$. We denote by $\phi$ the
 natural map fitting in the diagram  $$
 \begin{matrix}
 X^\theta_T&&\\
 \downarrow&\stackrel{\phi}{\searrow}&\\
 X\times_BT&\to&X\\
 \downarrow&&\downarrow\\
 T&\stackrel{\delta}{\to}&B.
 \end{matrix}
 $$

\end{defn}
Note that a boundary datum indeed lives over the boundary of $B$;
in the other direction, we can associate to any
component $T_0$ of the boundary of $B$ a finite number boundary data in this
sense: first consider a component $T_1$ of the
normalization of $T_0\times_B\mathrm{sing}(X/B)$, which already
admits a node-valued lifting $\theta_1$ to $X$, then further base-change
by the normal cone of $\theta_1(T_1)$ in $X$
(which is 2:1 unramified, possibly disconnected, over $T_1$),
to obtain a boundary datum as above. 'Typically', the curve corresponding to a general
point in $T_0$ will have a single node $\theta$ and then the degree of $\delta$ will be 1 or 2 depending on whether the branches along $\theta$ are distinguishable in $X$ or not (they always are distinguishable if $\theta$ is a separating node
and the separated subcurves have different genera). Proceeding
in this way and taking all components which arise, we obtain finitely many boundary data which 'cover', in an obvious sense,
the entire boundary of $B$. Such a collection, weighted so that
each boundary component $T_0$ has total weight $=1$ is called a
\emph{covering system of boundary data}.

 \newsection{The tautological module}
This section will provide a recursive procedure to compute
arbitrary powers
 of the discriminant polarization $\Gamma\spr m.$ on the
 Hilbert scheme $X\sbr m._B$ (see \S \ref{sec-discrim}, especially \eqref{disc-via-tauto-eq}).
 The computation will be
 a by-product of a stronger result determining the
 (additive)
 \emph{tautological module} on $X\sbr m._B$,
 to be described informally in this introduction,
and defined formally in the body of the chapter
(see Definition \ref{taut-mod-def}).
 \par
 The tautological module, with its associated cycle map
 $$T^m=T^m_R(X/B)\to A^\bullet(X\sbr m._B)_\Q$$
 is to be defined as the $A^\bullet_\Q(B)$-vector space generated by certain
basic formal symbols called
\emph{tautological classes} (as described below).
 On the other hand, let
 $$\Q[\Gamma\spr m.]\to A^\bullet(X\sbr m._B)_\Q$$
 be the polynomial ring
  generated by the discriminant polarization.
 Then the main result of this chapter is
 \begin{thm}[Module Theorem]\label{taut-module}
Compatibly with intersection product, $T^m$ is a
 $\Q[\Gamma\spr m.]$-module; moreover, the multiplication by
 $\Gamma\spr m.$ can be described explicitly.
 \end{thm}
 Because $1\in T^m$ by definition, this statement includes
 the nonobvious assertion that
 $$\Q[\Gamma\spr m.]\subset T^m;$$
 in other words, any polynomial in $\Gamma\spr m.$ is (explicitly) tautological.
 In this sense, the Theorem includes an 'explicit' (in the
 recursive sense, at least) computation of all the powers
 of $\Gamma\spr m.$.\par
 Now the aforementioned basic tautological classes come in two
 main flavors (plus some subflavors).\begin{enumerate}
 \item The (classes of)
 \emph{(relative) diagonal loci} $\Gamma\spr m.\subp{n_1,n_2,...}$: this
 locus  is
the closure of the set of schemes
 of the form $n_1p_1+n_2p_2+...$ where $p_1,p_2...$
 are distinct smooth points of the same (arbitrary) fibre.\par
 More generally, we will consider certain 'twists' of
 these, denoted \nl $\Gamma\subp{n_1,n_2,...}[\alpha_1,\alpha_2...]$,
 where the $\alpha.$ are 'base classes', i.e. cohomology classes on $X$.
 \item The \emph{node classes}. First, the
  \emph{node scrolls} $F_j^n(\theta)$:
 these are, essentially, $\P^1$-bundles over an
analogous diagonal locus $\Gamma\spr m-n._{(n.)}$ associated to
a boundary family $X^\theta_T$ of $X_B$,
whose general
 fibre can be naturally identified with the punctual
 Hilbert scheme component $C^{n}_j$ along the node $\theta$.\par
Additionally, there are the \emph{node sections}: these are simply
 the classes $-\Gamma\spr m..F$ where $F$ is a node scroll
 as above (the terminology comes from the fact that
 $\Gamma\spr m.$ restricts to $\O(1)$ on each fibre of a
 node scroll).\par Finally, node scrolls and node sections define
 correspondence operators, pulling back (tautological) classes from
 a Hilbert scheme $(X^\theta)\sbr m-n.$.

 \end{enumerate}
 \par
 Effectively, the task of proving Theorem \ref{taut-module}
 has two parts.\begin{enumerate}
 \item Express a
 product $\Gamma\spr m..\Gamma\subp{n.}$ in terms
 of other diagonal loci and node scrolls, see Proporsition
 \ref{disc.diag}.
 \item For
 each node $\theta$ and associated ($\theta$-normalized) boundary
 family $X_T^\theta$, determine a series of explicit line bundles
 $E^n_j(\theta), j=1,...,n$ on the relative Hilbert scheme
  $(X_T^\theta)\sbr m-n._T$
together with an identification
 $$F_j^n(\theta)\simeq \P(e_j^n(\theta)\oplus e_{j+1}^n(\theta)),$$
 such that the restriction of the discriminant
 polarization $-\Gamma\spr m.$ on $F^n_j(\theta)$ becomes the standard $\O(1)$
 polarization on the projectivized vector bundle. This is just the Node Scroll Theorem of
 \cite{structure}.
In fact, it transpires that $e_j^n(\theta)$ is just the sum of the polarization
 $\Gamma\sbr m-n.$ and a suitable base divisor, that is itself a tautological
 class in the sense of Mumford.
It then follows easily that the
 restriction of an arbitrary power $(\Gamma\spr m.)^k$ on
 $F$ can be explicitly expressed in terms of
tautological classes on Hilbert schemes of lower degree on
boundary (hence smaller-dimensional) families (which in the
stable case also have lower genus) : see Theorem
\ref{gamma-power-on-nodescroll}.
 \end{enumerate}

 \newsubsection{The small diagonal}\label{smalldiag}
We begin our study of diagonal-type loci and their
intersection product with the discriminant polarization with the smallest such
locus,
 i.e. the small diagonal.  In a sense this is
 actually the heart
 of the matter, which is hardly surprising, considering as the small
 diagonal is in the 'most special' position vis-a-vis the
 discriminant. The key result is Proposition \ref{excdiv-smalldiag-prop}
 below, which is the main ingredient in determining intersection
multiplicities.\par
 The next result is in essence a corollary to the Blowup Theorem of \cite{structure}.\par
  Let $\Gamma_{(m)}\subset
 X\sbr m._B$ be the small diagonal, which parametrizes schemes with
 1-point support, and which is the pullback of the small diagonal
  \beq D_{(m)}\simeq X\subset X\spr m._B. \eeq
  This corresponds to the distribution $\mu$ with the unique nonzero
  value $\mu(m)=1$.  The
 restriction of the cycle map yields a birational morphism
 \beq \frak c_m:\Gamma_{(m)}\to X \eeq which is an isomorphism except
 over the  nodes of $X/B$.
For the remainder of the paper, we fix a covering system of boundary data
$\{(T., \delta., \theta.)\}$ as in \cite{structure}.  and focus on its typical node
$\theta$. Thus, $\theta$ is a relative node of $X/B$, $\delta:T\to B$ is a
generically finite surjective map onto a boundary component, and
$X^\theta_T$ is the blowup of $X\times_BT$ in $\theta\times_BT$. Now define a 'local model' ldeal
\eqspl{}{ &J_m<\C[[x,y]], \\
&J_m=(x^{\binom{m}{2}}, ..., x^{\binom{m-i+1}{2}}y^{\binom{i}{2}},...,y^{\binom{m}{2}}).
} Because a formal neighborhood of $\theta$ in $X$ is locally
a pullback of a family of the form $xy=t$, there is an analogous ideal defined in a formal neighborhood of $\theta$,
and because this ideal is cosupported on $\theta$ and
independent of the choice of 'local coordinates' $x,y$, it
extends to an ideal
\[ J_m^\theta<\O_X.\]
 Then let \eqspl{J-ideal-def}{ J_m^{\theta.}=\bigcap\limits_i J_m^{\theta_i}\subset
 \O_X } be the ideal sheaf whose stalk at each fibre node
  $\theta_i$ is locally of
 type $J_m$.  Note that $J^{\theta.}_m$ is
 well-defined independent of the choice of local parameters and
 independent as well of the ordering of the branches at each
 node and invariant under permutation of the set of nodes,
  hence makes sense and is globally defined on $X$.
  \begin{prop}\label{small-diag-as-blowup-prop}
  Via $\frak c_m$, $\Gamma_{(m)}$ is
  equivalent to the blow-up of $J_m^{\theta.}.$
  If $\O_{\Gamma_{(m)}}(1)_J$ denotes the canonical blowup
 polarization, we have
 \beql{Gamma-smalldiag}{\O_{\Gamma_{(m)}}(-\Gamma\spr m.)=\omega_{X/B}^{\otimes
 \binom{m}{2}}\otimes\O_{\Gamma_{(m)}}(1)_J.}
 Furthermore, if $X$ is smooth at a node $\theta$, then $\Gamma\subp{m}$ has
 multiplicity $\min(i,m-i)$ along the corresponding divisor
 $C^m_i-\{Q^m_i, Q^m_{i+1}\}$ for $ i=1,...,m-1$.
 In particular, $\Gamma\subp{m}$ is smooth along
 $(C^m_1-Q^m_2)\cup(C^m_{m-1}-Q^m_{m-1})$.

 \end{prop}
 \begin{proof} We may work with the ordered versions
 of these objects, defined on the ordered Hilbert scheme
 $X\scl m._B$, then pass to $\frak G_m$-invariants.
We first work locally over a
 neighborhood of  a point on $\theta^m\in X^m_B$ where $\theta$ is a fibre node. As shown in \cite{structure}, \S6,
 $X\scl m._B$ is obtained from the relative Cartesian product $X^m_B$ by a suitable blowup, namely that of
 the big diagonal $OD^m$.
Because blowing up and the Hilbert scheme are both compatible with
base-change,  we may then assume $X$ is a smooth surface and $X/B$
is given by $xy=t$.
 Then the ideal of $OD^m$ is generated by $G_1,...,G_m$ and $G_1$
 has the Van der Monde form $v^m_x$, while the other $G_i$ are given
 by \cite{structure}, \S6.
  We try to restrict the ideal of $OD^m$ on the small
 diagonal $OD_{(m)}\simeq X.$ To this end, note to begin
 with the natural map
 $$\I_{OD^m}\to \omega^{\binom{m}{2}}, \omega:=\omega_{X/B}.$$
 Indeed this map is clearly defined off the singular locus of $X^m_B$, hence by reflexivity of $\I_{OD^m}$ extends everywhere, hence moreover factors through a map
 $$\I_{OD^m}.OD\subp{m}=\I_{OD^m}\otimes\O_{OD\subp m}/{\mathrm{(torsion)}}\to\omega^{\binom{m}{2}}.
 $$ To identify the image, note
 that
 \beq (x_i-x_j)|_{OD_{(m)}}=dx=x\frac{dx}{x} \eeq
 and $\eta=\frac{dx}{x}=-\frac{dy}{y}$ is a local generator of $\omega$ along $\theta$.
 Therefore
 \beq G_1|_{OD_{(m)}}=x^{\binom{m}{2}}\eta^{\binom{m}{2}}. \eeq
 From \cite{structure}, loc. cit.. we then deduce
 \beql{G-smalldiag}{ G_i|_{\Gamma_{(m)}}=x^{\binom{m-i+1}{2}}y^{\binom{i}{2}}
 \eta^{\binom{m}{2}}, i=1,...,m.}
 Since $G_1,...,G_m$ generate the ideal  $I_{OD^m}$
 along $\theta$, it follows
 that over a neighborhood of $\theta$, we have
 \beq I_{OD^m}.{OD_{(m)}}\simeq
 J_m^\theta\otimes\omega^{\binom{m}{2}}. \eeq
This being true for each node, it is also true globally.
 Consequently, passing to the $\frak S_m$-quotient, we also
 have \beq I_{D^m}.{D_{(m)}}\simeq
 J_m^\theta\otimes\omega^{\binom{m}{2}}. \eeq Then pulling back to
 $X\sbr m._B$ we get \refp{Gamma-smalldiag}..\par
 Finally, it follows from the above, plus the explicit
 description of the model $H_m$, that, along the
 'finite' part $C^m_i-Q^m_{i+1}$, $\Gamma\subp{m}$
 has equation $x^{m-i}-uy^i$ where $u$ is an affine
 coordinate on $C^m_i-Q^m_{i+1}$, from which our last
 assertion follows easily.
  \end{proof}

Let us now fix the node $\theta$ and analyze locally the blowup of
the ideal $J_m=J_m^\theta=(...,x^{\binom{m-i+1}{2}}y^{\binom{i}{2}},...)$.
\begin{lem}\label{product-ideal-lem}
\[J_m=\prod\limits_{i=1}^{m-1}(x^{m-i}, y^i)\]
\end{lem}
\begin{proof}
Consider for $i=1,...,m-1$ the cone $K_i$ in the 1st quadrant $\R^2_+$ generated by
$(m-i, 0)$ and $(0,i)$, i.e. \[K_i=(\bR^2_++(m-i,0))\cup(\bR^2_++(0,i)).\]
This cone corresponds to the ideal $J_{m,i}=(x^{m-i}, y^i)$ in the sense
that $J_{m,i}$ is generated by the monomials $x^ay^b$ with $(a,b)\in K_i$. In a similar way, the ideal $\prod\limits_i J_{m,i}$ corresponds to the cone
$\sum\limits_i K_i$. Now it is easy to see, e,g, by working with the partial sum $\sum\limits_{i=1}^n K_i$ and using induction on $n$, that the latter cone
is just equal to
\[\bigcup\limits_{i=0}^m ((\binom{m-i+1}{2},\binom{i}{2})+\bR^2_+)\]
which proves our claim.
\end{proof}
Now let $X_i$ be the blowup of $X$ in $J_{m,i}=(x^{m-i}, y^i)$, which
is the subscheme of $X\times C^m_i=X\times\P^1$ defined by
\[x^{m-i}u_i=y^iv_i\]
and contains the special points $Q^m_i=[1,0], Q^m_{i+1}=[0,1]$.
The pullback of $J_{m,i}$ on $X_i$ is an invertible ideal, generated by $x^{m-i}$ near $C^m_i\setminus Q^m_{i+1}$
and by $y^i$ near $C^m_i\setminus Q^m_{i}$.
The following is an immediate consequence of Lemma \ref{product-ideal-lem}: \begin{lem}The blowup of $J_m$ is isomorphic
to
the unique component dominating $X$ of the
fibre product (over $X$): $\prod\limits_{i=1}^{m-1}(X_i/X)
:=X_1\times_X...\times_XX_{m-1}$.
\end{lem}
As was analyzed in \cite{structure},
the special fibre of $\Gamma\subp{m}$, i.e. the blowup of $J_m$,
is a chain $C^m_1\cup...\cup C^m_{m-1}$ and the point $Q^m_i\in C^m_i$ is coupled in the cartesian product
$\prod\limits_{i=1}^{m-1}C^m_i$ with $Q^m_{j+1}\in C^m_j$ for $j<i$
and with $Q^m_j\in C^m_j$ for $j>i$. It follows that if we set
\[g_i:=x^{\binom{m-i}{2}}y^{\binom{i}{2}}.x^{m-i}=x^{\binom{m-i+1}{2}}y^{\binom{i}{2}}\]
then in $\Gamma\subp{m}$,
$J_m$ is locally generated by
by $g_i$
near $Q^m_i$ and by $g_{i+1}$ near $Q^m_{i+1}$.\par
Now note that the function $x$ has along $C^m_i$ multiplicity equal to the
length of $\C[x,y]/(x^{m-i}-y^i, x)$, i.e. $i$; similarly, $y$ has multiplicity
equal to $m-i$. Therefore, the multiplicity of the invertible ideal $J_m$ itself
along $C^m_i$ is equal to \[i(m-i)+i\binom{m-i}{2}+(m-i)\binom{i}{2}=i\binom{m-i+1}{2}+
(m-i)\binom{i}{2}=\frac{i(m-i)m}{2}.\]
Also, note that at $Q^m_{i+1}$, we have affine coordinates $u_i/v_i,
v_{i+1}/u_{i+1}$ on $C^m_i, C^m_{i+1}$ respectively. These have
respective zero-
sets $C^m_{i+1}, C^m_i$ and because
\eqspl{u-times-v}{(u_i/v_i)
(v_{i+1}/u_{i+1})=xy}
which has multiplicity $m$ along either $C^m_i$ or $C^m_{i+1}$,
it follows that $u_i/v_i$ (resp. $u_{i+1}/v_{i+1}$) has multiplicity
$m$ along $C^m_{i+1}$ (resp. $C^m_i$).\par
We summarize this discussion as follows:
\begin{prop}\label{excdiv-smalldiag-prop}
(i) The pullback ideal of  $J^\theta_m$  on $\Gamma\subp{m}$
defines a Cartier
divisor of the form
\eqspl{e-nu-def-eq}{
e^\theta_m=\sum\limits_{i=1}^{m-1}\nu_{m,i}C^m_i(\theta),\\
\nu_{m,i}:=\frac{i(m-i)m}{2}.
} Moreover $x$ and $y$ have along $C^m_i(\theta)$ multiplicity
equal to $m-i, i$ respectively.\par
(ii) Each $C^m_i$ is a $\Q$- Cartier divisor on $\Gamma\subp{m}$;
$mC^m_i$ is Cartier.
\end{prop}
\begin{cor}\label{C-m-i-intersect-lem} Notations as above, we have\par
 \[(i)\ \ C^m_i(\theta)C^m_{i+1}(\theta)=\frac{1}{m}Q^m_{i+1}(\theta);\]
\[(ii)\ \ (C^m_i(\theta)^2=-\frac{1}{m}Q^m_i(\theta)-\frac{1}{m}Q^m_{i+1}(\theta).\]
\end{cor}
\begin{proof} We will fix and suppress $\theta$.\par
(i)Locally at $Q^m_{i+1}$, $mC^m_i$ and $mC^m_{i+1}$ have respective
equations $u_{i+1}/v_{i+1}, v_i/u_i$, and these locally generate an
ideal of the form $(x^{m-i}, y^{i+1}, xy)$ (note \eqref{u-times-v}),
which has colength $m$.\par
(ii) With the above notations, the principal divisor
associated to $x$ has
the form $\sum\limits_j jC^m_j$, therefore
\[C^m_i\sum\limits_j jC^m_j=0.\]
Using similarly the divisor of $y$ yields
\[C^m_i\sum\limits_j (m-j)C^m_j=0,\]
hence numerically, \[C^m_i\sum\limits_j C^m_j\sim 0.\]
Because $C^m_i$ meets $C^m_j$ only for $|j-i|\leq 1$,  (i) yields the result.
\end{proof}

\begin{cor} \label{degree-Gamma-m-cor} With the above notations, we
 have
  \beql{disc.smalldiag}{\Gamma\spr m..\Gamma\subp{m}=
 \sum\limits_{\theta, i} \nu_{m,i}C^m_i(\theta)
 -\binom{m}{2}[\omega_{X/B}].} Moreover, if $\dim(B)=1$, we have:
  \beql{}{ e_m^2=-\sigma\nu_m, \nu_m:=\frac{m^2(m-1)(m+1)}{12},}
 where $\sigma$ is the number of nodes of $X/B$;
 \beql{}{\int\limits_{\Gamma_{(m)}}(\Gamma\spr
 m.)^2=-\sigma\nu_m+\binom{m}{2}^2\omega_{X/B}^2.}
\end{cor}
\begin{rem} The components $C_{i}^m(\theta), i=1,...,m-1$ of $e_m$
 are $\P^1$-bundles over the normalization $B(\theta)$ of
 the boundary component corresponding to the node $\theta$. These
 are special cases of the {{node scrolls}}, encountered in the previous section,
 which will be further discussed in
 \S \ref{node-scrolls-sec} below.
See \S \ref{punctual-transfer-sec} for further
discussion of the small diagonal and its intersection
theory. \end{rem}
For the remainder of the paper, we set \beq\omega=\omega_{X/B}.\eeq We will view this
interchangeably as line bundle or divisor class.\par

\newsubsection{Monoblock and polyblock  digaonals: ordered case}
\label{ordered-diags} Returning to our family $X/B$ of nodal curves, we
now begin extending the results of \S\ref{smalldiag} to the more general
diagonal  loci as defined above, first for those that live over all of $B$, and
subsequently for loci associated to the boundary.
In this section, We will work with ordered objects, chiefly as a tool
for understanding their unordered analogues, to be
considered in the next section.\par
We will work here with the ordered relative Hilbert scheme
of the nodal family $X/B$, defined as
\[ X\scl m._B=X\sbr m._B\times_{X\spr m._B}X^m_B\]
where $X\sbr m._B\to{X\spr m._B}$ is of course the cycle
map $c_m$ studied at length in \cite{structure}.
As discussed in \ref{sec-discrim}, $X\sbr m._B$ is endowed
with the (half) discriminant $\Gamma\spr m.$. We denote
by $\Gamma\scl m.$ the pullback of the latter on
$X\scl m._B$, which is effective, reduced and
Cartier and admits a splitting
as Weil divisor
\[\Gamma\scl m.=\sum\limits_{1\leq a<b\leq m} D_{a,b}\]
where the summands are pullbacks of diagonals in the 2-fold product and are not Cartier.
\par
We recall the ordered
polyblock diagonal loci $OD\subpar{I.}=OD_{(I.), X/B}$ discussed in
\S\ref{partitions}. Here we will use this notation to refer to the appropriate
loci in the \emph{relative} Cartesian product $X^m_B$. In particular, we
have  the ordered monoblock diagonal
\beql{}{OD_{I,X/B}^m=OD_I=\subset X^m_B,} and the big diagonal
\beql{}{OD^m=\sum\limits_{1\leq a<b\leq m}OD^m_{a,b}.} Similar loci
exist in the ordered Hilbert scheme: \beql{}{\Gamma_I=\Gamma_I\scl m. :=
oc\inv(OD_I)\subset X\scl m._B }  Note that $OD_I$, hence $\Gamma_I$,
are defined locally near a node by equations \beql{}{x_i-x_j=0=y_i-y_j, \
\forall i,j\in I.}
\par
Generally, for any b-partition
$$(I.)=(I_1,...,I_r)\subset [1,m],$$
we have an analogous locus ( ordered \emph{polyblock  diagonal })
\beql{}{ \Gamma_{I_1|...|I_r}=\Gamma _{I_1|...|I_r}\scl m.\subset X\scl
m._B} and note that \beql{}{
\Gamma_{I_1|...|I_r}=\Gamma_{I_1}\cap...\cap\Gamma_{I_r}} (transverse
intersection). Recall that $r$ is called the
 length of the b-partition $(I.)$ and denoted $\l(I.)$. Also \beql{}{ \Gamma\subp{I.}=oc\inv(OD\subp{I.}) } where
$OD\subp{I.}\subset X^m_B$ is the analogous polyblock diagonal.  We
may view $\Gamma\subpar{I.}$ as an operator
\[\Gamma\subpar{I.}[]:\bigotimes\limits^rR\to H_\bullet(X\scl m._B), r=\l(I.)\]
\[(\alpha.)\mapsto \Gamma\subpar{I.}\cap oc^*(\alpha.)\]
where as usual $H_\bullet$ denotes a homology theory coarser
than Chow.
Thus, the values of $\Gamma\subpar{I.}[]$ are \emph{homology} rather
than cohomology classes. However, their $\frak S_m$-symmetrized
versions will descend to the (unordered) Hilbert scheme $X\sbr m._B$,
which is typically smooth, so the distinction between homology and
cohomology will not matter.
\par

Now our first goal is to determine the
intersection action of discriminant operator on a
monoblock diagonal cycle, i.e. to determine the intersection cycle
$\Gamma\scl m..\Gamma_I$ . In this computation, a key technical question
is to determine the part of $OD_I$ and $\Gamma_I$ over the boundary of
$B$, or at least its irreducible components. Thus for each boundary datum
$(\theta, T, \delta)$, with the associated map $\phi:X^\theta_T\to X$ and cartesian power \[\phi^m:(X^\theta)^m_T\to X^m_B,\] we
need to determine $(\phi^m)^*(OD_I)$ and its inverse image in
$(X^\theta_T)\scl m.$  which we call the $(\theta, T, \delta)$ boundary of
$\Gamma_I$.  A priori, it is clear that any difference between the answers
in $\sym$ and $\Hilb$ will have to do with
 node-supported loci, i.,e. node
scrolls.\par To state the answer, we recall from \cite{structure} the ordered
node scroll $OF^I_j$, which is the portion of $\varpi\inv(F^n_j(\theta),
n=|I|$, where the $n$ points coalesced in $\theta$ lie in the $I$-indexed
coordinates. This maps to $(X^\theta_T)\scl [1,m]\setminus I.$ (i.e. a copy
of $(X^\theta_T)\scl m-n.$ indexed by $\{1,...,m\}\setminus I$,rather than $\{1,...,m-n\}$). Locally
near $\theta_x\cup\theta_y$,  $(X^\theta_T)\scl [1,m]\setminus I.$  breaks up into branches corresponding
to decompositions $\{1,...,m\}\setminus I=K_x\coprod K_y$
, where $K_x, K_y$ are the
indices of the points which lie in the $x$ or
$y$-branches, denoted $X', X"$ respectively. We denote the
corresponding branches of $OF^I_j$ by $OF_j^{I:K_x|K_y}(\theta)$.\par
For index-sets $I\subset K$, we will use the notation $K/I$
to denote the quotient set identifying $I$ to a singleton, i.e. $(K\setminus I)\coprod\{I\}$.
Correspondingly, $X^{K/I}$ will denote the subset of $X^K$
consisting of points whose components indexed by $I$ are mutually equal.\par We begin with a key technical Lemma analyzing the
boundary of the monoblock diagonal $\Gamma_I$.
 \begin{lem}\label{mono-spec}
Set-theoretically, the   $(\theta, T, \delta)$ boundary of $\Gamma_I$ is the
union of the following loci, each one itself a union of irreducible
components of the boundary:
\begin{enumerate}
\item for each index-set $K$, $[1,m]\supset K\supset I$, a locus
 $\tilde\Theta _{K/I}$, mapping
birationally to its image $\Theta_{K/I}\subset OD_I=(X')^{K/I}\times (X")^{K^c}$;
\item for each $K\subset I^c=[1,m]\setminus I$, a locus
 $\tilde\Theta _{K/I}$, mapping
birationally to its image $\Theta_{K/I}\subset OD_I=(X')^K\times (X")^{K^c/I}$;
\item for each $K$ straddling $I$ and $I^c$ (i.e. meeting both),
and each $j=1,...,|I|-1$, a component
$OF_j^{I:K-I|K^c-I}(\theta)\subset OF^I_j(\theta)$
projecting
as $\P^1$-bundle to its image in
$(X^\theta_T)\scl m-|I|.$, which lies over
 $(X')^{K-I}\times_T(X")^{K^c-I}=:(X^\theta_T)^{K-I|K^c-I}\subset
 (X^\theta_T)^{m-|I|}$.
\end{enumerate}\end{lem}
\begin{proof} The loci of type (i), (ii) are clearly there and any other component
must occur at the boundary. Hence, we may fix a node $\theta$ and work
locally over a neighborhood of $\theta$ in $X$. The main point is first  to
determine the boundary of $OD_I$ (in the symmetric product). But this is
easily determined as in the $\Theta$ decomposition of \cite{structure} \S4:
the boundary is given locally by
$$\bigcup\limits_{K\subset [1,m]} OD_I\cap \Theta_K.$$
Set $\Theta_{K/I}=OD_I\cap\Theta_K$. To describe these,
there are 3 cases
depending on $K$:\begin{enumerate}\item
if $I\subset K$, then
$$\Theta_{K/I}=(X')^{K/I}\times (X")^{K^c};$$
\item if $I\subset K^c$, then
$$\Theta_{K/I}=(X')^K\times (X")^{K^c/I};$$
\item otherwise, i.e. if $I$ meets both $K$ and $K^c$, then
 $$\Theta_{K/I}=\{y_i=0,\forall i\in K\cup I,
x_i=0, \forall i\in K^c\cup I\}$$
$$=(X')^{K-I}\times (X")^{K^c-I}\times 0^I=:X^{K-I|K^c-I}$$
.
\end{enumerate}
\par
Now is an elementary check that
the loci of type (i) and (ii) are precisely the irreducible
components of the special fibre of $OD_I$, while
 the union of the loci
$\Theta_{K/I}$ of type (iii) coincides with the intersection of $OD_I$ with the fundamental locus (=image
of exceptional locus) of the ordered cycle map
$oc_m$,
i.e. the locus of cycles containing the node with multiplicity
$>1$. Also, each $\Theta_{K/I}$ of type (iii) is of codimension 2 in $OD_I$. On the other hand, each such $\Theta_{K/I}=X^{K-I|K^c-I}$ is just a component of the inverse
image in $X^m_B$ of the locus denoted 
$X^{(a,b)}$  in \cite{structure}, \S5, where $a=|K-I|, b=|K^c-I|$, and
therefore by that Lemma, the ordered cycle map over it is a union of
$\P^1$ bundles, viz
\beql{}{oc_m\inv(X^{K-I|K^c-I})=\bigcup\limits_{j=1}^{|I|-1} OF_j^{I:K-I|K^c-I}}
where $OF_j^{I:K-I|K^c-I}$ is the pullback of $F_j^{(m-a-b:a|b)}$ over
$X^{K-I|K^c-I},$ which is a $\P^1$ bundle with fibre $C^{|I|}_j$. This
concludes the proof.
\end{proof}
\textdbend Notice that, given disjoint index-sets $K_1, K_2$
with $K_1\coprod K_2= I^c$, the
number of straddler sets $K$ such that $K-I=K_1, K^c-I=K_2$ is
precisely $2^n-2$ (i.e. the number of proper nonempty subsets
of $I$). Thus, a given $OF_j^{I:K_1|K_2}$ will lie on
this many components of $\tilde{\Theta}.$ This however
is a completely separate issue from the multiplicity of
$OF_j^{I:K_1|K_2}$ in the intersection cycle
$\Gamma\sbr m..\Gamma_I$, which has to do with the blowup structure
and will be determined below.
\par
From the foregoing
 analysis, we can easily compute the intersection of
a\nl monoblock diagonal cycle with the discriminant polarization, as
follows. We will fix a covering system of boundary data
$(T_s, \delta_s, \theta_s)$,and recall that each datum
must be weighted by $\frac{1}{\deg(\delta_s)}$ (cf. \S \ref{boundary-data-sec}.
\begin{prop}\label{cut-monom} We have an equality of divisor classes on
$\Gamma_I$:
\begin{eqnarray}\label{monom-diag}
\Gamma\scl m..\Gamma_I=&\sum\limits_{i<j\notin I}
\Gamma_{I|\{i,j\}}+|I|\sum\limits_{i\notin I} \Gamma_{I\cup \{i\}}-\binom{|I|}{2}p_{\min (I)}^*\omega\\ \nonumber
&+
\sum\limits_s\frac{1}{\deg(\delta_s)}
\sum\limits_{j=1}^{|I|-1}\nu_{|I|,j}
\delta^I_{s,j*}
OF_{j}^{I}(\theta_s),\end{eqnarray}
where ${I|\{i,j\}}$ and $I\cup \{i\}$ denote the evident diblock partition and uniblock, respectively,  the 4th term denotes
the class of the image of the node scroll on $\Gamma_I$,
$OF_{j}^{I}(\theta)=
\sum\limits_{K_1\coprod K_2=I^c}OF_j^{I:K_1|K_2}(\theta)$,
$\delta^I_{s,j}$ is the natural map of the latter to
$\Gamma_I\subset X\sbr m._B$ and the multiplicities $\nu_{n,j}$
are given by \eqref{e-nu-def-eq} ;
precisely put, the line bundle on $\Gamma_I$ given by
$\O_{\Gamma_I}(\Gamma\scl m.)\otimes p_{\min(I)}^*(\omega^{|I|})$ is represented by
an effective divisor comprising the 1st, 2nd and 4th terms
of the RHS of \refp{monom-diag}..
\end{prop}
\begin{proof}
To begin with, the asserted equality trivially holds
away from the exceptional locus of $oc_m$, where the 1st,
second and third summands come from components
$\Gamma_{i,j}$ of $\Gamma\scl m.$ having $|I\cap\{i,j\}|=0,1,2$,
respectively.\par
Next, both sides being divisors on $\Gamma_I$, it will suffice to check
equality away from codimension 2, e.g.
over a generic point of each (boundary) locus
$(X^\theta_T)^{K-I|K^c-I}$.
But there, our cycle map $oc_m$ is locally just $oc_r\times
{\mathrm{iso}}$, $r=|I|$, with
$$\Gamma\scl m.\sim \Gamma\scl r.+\sum\limits_{\{i,j\}
\not\subset I}\Gamma_{i,j}.$$ We are then reduced to the case of the small
diagonal, discussed in \S\ref{smalldiag}, especially Proposition \ref{excdiv-smalldiag-prop}.
\end{proof}
Now this result immediately implies an analogous one for the
\emph{operator}
\[\Gamma_I[]:\bigotimes\limits^{m-|I|+1}R\to A_\bullet(X\scl m._B),\]
whose arguments, as products of (co)homology
classes, can be represented by cycles in generic position.
Recall that by convention, the first $R$ factor is associated with the $I$ block. Also, as $\Gamma\scl m.$ is Cartier, it defines an endomorphism
\[\Gamma\scl m.=.\cup[\Gamma\scl m.]:A_\bullet(X\scl m._B)\circlearrowleft.\]
The result
 can be
written compactly using the 'formal discriminant' operator' of \S \ref{dsc-operator-sec},
as follows.
\begin{cor} Notations as above,
\eqspl{}{\Gamma\scl
m..\Gamma_I[\alpha.]=\Gamma_I&[\odsc(\alpha.)
-OU_\omega(\alpha.)]\\ +
&\sum\limits_s\frac{1}{\deg(\delta_s)} \sum\limits_{j=1}^{|I|-1}\nu_{|I|,j}
\delta^I_{s,j*} OF_{j}^{I}(\theta_s)[ou_{I,g}( \alpha.)] ,
\alpha.\in\ts_I(R)} where $ou_{I,g}$ is as in \eqref{ou-i-g} with
$g(\alpha)=\theta_s^*(\alpha)$,
where $\theta_s$ is viewed as a (partial) section $T_s\to X$ .
\end{cor} We note that in the node scroll operator above, we
are
viewing $H^.(B)\subset H^.(X^{\theta_s})$ via pullback
through $X^{\theta_s}\to T_s\to B$ .\par Our next goal is to
extend the foregoing result  from the monoblock to the polyblock case- still
in the ordered setting. While the extension in question is in principle
straightforward, it is a bit complicated to describe. Again, a key issue is to
describe the boundary of a polyblock diagonal locus $OD\subp{I.}$ in
terms of the $\Theta$ decomposition of \cite{structure}, \S4.. Fix a
boundary datum $(T,\delta, \theta)$. To simplify notations, we will assume,
losing no generality, that the partition $I.$ is full, i.e. $\bigcup I_j=[1,m]$. Now consider an index-set
$K\subset[1,m]$. As before, $K$ is said to be a \emph{straddler} with
respect to a block $I_\l$ of $(I.)$, and $I_\l$ is a \emph{straddler block} for
$K$,   if $I_\l$ meets both $K$ and $K^c$. The \emph{straddler number}
$\strad\subp{I.}(K)$ of $K$ w.r.t. $(I.)$ is the number of straddler blocks
$I_\l$. The \emph{straddler portion} of $(I.)$ relative to $K$ is by definition
the union of all straddler blocks, i.e. \beql{}{s_K(I.)=\bigcup\limits_{I_\l\cap
K\neq\emptyset\neq I_\l\cap K^c}I_\l.
 }
 The $x$- (resp. $y$-)\emph{-portion} of $(I.)$
(relative to $K$, of course)
  are by definition
 the partitions
\beql{}{x_K(I.)=\{I_\l:I_\l\subset K\}, y_K(I.)=\{I_\l:I_\l
\subset K^c\}.}
Finally the \emph{multipartition data} associated to $(I.)$
w.r.t. $K$ are
\beql{}{\Phi_K(I.)=(s_K(I.):x_K(I.)|y_K(I.)).}
In reality, this is a partition broken up into
3 parts: the \emph{nodebound} part $s_K(I.)$, a single block,
plus 2 \emph{at large} parts, an $x$ part and a $y$ part.
As before, we set
\beql{}{X^{\Phi_K(I.)}=(X')^{x_K(I.)}\times (X")^{y_K(I.)}}
and equip it as before with the map to $X_s^m$
obtained by inserting the node $\theta$ at the
$s_K(I.)$ positions.Now the analogue of Lemma
\ref{mono-spec} is the following
\begin{lem}For any partition $(I.)$ and boundary datum
$(T,\delta,\theta)$,, the corresponding boundary
portion of $\Gamma\subp{I.}$ is
\beql{}{\bigcup \limits_{\strad\subp{I.}(K)=0}
\tilde\Theta_{K,(I.)}\cup\bigcup\limits_{\l}
\bigcup\limits_{I'.\coprod I".=I.\setminus I_\l}
\bigcup\limits_{j=1}^{|I_\l|-1}OF_j^{(I_\l:I'.|I".)}
(\theta)}
\end{lem}
\begin{proof}
First, one easily verifies:
\beql{}{OD\subp{I.}\cap \Phi_K=X^{\Phi_K(I.)}=:\Theta_{K,(I.)}}
so that
\beql{}{OD\subp{I.}\cap (X^\theta)^m_T=\bigcup\limits_{K\subset[1,m]}
\Theta_{K,(I.)}.}
Now, an elementary observation is in order. Clearly, the codimension of $OD\subp{I.}$ in $X^m_B$ is $\sum\limits_\l (|I_\l|-1)$, and this also equals the codimension of
$OD\subp{I.}\cap (X^\theta)^m_T$ in $(X^\theta)^m_T$.
On the other hand,
we have
\begin{eqnarray}\dim(\Theta_{K,(I.)})=m-\left(\sum\limits_{I_\l \textrm{\ nonstraddler rel} K}(|I_\l|-1)+\sum\limits_{I_\l \textrm{\ straddler rel} K}|I_\l|\right)\\ \nonumber
= m-\sum_\l(|I_\l|-1)-\strad\subp{I.}(K).
\end{eqnarray}
It follows that\begin{itemize}
\item the index-sets $K$ such that $\Theta_{K,(I.)}$ is a
component of the boundary $OD\subp{I.}\cap (X^\theta_T)^m$
are precisely the nonstraddlers;
\item those $K$ such that $\Theta_{K,(I.)}$ is of codimension
1 in the special fibre are precisely those of straddle number 1 (unistraddlers).
\end{itemize}
\par
Next, what are the preimages of these loci
 upstairs in the ordered Hilbert scheme $X\scl m._B$?
 They can be analyzed as in the monoblock case:
\begin{itemize}\item
 if $K$ is a nonstraddler, a general cycle parametrized by
$\Theta_{K,(I.)}$ is disjoint from the node, so there will be
a unique component $\tilde\Theta_{K,(I.)}\subset oc_m\inv(\Theta_{K,(I.)}$ dominating $\Theta_{K,(I.)}$;
\item
if $K$ is a unistraddler (straddle number $=1$), the dominant
components of
$oc_m\inv(\Theta_{K,(I.)})$ will be the $\P^1$-bundles
$F_j^{\Phi_K(I.)}, j=1,,,s_K(I.)-1$; note that if $I_\l$
the unique block making $K$ a straddler, then $\Phi_K(I.)=(I_\l:
x_K(I.)|y_K(I.))$; moreover as $K$ runs through all unistraddlers, $\Phi_K(I.)$ runs through the date consisting of a choice of block $I_\l$ plus a partition of the set of remaining blocks in two ('$x$- and $y$-blocks');
\item because all fibres of $oc_m$ are at most 1-dimensional,
while every component of the boundary is of codimension 1 in
$\Gamma\subp{I.}$,
no index-set $K$ with straddle number $\strad\subp{I.}(K)>1$
(i.e. multistraddler) can contribute a component to
that special fibre.
\end{itemize}
This completes the proof.
 \end{proof}
The import of the Lemma is that the analysis leading to Proposition
\ref{cut-monom} extends with no essential changes to the polyblock case,
and therefore the natural analogue of that Proposition holds. This is the
subject of the next Corollary which for convenience will be stated in
operator form. The statement is nearly identical to the monoblock case,
except that the node scrolls appearing will themselves contain a
polydiagonal conditions on the variable points on $X^\theta$. We will write
$\Gamma_{(I.), Y}$ to indicate the appropriate polydiagonal locus
associated to a given family $Y$ (e.g. $Y=X/B, X^\theta_T/T$ etc.) then
define \eqspl{}{ OF^{I_\l/I.}_j(\theta)=OF^{I_\l}_j.\Gamma_{(I.)\setminus I_\l,
X^\theta_T/T}\subset (X^\theta_T)\scl [1,m]\setminus I_\l..
 }In words, this is the pullback of the appropriate
  polyblock diagonal from the base of the (ordered)
node scroll. As above, the pullback via $X^\theta_T\to X$ gives an inclusion
 $H^.(X)\to H^.(X^\theta_T)$ so for any subring $R\subset H^.(X)$ containing
 $\omega$ the operator
 \[OF^{I_\l/I.}_j(\theta)[]:\bigotimes\limits^{m-|I_\l|+1}R\to H.(X\scl m._B)\]
 is defined.
\begin{cor}
\label{disc.diag-ordered}
 For any block partition $I.=I_1|...|I_r$ on $[1,m]$, we have an equality of
 operators
\eqspl{cut-poly-ord}{\Gamma\scl
m.\circ\Gamma\subpar{I.}[]=\Gamma\subpar{I.}[]\circ&(\odsc-OU_\omega)\\
&+\sum\limits_s\frac{1}{\deg(\delta_s)}\sum\limits_{\l=1}^r
\sum\limits_{j=1}^{|I_\l|-1} \nu_{|I_\l|,j}\delta^{I_\l}_{s,j*}
OF_j^{I_\l/I.}(\theta_s)[]\circ ou_{\l,\theta_s^*} } where $u_{\l,\theta_s^*}$ is interior
multiplication by $\theta_s^*:R\to H_\bullet(B)$ as in \eqref{ou-i-g}.
 \qed


\end{cor}
What the Corollary means is that the discriminant action on the
ordered Hilbert scheme is the 'classical' action on the Cartesian product (cf.  Lemma \ref{classical-dsc-lem}), plus
a boundary term as above. In the next section we will derive (easily) and explore
the corresponding result for the (unordered) Hilbert scheme.

\newsubsection{Monoblock and polyblock  diagonals: unordered case}\label{unordered-sec}
Here we will transport the formulae of the last section
to the
(unordered) Hilbert scheme. This is essentially straightforward, and
is generally accomplished by applying to
the appropriate ordered formulae push-forward by the
symmetrization map
\[\varpi_m:X\scl m._B\to X\sbr m._B.\]
We begin with the monoblock case.
Recall first the the monoblock (unordered) diagonal operator
$\Gamma\subp{n}[]$ which may be  defined for $n>1$ as
\[\Gamma\subpar{n}[\alpha.]=\frac{1}{(m-n)!}\varpi_{m*}(\Gamma\subp{I})[\alpha.],
\alpha.=\alpha_1\otimes(\alpha_2...\alpha_{m-n})\in R\otimes\sym^{m-n}_{\Q}(R).\]
Thus $\alpha_1$ is associated to a block of size $n$ while each
of $\alpha_2,...,\alpha_{m-n}$ is associated to a singleton
block.
Generally for a distribution $\mu$ the polyblock diagonal operators
$\Gamma_\mu[]$ can be defined similarly by
\[\Gamma_\mu[]=\frac{1}{a(\mu)}\varpi_*\Gamma\subpar{I.}[]:
\bigotimes\limits_n\sym^{\mu(n)}(R)
=\ts_\mu(R)\to A_\bullet(X\sbr m._B)\]
where $(I.)$ is any b-partition with distribution $\mu$ and
$a(\mu)=\prod\limits_n\mu(n)!$ is the degree of the restricted
symmetrization map $\Gamma\subpar{I.}\to\Gamma_\mu$
and $\sym^{\mu(n)}(R)$ is viewed as subring of the tensor product $\bigotimes\limits^{\mu(n)}R$
. We will often
specify a distribution by specifying only its non-singleton blocks. Thus
$\Gamma_{(n),m}$ or $\Gamma{(n)}$ for $n>1$ is short for
$\Gamma_\mu$ with $\mu$ of weight $m$ with $\mu(n)=1, \mu(1)=m-n$;
similarly for $\Gamma_{(n|n'...)}$. Note
 that in the case of the trivial partition $(1^m)$,
 the corresponding operator
 \[\Gamma_{(1^m)}[]:\sym^m(R)\to A_\bullet (X\sbr m.)\]
 is just pullback by the cycle map $c_m$. This map
 admits 'transpose' (trace map)
 \[c_{m*}: A_\bullet (X\sbr m.)\to \sym^m(R).\]
  Also, corresponding to the
Cartier divisor $\Gamma\spr m.$, we have the endomorphism
\[\Gamma\spr m.=.\cup[\Gamma\spr m.]: A_\bullet(X\sbr m._B)\circlearrowleft.\]
Note that
\[\Gamma\spr m.\Gamma_{(1^m)}[]=\half\Gamma_{(21^{m-2})}=\half\Gamma_{(2)}[].\]
We will use $\Gamma_\mu\in A_\bullet(X\sbr m._B)$ to denote $\Gamma_\mu[1_\mu]$ where $1_\mu\in \ts_\mu(R)$ is the uniquely determined product of $1_R$ factors.
We will also use $\Gamma_\mu[(\alpha_.)]$ to denote the operator
\[\ts_\mu(R)\to A_\bullet(X\sbr m._B),\]
\[(\beta.)\mapsto \Gamma_\mu[(\alpha.\cdot_R\beta.)]\]
where the product refers to the product in $\ts_\mu(R)$ induced
by the product in $R$.
Taken together, the various $\Gamma_\mu []$ operators can be
assembled to a single operator
\eqspl{Gamma-bullet}{
\Gamma_{\bullet m}:\ts(R)\to A_\bullet(X\sbr m._B).
} The main result of this section, Theorem \ref{Gamma-spr m-on-Gamma-bullet} expresses the action of the discriminant on $\Gamma_{\bullet m}$.
\begin{rem}
Using the Nakajima-inspired notation introduced in \cite{R2},
we can write
\eqspl{}{\Gamma_{(n)}[\alpha]&=q_n[\alpha],\\
\Gamma_\mu[\alpha.]&={\prod\limits_{n}}{_\star\ (q_n[\alpha_{n,1}]\star
...\star q_n[\alpha_{n,\mu(n)}])}
}
Here all products are the 'external' or star products
(see \ref{diagonal operators on tensors}) and
$q_n[\alpha]=\Gamma_{(n)}[\alpha]$ is the Nakajima-like 'creation'
operator (evaluated on 1).
\end{rem}
 Note the following elementary
facts:\begin{enumerate}
\item \beql{}{\varpi_{m*}(\Gamma\scl m..\Gamma_I)=
\Gamma\spr m..\varpi_{m*}\Gamma_I} (projection formula,
because $\varpi_m^*(\Gamma\spr m.)=\Gamma\scl m.$; NB
$\varpi$ is ramified over the support of $\Gamma\spr m.$, still
no factor of 2 in $\varpi_m^*(\Gamma\spr m.)$,
by our definition of  $\Gamma\spr m.$ as $1/2$
times its support);
\item \beql{}{\varpi_{m*}(\Gamma_I[\alpha])=
(m-n)!\Gamma_{(n)}[\alpha]\ \ , n=|I|>1;}
\item\begin{equation} \varpi_{m*} (\Gamma_{I|\{i,j\}})
=\begin{cases}
(m-n-2)!\Gamma_{(n|2)},\ \  n\neq 2;\\
 2(m-n-2)!\Gamma_{(2|2)},\ \  n=2,\\
(1+\delta_{2,n})(m-n-2)!\Gamma_{(n|2)}, \ \ \forall n
\end{cases}
\end{equation}
($\delta$= Kronecker delta) here $\Gamma\subpar{n|2}$ is the diagonal
locus corresponding to the distribution (of weight $m$) with blocks of sizes
$n,2$ plus singletons;
\item \beql{}{\varpi_{m*}(\Gamma_{I\coprod \{i\}})=(m-n-1)!\Gamma_{(n+1)};}
\item \beql{}{\varpi_{m*}(OF_j^{I:K-I|K^c-I}(\theta))=
    a!b!F_j^{(n:a|b)}(\theta),\ \ a=|K-I|, b=|K^c-I|=m-n-a} where
    we recall that $F_j^{(n:a|b)}(\theta)$ is the unordered analogue of the node scroll
    $F_j^{(I:K'|K")}$;  moreover the number
    of distinct subsets $K-I$ with $a=|K-I|$,
 for fixed $I$ and $a$, is $\binom{m-n}{a}$ which easily implies that the push-forward, properly
 weighted, of the total of the ordered node scrolls by symmetrization equals the
 total of the unordered node scrolls.
\end{enumerate}
Putting these together with Proposition \ref{cut-monom}, we conclude
\eqspl{pre-disc.monomial}{ \Gamma\spr m..\Gamma_{(n)}\sim& \\
\frac{1+\delta_{2,n}}{2} \Gamma_{(n|2)}
+n\Gamma_{(n+1)}-&\binom{n}{2}\Gamma_{(n)}[\omega]+
\sum\limits_s\frac{1}{\deg(\delta_s)}
\sum\limits_{a=0}^{m-n}\sum\limits_{j=1}^{n-1}\
\nu_{n,j}\delta^{n}_{s,j*}F_j^{(n:a|m-n-a)}(\theta_s). }
Here and elsewhere, $\Gamma\subp{n}[\omega]$ is
short for $\Gamma\subp{n}[\omega\otimes 1^{m-n}],
\omega\otimes 1^{m-n}\in R\otimes\sym_\Q^{m-n}(R)$,
either as cycle or operator.\par
 Set
\[\sum\limits_{a=0}^{m-n}F_j^{(n:a|m-n-a)}=F_j^{n,m}\]
(when $m$ is understood, we will denote this by $F^n_j$).
 We note that in this
sum, the first 3 terms in \eqref{pre-disc.monomial} match up exactly with
\eqref{dsc-on-sym}, where the first term corresponds to uniting two
singleton blocks and the second to uniting a singleton block with the
$n$-block. Therefore the formula may be extended to the twisted case and
written more compactly as follows
\begin{prop}\label{disc.monomial}
 For any monoblock diagonal $\Gamma_{(n)},\: n>1$, we have
 \eqspl{}{ \Gamma\spr m..\Gamma_{(n)}[]\sim
\Gamma\subpar{n}\circ(\dsc\spr  m.-U_\omega)+ \sum\limits_s\frac{1}{\deg(\delta_s)}
\sum\limits_{j=1}^{n-1}\ \nu_{n,j}\delta^{n}_{s,j*}F_j^{n,m}(\theta_s)[] \circ
u_{n,\theta_s^*, (n)}.}
where $u_{n,\theta_s^*, (n)}$ is as in \eqref{u-n-g}
and $\nu_{n,j}=\half j(n-j)n$.
\end{prop}
For simplicity of notation, we will denote $u_{n,g,(n)}$ by $u_{n,g}$
(e.g. $u_{n,\theta^*}$).\par

 When
$n=2$, $\Gamma_{(n)}$ is just $2\Gamma\spr m.\Gamma_{(1^m)}$, hence
(cf. \eqref{u-op-multi-partition})
\begin{cor}\label{disc-square}
\eqspl{}{(&\Gamma\spr m.)^2\Gamma_{(1^m)[]}=
\\ &\frac{1}{2}\Gamma_{(2\to2)}[]\circ u_{((1^2);(1^2):(1^m))}+\Gamma_{(3)}[]\circ u_{((1^3):(1^m))} -\Gamma\spr
m.[\omega]+\sum\limits_s
\frac{1}{\deg(\delta_s)}\half\delta^{2}_{s,j*}F_1^{2,m}(\theta_s).
[]\circ u_{2, \theta_s^*}\circ u_{((1^2:1^m))}}\qed\end{cor} Here and elsewhere, we denote by $(n\to k)$
or, more traditionally $(n^k)$,
 the
distribution $\mu$ with $\mu(n)=k$ and zeros elsewhere; when $n=1$ it
will be omitted. Also $(\mu_1|\mu_2|...|\mu_r)$ denotes the sum as functions
$\mu_1+\mu_2+...+\mu_r$; the corresponding diagonal locus
is the external or star product $\Gamma_{\mu_1}\star...\star\Gamma_{\mu_r}$.
\begin{cor}
We have \eqspl{}{\Gamma\spr m..\Gamma\subp{2}[]\circ u_{2,\omega}=
\Gamma_{(2\to 2)}[]\circ u_{2,\omega}+2\Gamma_{(3)}[]\circ
u_{3,\omega} -\Gamma\subp{2}[] \circ u_{2,\omega^2}}
 \eqspl{}{\Gamma\spr m..\Gamma\subp{3}=
\half
    \Gamma\subp{3|2}+3\Gamma\subp{4}- 3\Gamma\subp{3}[\omega]+
    \sum\limits_s \frac{3}{\deg(\delta_s)}\sum\limits_{a=0}^{m-3}
    \delta^{3}_{s,1*}(F_1^{3,m}(\theta_s)+
    \delta^{3}_{s,2*}F_2^{3,m}(\theta_s)). }
\end{cor}
\begin{proof} The first formula follows from the fact
that the boundary term involves $\theta^*$,
so drops out after multiplying by
$\omega$ because $\theta^*(\omega)$ is trivial.
The second formula is straight substitution. NB this formula,
and similar ones below in cycle form, imply analogous ones in operator from that we will leave to
the reader to explicate.
\end{proof}
\begin{cor}\label{cube}
We have for $m=2$: \eqspl{Gamma2'k}{(\Gamma\spr
2.)^k\Gamma_{(1^2)}&=\half\Gamma\subp{2}[(-\omega)^{k-1}]
+ \sum\limits_s
\frac{1}{\deg(\delta_s)}\half\delta^{2}_{s,1*} (\Gamma\spr
2.)^{k-2}.F^{2,2}_1(\theta_s)), k\geq 3;}  {{if}} $m=2, \dim(B)=1$, \eqspl{gamma-cube-degree-2-eq}{
\int\limits_{X\sbr 2._B}[\Gamma\spr 2.]^3&=
\half\omega^2-\half\sigma,\qquad \sigma=|\{\text{singular fibres}\}| ;}
 for $m=3$ \eqspl{gamma3-cube-eq}{
 (\Gamma\spr 3.)^3\Gamma_{(1^3)}&=-4\Gamma\subp{3}[\omega] +
\Gamma\spr 3.[\omega^2]\\
 +\sum\limits_s
\frac{1}{\deg(\delta_s)}(3(\delta^{3}_{s,1*}F^{3,3}_1(\theta_s)
&+\delta^{3}_{s,2*}F^{3:3}_2(\theta_s))+\half \delta^{2}_{s,1*}\Gamma\spr
3.(F_1^{2,3}(\theta_s)))
 }
\end{cor}
[We have used the elementary fact that $\omega.\theta_s=0$, hence
$\omega^i.F^{2:*}_1(\theta_s)=0, \forall i>0$, because this node scroll
maps to $\theta_s$, more precisely to $2[\theta_s]\subset X\spr
2._B$.] Note that the last term in the last equation
is minus half a node section over $X^{\theta_s}$, therefore
its support maps birationally to $X^\theta_s$. Despite the $1/2$ factor, the cycles in question are all integral
because $\Gamma\spr m.$ is integral and Cartier
(albeit non-effective). In particular, \eqref{gamma-cube-degree-2-eq}
 implies that $\omega^2-\sigma$ is even.\par
 To simplify notation we shall henceforth denote
$\frac{1}{\deg(\delta_s)}\sum\limits_sF^\bullet_\bullet(\theta_s)$ simply as
$F^\bullet_\bullet$.
\begin{example}{\rm\label{P1}This is presented here mainly as a check
on some of the coefficients in the formulas above.
For $X=\P^1$, \ \  $X\spr m.=\P(H^0(\O_X(m)))=\P^m$,
and the degree of $\Gamma\spr m._{(n)}$ is $n(m-n+1)$.
Indeed this degree may be computed as the degree of the degeneracy
locus of a generic map $n\O_X\to P^{n-1}_X(\O_X(m))$ where $P^k_X$ denotes
the $k$-th principal parts or jet sheaf.
It is not hard to show that $P^{n-1}_X(\O_X(m))\simeq n\O_X(m-n+1)$.\par
For example, $\Gamma\spr 3._{(2)}$ is a quartic scroll
equal to  the tangent developable of its
cuspidal edge, i.e. the twisted cubic $\Gamma\spr 3._{(3)}$.
The rulings are the lines $L_p=\{2p+q: q\in X\}$,
tangent to the $\Gamma\spr 3._{(3)}$,
each of which has class $-\half\Gamma\spr 3.[\omega]$.
Therefore by Corollary \ref{disc-square}, the self-intersection of
$\Gamma\spr 3.$  in $\P^3$ (or half the intersection
of $\Gamma\spr 3.$ with $\Gamma\spr 3._2$, as a class on
$\Gamma\spr 3._2$) is represented by $\Gamma\spr 3._{(3)}$ plus one ruling $L_p$.\par
If $m=4$ then $\Gamma\spr 4.$ is formally a cubic (half a
sextic hypersurface) in $\P^4$, whose self-intersection,
as given by Corollary \ref{disc-square}, is
half the Veronese $\Gamma\subp{2|2}$ plus the (sextic)
tangent developable $\Gamma\subp{3}$, plus one osculating
plane to the twisted quartic $\Gamma\subp{4}$,
representing $-\Gamma\spr 4.[\omega]$ .
\qed
}
\end{example}
Next we extend Proposition \ref{disc.monomial}
 to the polyblock  case, in other words work out the unordered analogue
 of Corollary \ref{disc.diag-ordered}.
Consider a distribution $\mu$ of weight $m$ and associated polyblock
diagonal loci and operators $D_\mu, \Gamma_\mu, \Gamma_\mu[]$ where,
e.g.
\[\Gamma_\mu[]:\ts_\mu(R)\to A_\bullet(X\sbr m._B)\]
$A_\bullet$ could be replaced by for any suitable homology theory $H_\bullet$
such as singular  (over $\Q$).
The group becomes a ring
whenever $X$ is smooth, hence so is $X\sbr m._B$.
Now the node scroll $F^n_j(\theta)$ (see the next section for more detail)
is a $\P^1$-bundle over $(X^\theta_T)\sbr m-n.$, whence operators, for
any distribution $\nu$ of weight
$m-n$: \eqspl{}{ F^n_{j,\nu}(\theta)[]:\ts_\nu(R)\to A_\bullet(X\sbr m._B)\\
\alpha.\mapsto [F^n_j(\theta)]\cap p_{[m-n]}^*\circ
\Gamma_{\nu,X^\theta_T}\circ\phi^*(\alpha.)
 } where $\phi:X^\theta_T\to X$ is the natural map.
 Clearly, given a distribution $\mu$ of weight $m$, the $\nu$-s corresponding
 to it via the unordered analogue of Corollary \ref{disc.diag-ordered} will have the form
 $\nu=\mu-1_n$ with $\mu(n)\geq 1$.
 \par A convenient way to represent the classes $\Gamma_\mu[\alpha]$
 and $F^n_{j,\nu}(\theta)[\alpha]$, adopted in the
 macnodal program (see \S\ref{macnodal}) is as matrices where
the first row represent the partition $\mu$
and each column has header $n$ and beneath it a
vector representation of the corresponding
class.

 Now the following result, which is the proper
 Hilbert scheme analogue of Lemma \ref{classical-dsc-lem},
 follows directly from
Corollary \ref{disc.diag-ordered} by adjusting for the degrees of the
various symmetrization maps. [NB The factor of $\frac{1}{\mu(n)}$ in
the boundary term on the RHS is due to the fact that
In the boundary term, the relevant ordered node scrolls map to
their unordered versions
with degree  $\prod\limits_{p\neq n}\mu(p)!(\mu(n)-1)!$, whereas $\Gamma_{(I.)}$ maps to  its unordered version
$\Gamma_\mu$ with
degree $\prod\limits_{p\neq n}\mu(p)!(\mu(n))!$. This introduces a factor of $1/\mu(n)$, which gets canceled
as there are $\mu(n)$ terms of this type.]
\begin{prop}\label{disc.diag}
For a distribution $\mu$ of weight $m$, we have an equality of operators
$\ts_\mu(R)\to A_\bullet(X\sbr m._B)$:
\eqspl{disc.diag-display}{ \Gamma\spr
m..\Gamma_\mu[]=\Gamma_\mu\circ(\dsc\spr m.-U_\omega)+
\sum\limits_s\sum\limits_{\mu(n)>0}
\sum\limits_{j=1}^{n-1}
\half j(n-j)nF^{n,m}_{j,\mu-1_n}(\theta_s)[]\circ u_{n,\theta_s^*, \mu}}
\end{prop}
\begin{example}
\eqspl{}{\Gamma\spr m..\Gamma\subp{2\to 2}[]=
\frac{3}{2}\Gamma\subp{2\to 3}[]\circ
u_{1,1}+2\Gamma\subp{4}[]\circ u_{2,2} +2\Gamma\subp{3|2}[]\circ
u_{2,1}-\Gamma\subp{2\to 2}[]\circ u_{2,\omega}+ 
\sum\limits_s F_{1}^{2,m}[]\circ u_{2,\theta^*_s}\qed }
\end{example}
\begin{example}

\eqspl{Gamma-m-cubed}{ (\Gamma\spr m.)^3\Gamma_{(1^m)}[]=
 \frac{3}{4}\Gamma\subp{2\to 3}+4\Gamma\subp{4}
+\frac{3}{2}\Gamma\subp{3|2}-\Gamma\subp{2\to 2}\circ u_{2,\omega}
-4\Gamma\subp{3}\circ u_{3,\omega}\\ + \half\Gamma_{(2)}\circ
u_{2,\omega^2}  + \frac{1}{4} F^{2, m}_{1,(2)}
 +3
(F_1^{3,m}+F_2^{3,m}) +\half\Gamma\spr m..F_1^{2,m}\qed}
\end{example}
Now assembling the various $\Gamma_\mu$ together, we
obtain our definitive result on
multiplying generalized (twisted) diagonal cycles by the discriminant polarization:
\begin{thm}\label{Gamma-spr m-on-Gamma-bullet}
We have an equality of operators $\ts(R)\to A_\bullet(X\sbr m._B)$:
\eqspl{Gamma-spr m-on-Gamma-bullet-display}{ \Gamma\spr
m..\Gamma_{\bullet m}[]=\Gamma_{\bullet m}\circ(\dsc\spr m.-U_\omega)+
\sum\limits_\mu
\sum\limits_s\sum\limits_{\mu(n)>0}
\sum\limits_{j=1}^{n-1}
\half j(n-j)nF^{n,m}_{j,\mu-1_n}(\theta_s)[]\circ u_{n,\theta_s^*, \mu}}
\end{thm}
Recall from Lemma \ref{classical-dsc-lem} that $\dsc\spr m.-U_\omega$ represents the action of the discriminant on
the the various diagonals put together. Therefore the Theorem can be viewed as a 'commutation relation' for this action:
the failure of commutativity is measured by the node scrolls. The nontrivial part is determining the multiplicities
with which they occur.\par
Because we want the Tautological Module $T^m$ (yet to be defined)
to include the $\Gamma_\mu[]$, it must also include multiples of  these by powers of the polarization $\Gamma\spr m.$.
Therefore by the above, $T^m$ must also include
the (twisted) node scrolls $F=F^n_{j,\nu}(\theta)[]$ and their multiples by powers of $\Gamma\spr m.$. Fortunately, it turns out that including the twisted scrolls $F$ and their
first-degree multiples
$\Gamma\spr m.F$ already leads to closure;
moreover, products of all these by arbitrary powers of $\Gamma\spr m.$
can be computed.
In essence, this is accomplished by the Node Scroll Theorem
of \cite{structure}.
The details are taken up in the next section.
\newsubsection{Polarized node scrolls}\label{node-scrolls-sec}
Before taking up the node scrolls, we mention an elementary analogue.
Suppose the family $X/B$ admits a relative Carter divisor $\rho$, which is
flat over $B$ of degree $k$. Then there is an induced 'incrementation'
map \eqspl{incrementation}{ \rho_+:X\sbr m._B\to X\sbr m+k._B} which
send an ideal $z$ to $z.\O(-\rho)$. In particular, if $X/B$ admits a section
$\theta$- necessarily supported in smooth points- we get maps
$(k\theta)_+$. The node scrolls are analogues of this construction where $k\theta$ is replaced by a subscheme
supported on a relative node of $X/B$.

 We recall from \cite{structure} that the node scroll
$F^{n,m}_j(\theta)$ (fixing $m$) are correspondences
\eqspl{}{\begin{matrix}&F^{n,m}_j(\theta)&\stackrel{p_{[m]}}{\to}&X\sbr m._B\\
p_{[m-n]}&\downarrow&&\\
&(X_T^\theta)\sbr m-n.&&\end{matrix}} where $p_{[m]}$ is generically
finite onto a component the locus of schemes having length at least $n$ at
$\theta$, while $p_{[m-n]}$ is a $\P^1$-bundle projection. Note that $F^{n,m}_j(\theta)$
defines an operator
\eqspl{}{
A_\bullet((X^\theta)\sbr m-n._T)\to A_\bullet(X\sbr m._B),\\
\beta\mapsto p_{[m]*}p_{[m-n]}^*(\beta) } We will however view
$F^{n,m}_j(\theta)$ as acting just on the \emph{tautological module}
$T^{m-n}(X^\theta_T)$, which we may assume defined by induction on Hilb
degree (more on this shortly), and as such its
image will be in $T^m(X/B)$. We will call $F^{n,m}_j(\theta)[\beta]$ for
$\beta\in T^{m-n}(X^\theta_T)$ a \emph{twisted node scroll class} (of Hilb
degree $m$ on $X/B$).  The \emph{polarized structure}
  of the node scroll $F^{n,m}_j(\theta)$, refers to its description
as projectivization of a
particular rank-2 vector bundle (in fact, a direct sum of two
explicit line bundles) on the degree-$(m-n)$ Hilbert scheme
$(X^\theta_T)\sbr m-n.$, with the property that the associated $\O(1)$
relative polarization coincides with $-p_{[m]}^*(\Gamma\spr m.)$. This was
worked out in \cite {structure} and can be described as follows. \par Fix a
boundary family $X^\theta_T$ and let $\theta_x, \theta_y$ be the sections
of $X^\theta_T$ mapping to the node $\theta$, and let
\[\psi_x=\theta_x^*(\omega_{X^\theta/T}),\]
considered as a line bundle on $T$ (and by pullback, on any space
mapping to $T$). As in \S\ref{norm section}, let $[k]_*L$ be the $k$-th norm
associated to a line bundle $L$ on $X$ (which is a divisor
 class on $X\sbr k._B$). Then set \eqspl{OD-n-j}{D^{n,m}_j(\theta)=
\binom{n-j+1}{2}\psi_x+\binom{j}{2}\psi_y-(n-j+1)[m-n]_*\theta_x-j[m-n]_*\theta_y
} (confusing divisors and line bundles on $(X^\theta)\sbr m-n._T$). The
Node Scroll theorem of \cite{structure} yields an isomorphism
\[F^{n,m}_j(\theta)\simeq\P(\O(D^{n,m}_j(\theta))
\oplus\O(D^{n,m}_{j+1}(\theta)))\]
under which
\[-p_{[m]}^*(\Gamma\spr m.)+p_{[m-n]}^*(\Gamma\spr m-n.)\leftrightarrow\O(1).\]

 To make use of this, set
\eqspl{e-operator-eq}{e^{n,m}_j(\theta)=[D^{n,m}_j(\theta)]-\Gamma\spr m-n.\in
A^1((X^\theta)\sbr m-n._T).}  Of course, $\Gamma\spr m-n.=0$ if $m-n\leq
1$. Thus, the $\Gamma\spr m-n.$ term begins to appear only for $m\geq 4$.  We will identify
this class with its pullback on $F^{n,m}_j(\theta)$. Then
the  $e^{n,m}_j(\theta)$, and polynomials in them, also define
operators on classes on $X^\theta_T$. Also, set formally
\eqspl{s-poly}{s_k(a,b)=a^{k}+a^{k-1}b+...+b^{k}
\quad('='\frac{a^{k+1}-b^{k+1}}{a-b})} Thus,
\eqspl{}{ s_k(e_j^{n,m},
e_{j+1}^{n,m})=\frac{(e_{j+1}^{n,m})^{k+1}-(e_j^{n,m})^{k+1}}
{-(n-j)\psi_x+j\psi_y+\theta_x-\theta_y}. }
 Then the Node Scroll Theorem plus the usual
relation of Chern and Segre classes yield immediately
\begin{thm} \label{gamma-power-on-nodescroll}For any twisted node scroll class
$F^{n,m}_{j}(\theta)[\beta]$, we have \eqspl{gamma-powers-on-nodescroll}{
(-\Gamma\spr m..)^\l.F^{n,m}_{j}(\theta)[\beta]=(-\Gamma\spr
m.)F^{n,m}_j(\theta)[s_{\l-1}(e^{n,m}_j,e^{n,m}_{j+1})\beta]-
F^{n,m}_j(\theta)[e^{n,m}_je^{n,m}_{j+1}s_{\l-2}(e^{n,m}_j,
e^{n,m}_{j+1})\beta]}
\end{thm}
[The first term is just the definition of Segre class; to get the second term, work inductively and use the $(\l-1)$ case and
the Grotendieck formula (i.e. the $\l=2$ case).]\ls
The class $-\Gamma\spr m..F^n_{j}(\theta)$, called a  \emph{node
section}, projects with degree 1 to  $(X^\theta)\sbr m-n._T$. Evaluating the
rest of the RHS of \ref{gamma-powers-on-nodescroll} involves, essentially,
the tautological module in lower degree \emph{and}, in case $X/B$ is a
family of stable curves, lower genus as well, albeit for a family of
\emph{pointed} curves $X^\theta_T$, with distinguished sections
$\theta_x, \theta_y$. To evaluate the terms involving these, we may note
the following elementary formulas, in which $\theta$ denotes any section
and $\psi=\pi_*(\omega|_\theta)$: \eqspl{}{(\theta)^r=(-\psi)^{r-1}\theta,
r\geq 1; \theta_x\theta_y=0;} \eqspl{k-theta't}{
([k]_*\theta)^t=\sum\limits_{s=1}^{\min(k,t)}(s^t-(s-1)^t)(-\psi)^{t-s}[k]_*^s(\theta)}
where we recall, cf. \eqref{norm's}, that $[k]_*^s(\theta)$ denotes the
symmetrization of $\theta^{\times s}$ and its pullback on the Hilbert
scheme.\nl
\begin{proof}[proof of \eqref{k-theta't}]: clearly,
\[([k]_*\theta)^t=\sum\limits_s\sum\limits_{\substack{r_1+...+r_s=t\\
 r_i\geq 1 \forall i}}\binom{t}{r_1,...,r_s}(-\psi)^{t-s}[k]^s_*(\theta).\]
To evaluate the
numerical coefficient, say $a_s$, note that
\[a_1+...+a_s=\sum\binom{t}{r_1,...,r_s}=s^t,\]
hence $a_s=s^t-(s-1)^t$.
\end{proof}
\par The pullback of \eqref{k-theta't} on a
polyblock diagonal $\Gamma_\nu$ is given by the $D_\nu^\dag$ operator
defined in \S\ref{diagonal operators on tensors}, viz. \eqspl{}{
\Gamma_\nu.([k]_*\theta)^t=
\sum\limits_{s=1}^{\min(k,t)}(s^t-(s-1)^t)(-\psi)^{t-s}\Gamma_\nu[D_\nu^\dag([k]_*^s(\theta))]}
Similarly, on the operator level,
 \eqspl{}{
\Gamma_\nu.([k]_*\theta)^t[\beta]=
\sum\limits_{s=1}^{\min(k,t)}(s^t-(s-1)^t)(-\psi)^{t-s}\Gamma_\nu[D_\nu^\dag([k]_*^s(\theta).\beta)]}
for $\beta\in \ts_\nu (H^.(X^\theta_T))$ (where the $.\beta$ means formal
symmetric multiplication). In particular, using the inductive case of
Theorem \ref{taut-module}, it follows that the bracketed expressions
appearing in Theorem \ref{gamma-power-on-nodescroll} are all
tautological classes, therefore
\begin{cor}\label{node-scroll-class-cor} Notations as above,
$(-\Gamma\spr m..)^\l.F^n_{j}(\theta)[\beta]$ is a twisted node scroll class.
\end{cor}
\begin{rem}
Given the canonical, mutually disjoint  sections
\[Q^{n,m}_i=\P(\O(e^{n,m}_i)), Q^{n,m}_{i+1}\subset F^{n,m}_i\]
we can write the node section in the form
\[-\Gamma\spr m.F^{n,m}_i=Q^{n,m}_i+\pi^*(e^{n,m}_{i+1}).\]
Consequently,
\eqspl{gamma-from-q1}{
-\Gamma\spr m.F^{n,m}_i=Q^{n,m}_1+\sum\limits_{1\leq j<i}
F^{n,m}_j[[m-n]_*(\theta_x-\theta_y)-(n-j)\psi_x+j\psi_y]
+F^{n,m}_i[e^{n,m}_i].
}
Consequently, rather than work with the $n-1$ node sections
$-\Gamma\spr m.F^{n,m}_i, i=1,...,n-1$, one could instead work with
a single canonical section like $Q^{n,m}_1$, together with
various twisted node scrolls.
\end{rem}

\begin{example}\label{nodescroll-dimB=1}
Note that when $\dim(B)=1$, so $T$ is a point, we have
$\theta_x\sim\theta_y\sim\theta_0$, a point on $X^\theta$. Therefore $D^{n,m}_j(\theta)$ is independent of $j$ up to numerical equivalence, hence $F^{n,m}_j(\theta)$ is also, for all $j$, deformation-equivalent, hence has the same intersection theory, as the trivial $\P^1$-bundle $(X^\theta)\spr m-n.\times\P^1$,
so that $\Gamma\spr m.\sim \Gamma_{(X^\theta)\spr m-n.}+
(n+1)[m-n]_*(\theta_0)-h_{\P^1}.$ Therefore
\[(\Gamma\spr m.)^k.F^{n,m}_j\sim (\Gamma_{(X^\theta)\spr m-n.}+
(n+1)[m-n]_*(\theta_0))^k-k(\Gamma_{(X^\theta)\spr m-n.}+
(n+1)[m-n]_*(\theta_0))^{k-1}.h_{\P^1}.
\] See Example \ref{disc-degree-single-example} for an evaluation of these
cycles.
\end{example}
\begin{example}\label{deg-pol-ord}
We have \beql{}{F_{1}^{(2,3)}(\theta)=
\P_{X^\theta}(\O(-2\theta_x-\theta_y)\oplus\O(-\theta_x-2\theta_y)) }
Consequently, if the boundary is finite, \beql{}{(-\Gamma\spr
3.)^2.F_{1}^{(2,3)}=-6. }
\end{example}
 Note that in the 'extreme' case $m=n$, the
$e^n_j(\theta)$ and the node scroll $F^n_j(\theta)$ live on the base itself
$T$ of the boundary datum and we have \eqspl{E-extreme}{
e^n_j(\theta)=\binom{m-j+1}{2}\psi_x+ \binom{j}{2}\psi_y:=\psi^m_j. }
\begin{example}
For $m=n=2, F=F^2_1(\theta)=\P(\psi_x\oplus\psi_y)$, we have
\eqspl{Gamma'k-on-F}{ (-\Gamma\spr
2.)^k|_F=(\psi_x^{k-1}+\psi_x^{k-2}\psi_y+ ...+\psi_y^{k-1})(-\Gamma\spr
2.)-\psi_x\psi_y (\psi_x^{k-2}+\psi_x^{k-3}\psi_y+ ...+\psi_y^{k-2}) .} In
particular, for $k=\dim(B)=\dim(F)=1+\dim(T)$, which is when the class
becomes 0-dimensional, we have for its degree \eqspl{}{ (-\Gamma\spr
2.)^k.F=\int\limits_T(\psi_x^{k-1}+\psi_x^{k-2}\psi_y+ ...+\psi_y^{k-1}). }
Note that if $B=\overline{\mathcal M}_g$ and $T=\overline{\mathcal
M}_{i,1}\times\overline{\mathcal M}_{g-i,1}$, $1\leq i\leq g/2$ (the usual
$i$-th boundary component), only one summand contributes to the latter
integral, which  reduces to
\[\int\limits_{\overline{\mathcal M}_i}\psi_x^{3i-2}
\int\limits_{\overline{\mathcal M}_{g-i}}\psi_y^{3(g-i)-2}
\]
\end{example}
 Note that \eqref{Gamma'k-on-F} and \eqref{Gamma2'k} together imply
\begin{cor}\label{polpowers,m=2}(i) The powers of the polarization on $X\sbr 2._B$ are
\eqspl{Gamma2-kpower}{
(-\Gamma\spr 2.)^k=-\Gamma[\omega^{k-1}]&+\\
\half\sum\limits_s\delta_{s*}( (\psi_x^{k-3}+\psi_x^{k-4}\psi_y+
...+\psi_y^{k-3})(-\Gamma\spr 2.)&-\psi_x\psi_y
(\psi_x^{k-4}+\psi_x^{k-5}\psi_y+ ...+\psi_y^{k-4}) |_{F^{2,2}_1(\theta_s)})} (ii) The image of the
latter class on the symmetric product $X\spr 2._B$ equals
\eqspl{Gamma2-kpower-on-sym}{ -\Gamma[\omega^{k-1}]+
\half\sum\limits_s\delta_{s*}( (\psi_x^{k-3}+\psi_x^{k-4}\psi_y+
...+\psi_y^{k-3})} (iii) The image of the latter class on $B$ equals
$-\kappa_{k-2}+
\half\sum\limits_s\delta_{s*}( (\psi_x^{k-3}+\psi_x^{k-4}\psi_y+
...+\psi_y^{k-3})$.
\end{cor}
\begin{proof}
(i) has been proved above; (ii)  follows because in the last summation in
\eqref{Gamma2-kpower},  the terms without $\Gamma\spr 2.$, i.e. the
twisted node scroll, collapses under the cycle map to $X\spr 2._B$; (iii)
follows similarly.
\end{proof}
\begin{rem}
It is interesting to compare the above boundary term with the boundary term in Mumford's formula \cite{Mu} for the Chern character of the Hodge bundle; our $\psi_x, \psi_y$ are his $K_1, K_2$, and $\psi_x\oplus\psi_y$ is the conormal bundle to $\theta$ in $X$;
so our term is essentially the Segre class of $\theta$ in $X$, while
Mumford's term is a Todd class of the same.
\end{rem}
\begin{example}\label{deg-pol}$m=3, n=2, \dim(B)=1$:
\beql{}{(-\Gamma\spr 3.)^2.F_{1}^{(2:3)}(\theta)=-6 } (see Example
\ref{deg-pol-ord}). Consequently, in view of Corollary \ref{cube}, we
conclude that if $X/B$ is a 'good pencil': i.e. smooth
 total space and base, all singular fibres 1-nodal, then,
 where $\sigma$ denotes the number of singular fibres, we
 have
 \eqspl{fourth}{ \int\limits_{X\sbr 3._B}(\Gamma\spr
3.)^4=13\omega^2-9\sigma } (recall that each $F^{(3,3)}_i, i=1,2$ is a line
with respect to the discriminant polarization $-\Gamma\spr 3.$).
\end{example}
\begin{example}\label{degree-mth-node-scroll}
In general, each $F^{m,m}(\theta)$ is a cycle of $\P^1$-bundles over the appropriate boundary component, whose total degree (with respect to $-\Gamma\spr m.$) is given by $\nu_m=\frac{m^2(m^2-1)}{12}$ (cf. Corollary \ref{degree-Gamma-m-cor}). In particular, for a good pencil, we get
\eqspl{}{
(-\Gamma\spr m.)F^{m,m}=\sigma\frac{m^2(m^2-1)}{12}
}
\end{example}
In the ensuing examples, we will explore a non-recursive
approach to some questions in tautological enumerative geometry based on graph enumeration. Some closely related arguments
were discovered independently and earlier by Cotteril \cite{cotteril2}.
\begin{example}\label{disc-degree-single-example}
Let $X/B$ be a single smooth curve of genus $g$ over a point,
and set
\[f_{m,g}=\int\limits_{X\spr m._B}(\Gamma\spr m.)^m=\frac{1}{m!}\int\limits_{X^m_B}(\Gamma\spr m.)^m\]
Because for genus 0 the discriminant is a hypersurface of degree $m-1$ in $(\P^1)\spr m.=\P^m$, we have
\[ f_{m,0}=(m-1)^m.\]
Consider the generating function
\[f_g(z)=1+\sum\limits_{m=2}^\infty {f_{m,g}}z^m.\]

In particular,
\[f_0(z)=1+\sum\limits_{m=2}^\infty \frac {(m-1)^m}{m!}z^m.\]
Let $W_b(z)$ be as in \S\ref{graph}.
Then by viewing $f_{m,g}$ as
$\frac{1}{m!}\int\limits_{X^m_B}(\sum\limits_{i<j}D_{i,j})^m$,
we will show that
\eqspl{f_g(z)}{f_g(z)=\exp(2(1-g)W_1(z)).}
In particular,
\eqspl{}{f_g(z)=(f_0(z))^{1-g}=(1+\sum\limits_{m=2}^\infty \frac{(m-1)^m}{m!} z^m)^{1-g}.}
Also,
\[W_1(z)=\half\log(f_0(z)).\]
To prove \eqref{f_g(z)}, expand $(\sum D_{i,j})^m$ multinomially, and attach to each monomial $M$ an edge-weighted graph with vertex-set the set of indices occurring in $M$ and
 with an edge for each $D_{i,j}$ occurring in $M$
(i.e. the multiplicity $m_{i,j}$ of the edge $(i,j)$ equals the exponent of $D_{i,j}$ in $M$).
 This graph is assigned a weight of $\frac{1}{\prod m_{i,j}!}$
due to the multinomial coefficient $\frac{m!}{\prod m_{i,j}!}$.  A 'connected' monomial, i.e. one with connected graph, with
$n$ vertices and $n$ edges, will contribute $2(1-g)w_{n,n}$ to the degree,
and a general monomial will contribute the product of the
contributions of its connected components. Then standard
generating function techniques yield the above formula.
\par More generally, consider the divisor class $\theta=\theta_m=[p+X\spr m-1.]=[m]_*(pt)\subset X\spr m.$ and let
\[h_m(u,g)=\int\limits_{X\spr m.}(\Gamma\spr m.+u\theta)^m.\]
Now because $\Gamma\spr m..X\spr m-1.=\Gamma\spr m-1.+\theta_{m-1}$, we have
\[h_m(u,g)=f_{m,g}+uf_{m-1,g}(u+1).\]
Therefore
\[h_m(u,g)=\sum\limits_{i=0}^{m-1} f_{m-i,g}u...(u+i-1)\]
In particular, when $u$ is an integer we get
\eqspl{}{h_m(u,g)=\int\limits_{X\spr m.}(\Gamma\spr m.+u\theta_m)^m
=\sum\limits_{i=0}^\infty i!\binom{u+i-1}{i}f_{m-i,g}
} Such formulas could persumably
 also be obtained by Macdonald's intersection theory (see Remark \ref{Theta}); however his methods don't seem adaptable to the singular case .\par
 Referring back to Example \ref{nodescroll-dimB=1}, it follows that in the good pencil case,
 the degree of the node scroll $F^{n,m}_j(\theta)$ is computed by
 \eqspl{}{
 \int\limits_{F^{n,m}_j(\theta)}(-\Gamma\spr m.)^{m-n+1}
 =(-1)^{m-n}(m-n+1)h_{m-n}(n+1, g-1).
 }
 See Example \ref{degree-nodescroll-partition-example} for a generalization.

\end{example}
\begin{example}\label{coefficient-omsq-Gamma-deg-pencil-example}
In the good pencil case, Theorem \ref{Gamma-spr
m-on-Gamma-bullet} shows that the degree of the discriminant
polarization has the form
\[\int\limits_{X\sbr m._B}(\Gamma\spr
m.)^{m+1}=f^1_{m,g}\omega^2+b_{m,g}\sigma\] with coefficients
universal rational numbers. The coefficient $f^1_{m,g}$, which
comes from the first, 'classical' summand in \eqref{Gamma-spr
m-on-Gamma-bullet-display}, can be determined via the
expression \eqref{dsc-as-sum-Dij} (or
 by working with
any particular smooth pencil with $\omega^2\neq 0$). Writing
formally
\[f^1_{m,g}=\frac{1}{m!}\int\limits_{X^m_B}(\sum\limits_{i<j}D_{i,j})^{m+1}\]
with the $D_{i,j}$ as in \eqref{dsc-as-sum-Dij} (where in
operator terms the integral signifies 'apply on $1^{\otimes
m}$), we can again relate these to graph numbers. Consider the
generating function
\[f^1_g(z)=1+\sum\limits_{m=2}^\infty f^1_{m,g}z^{m}.\]
Then as in Example \ref{disc-degree-single-example}, we have
\[W_2(z)\exp(2(1-g)W_1(z))=f^1_g(z)\]
i.e.
\eqspl{f1}{f^1_g(z)=W_2(z)(f_0(z))^{1-g}}
(again $W_2(z)$ is as in \S \ref{graph}).
This is because the monomials computing $f^1_{m,g}$ correspond
to graphs with one connected component of Betti number 2 and all others of Betti number 1.
\end{example}
\begin{example}\label{gamma-power-diag-part-example}

More generally we can compute, for arbitrary base dimension, the
 polyblock diagonal $\Lambda_{m,k}$ (modulo node classes ) portion of
 $\frac{1}{k!}(\Gamma\spr m.)^k$.
To state the result, we need the notion of \emph{pair-partition} or
\emph{ppartition}. By definition, a ppartition $\rho$ is an
unordered collection of ordered pairs $(n,b), n>0, b\geq 0$. Formally, a ppartition
is a function $\rho:\N\times \Z_{\geq 0}\to Z_{\geq 0}$ where
$\rho(n,b)$ counts the
frequency of a block pair of sizes $n, n+b$. The underlying partition
of $\rho$ is by definition the partitiion
$\mu=\rho_\dag$ defined by \[ \mu(n)=\sum\limits_b \rho(n,b).\]
Heuristically, we think of a ppartition $\rho$ as consisting of a partition $\rho_\dag$, together with a choice of nonnegative
'exponent' $b$ for each block, so that $\rho(n,b)$ counts the
size-$n$ block with exponent $b$.
The 'upper partition' $\mu'=\rho^\dag$ of $\rho$ is defined as \[ \mu'(b)=\sum\limits_{n}\rho(n,b).\]
we define the \emph{degree} of $\rho$ as
\[\deg(\rho)=\sum (n+b-1) \rho(n,b).\]
An ordinary partition is naturally viewed as a pppartition
with all exponents $b=0$.
Our motivation for this definition is as follows.
Given a family $X/B$, we can associate to $\rho$ a polyblock diagonal
class on $X\sbr m._B, m=|\rho_\dag|$:
\[\Gamma_\rho[(-\omega)]={\prod\limits_{n,b}}{_ \star \Gamma}_{(n)}[(-\omega)^b]^{\star
\rho(n,b)}\]
where $\star$ is star multiplication, where
$\Gamma_\mu[\alpha]\star\Gamma_\nu[\beta]=\Gamma_{\mu+\nu}[\alpha\beta]$
(see \S \ref{diagonal operators on tensors}).
Note that a connected monomial in the $D_{i,j}$ yields $\Gamma_{(n)}[(-\omega)^b]$ where $n+b$ is the degree of the monomial (aka number of edges)
and $n$ is the number of distinct indices (aka number of
vertices). Thus a possibly disconnected monomial yields $\Gamma_\rho[(-\omega)]$ for some ppartition
$\rho$.

Let $s=\dim(B)$.
Set formally, and analogously as in \S \ref{graph},
\[\tilde W_b=\sum\limits_n \frac{w_{n,n-1+b}}{n!}\Gamma\subp{n}[(-\omega)^b],\]
\[\Lambda=\sum\limits_{m,k}\Lambda_{m,k}=\sum\limits_{\rho} {L(\rho)}\Gamma_\rho[(-\omega)]:=\sum\limits
_{m,k}\frac{1}{k!}(\Gamma\spr m.)^k\Gamma_{(1^m)}
=\sum\limits_{m,k}\frac{1}{m!k!} (\sum\limits_{1\leq i<j\leq
m}D_{i,j})^k\] (equality modulo node classes). We will also use
the notation $\Lambda(X/B)$ etc. when the family needs
specification.
Note that $(-\omega)^{b}=\tilde W_b=0$ for $b>s+1$.
Then
\eqspl{L-Gamma-power}{\Lambda=\exp_\star(\sum\limits_{b=0}^{s+1}
\tilde W_b),
}

This formula results from the fact that
a connected monomial of $n$ vertices and $n-1+b$ edges yields
$\Gamma\subp{n}[(-\omega)^b]$.\par
This approach to deriving explicit formulas for powers of $\Gamma$,
and hence also for Chern polynomials in tautological bundles, can be carried
substantially further. We will return to this elsewhere (see \cite{explicit}).
\end{example}
\begin{example}
We extend the above example to compute the action of $(\Gamma\spr m.)^k$
on a twisted polyblock diagonal $\Gamma_{(n.)}[a.]$,
 for any partition $(n.)=(n_1\geq...\geq n_r)$
 of weight $m$ and corresponding collection
of classes $a_i\in R\subset H^.(X)$ (the $a_i$ can also be taken to be indeterminates).
Thus we will compute
\[(\Gamma\spr m.)^k.
\Gamma_{(n.)}[a.]=:\Lambda_{m,k,(n., a.)}+\sum F^{n,m}_j(\theta)[B_{m,k,(n., a.),n,j}^\theta]+\sum
Q^{n,m}_j(\theta)[C_{m,k,(n., a.),n,j}^\theta].\]
Set $p_i=n_ia_i.$
The computation is based on the weighted counts $w_{n,p.,m}$ as in
\eqref{w-n-p.-m}. For an index-set $S\subset [r]$, set
\[ n_S=\sum\limits_{i\in S} n_i, z^S=\prod\limits_{i\in S}z_i\]
where $z_i$ are indeterminates with $z_i^2=0.$
\[\tilde W_{(n., a.),k}=
\sum\limits_{S}\frac{w_{|S|, (p.), k}}{|S|!}\Gamma_{(n_S)}[(-\omega)^{k+1-|S|}]z^S.\]
Then we have, with a similar proof as above,
\eqspl{}{
\exp(\Gamma\spr m.)\Gamma_{(n.)}[a.]=
\sum\limits_{k}\frac{1}{k!}\Lambda_{m,k,(n.c.)}z_1...z_r=
{\exp}_\star(\sum\limits_{k)}\tilde W_{(n.,a.),k})
}
\end{example}
\begin{example}
Consider again a single smooth curve $X/B$ of genus $g$.
Extending  Example \ref{disc-degree-single-example}, we now determine the degree \[f_{m,\mu,g}=\int_{\Gamma_\mu}(\Gamma\spr m.)^{m-\l}\]
Here $\mu$ is a partition of weight $m$ and length
$\l=\sum\mu(n)$. Indeed in the case $g=0$ it is easy to see
that $\Gamma_\mu$, which is a subset of $(\P^1)\spr m.=\P^m$,
has normalization $\prod\limits_n\P^{\mu(n)}$ and the pullback of $\O_{\P^m}(1)$ is $\boxtimes\O_{\P^{\mu(n)}}(n)$. It follows that
\eqspl{}{
f_{m,\mu,0}=(m-1)^\l\l!\prod\frac{n^{\mu(n)}}{\mu(n)!}
} Then again if we form the generating function
\[ f_{\mu,g}(z)=1+\sum\limits_{m=2}^\infty
\frac{f_{m,\mu,g}}{m!}z^m\]
then we have
\eqspl{}{
f_{\mu,g}(z)=(f_{\mu,0}(z))^{1-g}=(1+\sum\limits_{m=2}^\infty
\frac{(m-1)^\l\l!}{m!}\prod\limits_n \frac{n^{\mu(n)}}{\mu(n)!})^{1-g}.
}
As in Example \ref{disc-degree-single-example}, we now extend this to compute
\[ h_{\mu,g}(u):=\int\limits_{\Gamma_\mu}(\Gamma\spr m.+u\theta)^{m-\l}, m:=|\mu|.\]
 Indeed note that
\[\theta.\Gamma_\mu=\sum\limits_n n\Gamma_{\mu-1_n}, \Gamma_{\mu-1_n}\subset X\spr m-n.\]
where $\mu-1_n$ is the partition obtained by eliminating a
block of size $n$ from $\mu$ if it has one, and otherwise
$\Gamma_{\mu-1_n}=0$; equivalently, $\theta$ acts as a
derivation for $\star$ product and equals $n$ on
$\Gamma_{(n)}$. This formula is easy to verify, e.g. by working
on the cartesian product and using the projection formula. Then
using that $f_{\mu,g}=h_{\mu,g}(0)$,
 it follows that
\eqspl{}{
h_{\mu,g}(u)=f_{\mu,g}+u\sum\limits_n nh_{\mu-1_n}(u+n,g).
} Because $f_{\mu,g}$ has been computed above and $h_\emptyset =1$, this recursion computes all the $h_{\mu,g}(u).$
Explicitly,
\eqspl{}{
h_{\mu,g}(u)=\sum\limits_{\nu\leq \mu}
\left(\prod\limits_{\{n|\nu(n)>0\}} n^{\nu(n)}(u+(\nu(n)-1)n)\right)
f_{\mu-\nu, g}
}
(sum over all partitions $\nu$ dominated by $\mu$, product over all block sizes $n$ occurring in $\nu$).
\end{example}

\begin{example}\label{degree-nodescroll-partition-example}
Assume $B$ is 1-dimensional and
 the boundary locus corresponding to $\theta$ is a point, so the
associated boundary family is a single
smooth  curve $X^\theta$
of genus $g-1$. Then using \eqref{gamma-powers-on-nodescroll}, we get
\eqsp{&\int\limits_{X\sbr m._B}(-\Gamma\spr m.)^{m-n+1}F^{n,m}_j(\theta)=\int\limits_{F^{n,m}_j(\theta)}
(-\Gamma\spr m.)^{m-n+1}
=\\
(-1)^{m-n}\sum\limits_{a=0}^{m-n} &\int\limits_{(X^\theta)^{m-n}}(\Gamma\spr m-n.+(n-j+1)\theta_x+j\theta_y)^a(\Gamma\spr m-n.+(n-j)\theta_x+(j+1)\theta_y)^{m-n-a}\\
=&(-1)^{m-n}(m-n+1)\int\limits_{(X^\theta)^{m-n}}(\Gamma\spr m-n.+(n+1)\theta_{m-n})^{m-n}=(-1)^{m-n}(m-n+1)h_{m-n}(n+1,g-1)
} (see Example \ref{disc-degree-single-example}).
Note that this number is independent of $j$ for $1\leq j\leq n-1$.
This formula extends easily to the case of a node scroll
$F^{n,m}_{j,\mu}$ constrained by a partition $\mu$ of
weight $m-n$ and degree $d=d(\mu)=\sum(n-1)\mu(n)$, i.e.
$F^{n,m}_j(\theta)$ restricted over $\Gamma_\mu(X^\theta)$:
we have
\eqspl{deg-nodescroll-with-partition-point-boundary}{
&\int\limits_{X\sbr m._B}(-\Gamma\spr m.)^{m-n+1-d}F^{n,m}_{j,\mu}(\theta)=
(-1)^{m-d-n}(m-d-n+1)\int\limits_{\Gamma_\mu(X^\theta)}
(\Gamma\spr m-n.+(n+1)\theta)^{m-n-d}\\
&=(-1)^{m-d-n}(m-d-n+1)h_{\mu}( n+1,g-1).
}
These results are readily combined with Proposition  \ref{disc.diag}. To state the result, define,
for a partion $\mu$,
\eqspl{N(mu)}{N(\mu)=
(-1)^{m-d(\mu)+1}(m-d(\mu))\sum\limits_{\mu(n)>0}\frac{\nu_n}{\mu(n)}
h_{\mu-1_n}(n+1, g-1)
}
where
\[\nu_n=\frac{n^2(n-1)(n+1)}{12}=\sum\frac{j(n-j)n}{2}\] and $h$ is as above.
Then we obtain,
for a good pencil with $\sigma$ singular fibres:
\eqspl{deg-gen-diag}{
&\int\limits_{\Gamma_\mu}(-\Gamma\spr m.)^{m-d(\mu)}=
\int\limits_{\Gamma_\mu[\dsc\spr m.-U_\omega]}(-\Gamma\spr m.)^{m-d(\mu)-1}+\sigma N(\mu).
}

Thus, the degree of a diagonal locus decomposes as a sum of a singular contribution plus the degree of a smaller-dimensional
(twisted) diagonal locus; of course, the latter terms itself decomposes in a like manner, etc.\qed
\end{example}
Combining the above examples, we obtain a closed-form formula for the degree of the discriminant polarization in the pencil case:
\begin{prop}\label{disc-deg-pencil-prop}
For a pencil $X/B$ with $\sigma$ singular fibres and canonical degree $\omega^2$, the degree of the $m$-th discriminant polarization is
\eqspl{disc-degree-pencil}{
\int_{X\sbr m._B}(-\Gamma\spr m.)^{m+1}=(-1)^{m+1}f^1_{g,m}\omega^2
+\sigma\sum\limits_{|\mu|=m}(-1)^{d(\mu)}L(\mu)N(\mu)
} where $f^1_{g,m}$ is given by \eqref{f1}, $L(\mu)$
is given by \eqref{L-Gamma-power} with $\mu$ viewed as a ppartition with
zero exponents, and $N(\mu)$ is given by \eqref{N(mu)}.
\end{prop}
\begin{rem}
In general, for a pencil $X/B$ (respectively,
 a single smooth curve over a point) and line bundle $L$ on $X$, the value of \[\int\limits_{X\sbr m._B}(-\Gamma\spr m.)^a[m]_*(L)^b, \forall a,b, a+b=m+1, b\geq 0,\]
 for fixed m, is a polynomial
in $\omega^2, L.\omega, L^2$, fibre degree $d$, fibre genus $g$, and the number $\sigma$ of singular points, linear in $\sigma$. For a single smooth curve, the analogous number is a polynomial in $d,g$ only.
Indeed the latter assertion is clear. For the former, the case $b=0$ has been discussed above. For $b>0$ we use induction on $m$ and
the flaglet Hilbert scheme discussed in \S \ref{flaglet-sec}. We can write
\[\int\limits_{X\sbr m._B}(-\Gamma\spr m.)^a[m]_*(L)^b
=\frac{1}{m}\int\limits_{X\sbr m,m-1._B}
p_{m-1}^*((-\Gamma\spr m.)^a[m]_*(L)^{b-1})a^*(L)+
p_{m-1}^*((-\Gamma\spr m.)^a[m]_*(L)^{b-2})a^*(L^2).\]
The first summand is just $d$ (=fibre degree) times an
analogous number for $m-1$, while the second summand,
which does not occur for $b=1$, is just
$L^2$ times an analogous summand for $m-1$ on a general fiber.
\end{rem}
Note that the remark implies that the values in question are independent of the distribution of genus or $L$-degrees
in reducible fibres. This had been pointed out by Gwoho Liu.

\subsection{Tautological module}
We are now in position to give the formal (recursive) definition of the
tautological module $T^m(X/B)$ and the proof of Theorem
\ref{taut-module}.
\begin{defn}\label{taut-mod-def}Let $X/B$ be a family of (possibly pointed) nodal curves.
Given a (co)homology theory $H_\bullet, H^\bullet$
 admitting a natural map from $A_\bullet, A^\bullet$ and a $\Q$-subalgebra $R\subset
H^\bullet(X)_\Q$ containing the canonical class $\omega$ and
the classes of all marked points, the tautological
module $T^m_R(X/B)$ is the $R$-submodule of
$\Hom(\ts(R), H_\bullet(X\sbr m._B))$
generated by\begin{enumerate}\item the twisted polyblock diagonal
classes $\Gamma_{\mu}[]$, $ w(\mu)=m$;
\item the direct images on $X\sbr m._B$ the twisted node scroll classes
    $F_{j}^n(\theta)[\beta]$ and the twisted node scroll sections
    $-\Gamma\spr m..F_{j}^n(\theta)[\beta]$ as $(T,\delta, \theta)$
    ranges over a fixed covering system of boundary data for the family
    $X/B$, $\beta\in T^{m-n}_R(X^\theta_T)$ and $2\leq n\leq m$.
\end{enumerate} For the default choice $R=\Q[\omega, p_1,...,p_k]$, where $p_1, ...,p_k$ are the markings,
$H^\bullet=A^\bullet, H_\bullet=A_\bullet$,
we denote
$T^m_R$ by $T^m$.
\end{defn}
\begin{proof}[Proof of Theorem \ref{taut-module}].
We wish to compute the product of a tautological class $c$ by
$\Gamma\spr m.$. If $c$ is a (twisted) diagonal class
$\Gamma_\mu[\alpha.]$, this is clear from Proposition \ref{disc.diag}. If
$c$ is a twisted node scroll class $F_{j}^n[\alpha.]$, it is obvious. Finally if
$c$ is a node scroll section $-\Gamma\spr m..F_{j}^n(\theta)[\alpha.]$, it is
clear from the case $\l=2$ of Theorem \ref{gamma-power-on-nodescroll}.
\end{proof}
\begin{rem}
It is perhaps advantageous to view $T^m_R(X/B)$ as a functor on
the category of $B$-schemes, associating to a map $T\to B$ the
module $T^m_R(\tilde X_T/T)$, where $\tilde X_T$ is a desingularization of
$X\times_BT.$ We will not pursue this formally though.
\end{rem}
\begin{rem}
In the important special case of computing a power $(\Gamma\spr m.)^k$
it is probably more efficient not to proceed by simple recursion, but rather
to apply just Proposition \ref{disc.diag} repeatedly to express
$(\Gamma\spr m.)^k$ in terms of twisted diagonals plus classes
$(\Gamma\spr m.)^t.F$ for various $t$'s and various $F$'s; then each of
the latter classes can be computed at once using Theorem
\ref{gamma-power-on-nodescroll}.
\end{rem}
\newsection{Tautological transfer and Chern numbers}\label{transfer}
In this chapter we will complete the development of our intersection calculus. First we study the \emph{transfer} operation $\tau_m$, taking cycles on $X\sbr m-1._B$ to cycles
(of dimension larger by 1) on $X\sbr m._B$, via the flag Hilbert
scheme $X\sbr\mm._B$.
Our main aim is to show that transfer can be effected
within the Tautological modules.
In the Transfer Theorem \ref{taut-tfr} we will show
in fact that for any basic
tautological class $u$ on $X\sbr m-1._B$,
the image $\tau_m(u)$ is a simple linear combination of
basic tautological classes on $X\sbr m._B$. We then review a
splitting principle established in \cite{R}, which expresses the Chern classes of the tautological bundle $\Lambda_m(E)$,
pulled back on $X\sbr \mm._B$, in terms of
those of $\Lambda_{m-1}(E)$, the discriminant polarization
$\Gamma\spr m.$, and base classes. Putting this result together with the Module Theorem and the Transfer Theorem yields the calculus for arbitrary polynomials in the Chern classes of $\Lambda_m(E).$
\newsubsection{Flaglet geometry and the transfer theorem}
\label{flaglet-sec}
In this section we study the $(m,m-1)$ flag (or 'flaglet') Hilbert scheme,
which we view as a correspondence between the Hilbert schemes for
lengths $m$ and $m-1$ providing a way of transporting cycles, especially
tautological ones, between these Hilbert schemes. We will make strong,
chapter-verse use of the results of \cite{Hilb}.\par Thus let
$$X\sbr m,m-1._B\subset X\sbr m.\times _B X\sbr m-1.$$
denote the flag Hilbert scheme, parametrizing
pairs of schemes $(z_1, z_2)$ satisfying
$z_1\supset z_2$. This comes equipped with a (flag) cycle
map
$$c_{m,m-1}:X\sbr m,m-1._B\to X\spr m,m-1._B,$$
where $X\spr m,m-1._B\subset X\spr m._B\times_B X\spr m-1._B$ is the
subvariety parametrizing cycle pairs $(c_{m}\geq c_{m-1}).$ Note
that the \emph{normalization} of $X\spr \mm._B$ may be identified with
$X\spr m-1._B\times _B X$; however the normalization map, though bijective,
is not an isomorphism. Note also that we also have an ordered
version $X\scl m,m-1._B$, with its own cycle map
$$oc_{m,m-1}:X\scl m,m-1._B\to X^m_B.$$
In addition to the obvious projections
\beql{}{\begin{matrix}&&X\sbr\mm._B&&\\
&p_m\swarrow&&\searrow p_{m-1}\\
X\sbr m._B&&&&X\sbr m-1._B
\end{matrix}
}
with respective generic fibres $m$ distinct points (corresponding to removing a point from a
given $m$-tuple)
and a generic fibre of $X/B$ (corresponding to adding a point
to a given $m-1$-tuple), $X\sbr \mm._B$ admits
a natural map
\beql{}{a:X\sbr\mm._B\to X,}
$$(z_1\supset z_2)\mapsto \mathrm{ann}(z_1/z_2)$$
(identifying $X$ with the Hilbert scheme of colength-1 ideals).
Therefore $X\sbr \mm._B$ admits a 'refined cycle map'
(factoring the flag cycle map) \beql{}{c:X\sbr \mm._B\to X\times_B X\spr m-1._B
}
$$c=a\times (c_{m-1}\circ p_{m-1}).$$
\par
Now in \cite{Hilb} (Theorem 5 et seq., especially
Construction 5.4 p.442) we worked out a
complete model for $X\sbr m,m-1._B$, locally over
$X\spr\mm._B$. Let
\begin{eqnarray}H_m&\subset& X\spr m._B\times\tilde C^m
_{[u.,v.]}
\subset X\spr m._B\times\P^{m-1}_{Z.},\\
H_{m-1}&\subset& X\spr m-1._B\times \tilde C^{m-1}
_{[u'.,v'.]}
\subset X\spr m-1._B\times\P^{m-2}_{Z'.}
\end{eqnarray}
be respective local models for $X\sbr m._B, X\sbr m-1._B$
as constructed in \S1 above, with coordinates
as indicated. Consider  the subscheme
\beql{}{H_{\mm}\subset H_m\times_BH_{m-1}
\times_{X\spr m._B\times X\spr m-1._B}X\spr \mm._B
}
defined by the equations
\beql{}{u_i'v_i=(\sigma^x_1-\sigma^{'x}_1)u_iv'_i,\ \
v'_iu_{i+1}=(\sigma^y_1-\sigma^{'y}_1)v_{i+1}u'_i,\ \
 1\leq i\leq m-2
}
or alternatively, in terms of the $Z$ coordinates,
\begin{eqnarray}\label{ZZ'} Z_iZ'_j&=&(\sigma^x_1-\sigma^{'x}_1)
Z_{i+1}Z'_{j-1},\ \  i+1\leq j\leq m-1\\ \nonumber
   &=& (\sigma^y_1-\sigma^{'y}_1)Z_{i-1}Z'_{j+1},l
\ \ 1\leq j\leq i-2.
\end{eqnarray}
To 'explain' these relations in part, note that
in the ordered model over $X^m_B$, we have
$$\sigma^x_1-\sigma^{'x}_1=x_m, \sigma^y_1-\sigma^{'y}_1
=y_m$$ and then the analogue of \refp{ZZ'}. for the $G$ functions is
immediate from the definition of these in \cite{structure}, \S4. Then the
result of \cite{Hilb}, Thm. 5, is that $H_{\mm}$, with its map to $X\spr
\mm._B$ is isomorphic to a neighborhood of the special fibre over $(mp,
(m-1)p)$ of the flag Hilbert scheme $X\sbr \mm._B$. In fact the result of
\cite{Hilb} is even more precise and identifies $H_{\mm}$ with a
subscheme of $H_m\times_B H_{m-1}$ and even of $H_{m-1}\times_B
\tilde C^m\times_B X$, where the map to $X$ is the annihilator map $a$
above.\par As noted in \cite{Hilb}, Thm 5, the special fibre of the flag cycle
map on $H_{\mm}$, aka the punctual flag Hilbert scheme, is a
normal-crossing chain of $\P^1$'s: \beql{C-mm-1}{C^{\mm}=\td
C^m_1\cup \td C^{m-1}_1\cup \td C^m_2\cup... \cup \td
C^{m-1}_{m-2}\cup \td C^m_{m-1}\subset
 C^m\times C^{m-1}.
}
where the embedding is via
$$\td C^m_i\to C^m_i\times \{Q^{m-1}_i\},\ \
 \td C^{m-1}_i\to \{Q^m_{i+1}\}
\times C^{m-1}_i$$
and in particular,
\beql{}{\td C^m_i\cap \td C^{m-1}_i=\{(Q^m_{i+1}, Q^{m-1}_i)\},
\td C^{m-1}_i\cap \td C^m_{i+1}=\{(Q^m_{i+1}, Q^{m-1}_{i+1})\}
}
where $Q^m_i=(x^{m-i+1}, y^i)$ as usual.
\begin{thm}\label{flag-blowup}
The cycle map $c_{\mm}$ exhibits the flag Hilbert scheme
$X\sbr\mm._B$ as the blow-up of the sheaf of ideals
$\I_{D^{\mm}}:=\I_{D^{m-1}}.\I_{D^m}$ on $X\spr\mm._B$.
\end{thm}
We shall not really need this result, just the explicit constructions above, so
we just sketch the proof, which is analogous to that of the Blowup Theorem
of \cite{structure}. To begin with, it is again sufficient to prove the ordered
analogue of this result, for the 'ordered flag cycle map'
$$X\scl\mm._B\to X^m_B.$$
Here $X\scl\mm._B$ is embedded as a subscheme of
$X\scl m._B\times_{X^m_B}(X\scl m-1._B\times _B X)$,
and we have already observed that as such, it satisfies
the equations \refp{ZZ'}..\par
Now we will use the following construction. Let $\I_1, \I_2$
be ideals on a scheme $Y$. Then the surjection of
graded algebras
$$(\bigoplus\limits_n \I_1^n)\otimes
(\bigoplus\limits_n \I_2^n)\to \bigoplus\limits_n (\I_1\I_2)^n$$ yields a
closed immersion \beql{}{\Bl_{\I_1\I_2}Y\hookrightarrow
\Bl_{\I_1}Y\times_Y \Bl_{\I_2}Y; } the latter is in turn a subscheme of the
Segre subscheme \beql{}{\P(\I_1)\times_Y
\P(\I_2)\subset\P(\I_1\otimes\I_2) .} In our case, the Blowup Theorem of
\cite{structure} allows us to identify
$$OH_m\simeq \Bl_{\I_{OD^m}}X^m_B,\ OH_{m-1}\times_B X
\simeq\Bl_{\I_{OD^{m-1}.X^m_B}}X^m_B$$ (where
$OH_m=H_m\times_{X\spr m._B}X^m_B$ etc.), whence an embedding
\beql{}{\Bl_{\I_{OD^{\mm}}}X^m_B \to OH_m\times_{X^m_B}
(OH_{m-1}\times_BX) } As observed above, the generators $G_i\cdot
G_j'$ satisfy the analogues of the relations \refp{ZZ'}., so the image is
actually contained in $OH_{\mm}$, so we have an embedding
\beql{}{\Bl_{\I_{OD^{\mm}}}X^m_B \to OH_{\mm}. } We are claiming that
this is an isomorphism. This can be verified locally, as in the proof of the
Blowup Theorem in \cite{structure}. \qed\par One consequence of the
explicit local model for $X\sbr\mm._B$ is the following
\begin{cor}\label{q-flat}
\begin{enumerate}
\item The projection $q_{m-1}$ is flat, with 1-dimensional fibres;\item
Let $z\in X\sbr m-1._B$ be a subscheme of a fibre $X_s$, and let $z_0$
be the part of $z$ supported on nodes of $X_s$, if any. Then
 if $z_0$ is principal (i.e. Cartier)  on $X_s$, the
fibre $q_{m-1}\inv(z)$ is birational to $X_s$ and its general
members are equal to $z_0$ locally at the nodes.
\end{enumerate}
\end{cor}
\begin{proof}
(i) is proven in \cite{R}, and also follows easily from our
explicit model $H_{\mm}$. As for (ii), we may suppose, in the
notation of \cite{Hilb}, that $z_0$ is of type $I^n_i(a)$. Now if $z'\in q_{m-1}\inv(z)$, then the part $z'_0$ of $z'$ supported
on nodes must have length $n$ or $n+1$. In the former case
$z'_0=z_0$, while in the latter case
$z'_0$ must equal $Q^{n+1}_{i+1}$ by
 \cite{Hilb}, Thm. 5 p. 438, in which case $z'$ is unique, hence
 not general.
\end{proof}
Another consequence of the Theorem is a description of the 'flaglet'
small diagonal. Denote by \[\tilde\Gamma\subp{ m}\subset X\sbr m,m-1.\]
the inverse image of the small diagonal $\Gamma\subp{ m}\subset X\sbr m._B$. In Proposition \ref{small-diag-as-blowup-prop} we described the
usual small diagonal as an explicit blowup. Thanks to the Theorem,
we obtain an analogous description for the flaglet small  diagonal:
\begin{cor}\label{flaglet-smalldiag-cor}
The cycle map $\tilde \Gamma\subp{ m}\to X$,  is the blowing up of the ideal
$J^{\theta.}_mJ^{\theta.}_{m-1}$ (cf. \eqref{J-ideal-def}).
\end{cor}
See \S \ref{punctual-transfer-sec} for further discussion of the flaglet
small diagonal and the associated 'punctual' transfer operation.\par
Next we define the fundamental {transfer} operation. Essentially, this takes cycles from $X\sbr m-1._B$ to $X\sbr m._B$, but we  also allow the additional flexibility of
 twisting by base classes
via the $m$-th factor. Thus the \emph{twisted transfer map}
$\tau_m$ is defined by
\beql{}{\tau_m=\tau_{m,X/B}:A.(X\sbr m-1._B)\otimes A.(X)\to A_{.}(X\sbr m._B)_\Q,
}
$$\tau_m=q_{m*}(q_{m-1}^*\otimes a^*).$$
Note that this operation raises dimension by 1 and preserves
codimension. Suggestively, and a little abusively, we will
write a typical decomposable element of the source of $\tau_m$
as $\gamma\beta\subp{m}$ where $\gamma\in A.(X\sbr m-1._B), \beta\in A.(X)$. Note also that when $B$ is a point (and $X$
is smooth), $\tau_m$ is just the 'identity' in the above notation.\par
The following Transfer Theorem, which computes the transfer  $\tau_m$ on the Tautological Module, is a key to our inductive computation of Chern numbers. The formula involves, in the case of node classes, the analogous transfer operation in lower degree on boundary families $X^\theta_T$, and may be viewed as a commutation relation between the transfer operator  and the node scroll and
node section \emph{operators}. The relation involves the operators
$e^{n,m}_j$ of \eqref{e-operator-eq}.
\begin{thm}
\label{taut-tfr} ({\bf{Tautological transfer}}) $\tau_m$ takes
tautological classes on $X\sbr m-1._B$ to tautological classes on
$X\sbr m._B$. More specifically we have, for any class $\beta\in A.(X)$:
\begin{enumerate}\item
for any twisted polyblock diagonal class $\Gamma_\mu[\alpha.],
\alpha\in \ts_\mu(H^.(X)), w(\mu)=m-1$, \beql{}{
\tau_m(\Gamma_\mu[\alpha.]\beta\subp{m})=\Gamma_{\mu+1_1}
[\alpha.\beta] } where $1_1$ is the distribution of weight 1 and support
$\{1\}$ and $\alpha.\beta$ is the formal symmetric multiplication;
\item for any twisted node scroll class
    $F_{j}^{n,m-1}(\theta)[\alpha]$, $\alpha\in
    T^{m-n-1}(X^\theta_T)$,
\eqspl{tfr-scroll}{\tau_m(F^{n,m-1}_j(\theta)[\alpha]\beta\subp{m})=\
F^{n,m}_j(\theta)[\tau_{m-n,X^\theta_T/T}(\alpha\otimes(\beta|_{X^\theta_T})]}
\item for any twisted node section $-\Gamma\spr
    m-1..F^{n,m-1}_j(\theta)[\alpha]$,

\eqspl{tfr-section}{\tau_m(-\Gamma\spr m-1.F^{n,m-1}_j(\theta)[\alpha]\beta_{(m)})&=\\
&\theta^*(\beta) F^{n+1,m}_j(\theta)[\alpha]+
(-\Gamma\spr m.)F^{n,m}_j(\theta)[\tau_{m-n,X^\theta_T/T}(\alpha.\beta|_{X^\theta_T})\\
&-F^{n,m}_j(\theta)[e^{n,m}_{j+1}(\theta)(\tau_{m-n,X^\theta_T/T}(\alpha.\beta|_{X^\theta_T}))]
+F^{n,m}_j(\theta)[\tau_{m-n,X^\theta_T/T}(e^{n,m-1}_{j+1}(\theta)(\alpha).\beta|_{X^\theta_T})]
}
\end{enumerate}\end{thm}

\begin{proof}
Part (i) is obvious. As for Part (ii), the flatness of $q_{m-1}$ allows us to
work over a general $z\in F$ and then Corollary \ref{q-flat}, (ii) allows us to
assume that the added point is a general point on the fibre $X_s$, which
leads to \refp{tfr-scroll}..\par As for (iii), we recall Corollary 8.4 of
\cite{structure}, which states (using our current notation, not consistent with notation there) that on $F^{n,m-1}_j(\theta)$,
if we denote by $Q_j^{n,m-1}$ the canonical cross-section
$\P(D_j^{n,m-1})$, and by $\Gamma\spr m-n-1.$ the pullback of the discriminant from $(X^\theta)_T\sbr m-n-1.$, then
we have
\[-\Gamma\spr m-1.+\Gamma\spr m-n-1.\sim Q^{n,m-1}_j + p_{[m-n-1]}^*(D^{n,m-1}_{j+1})\sim
Q^{n,m-1}_{j+1} + p_{[m-n-1]}^*(D^{n,m-1}_{j}).\]
Hence,,
\[-\Gamma\spr m-1.\sim Q^{n,m-1}_j + e^{n,m-1}_{j+1}.\]

Similarly, on $F^{n,m}_j(\theta)$, we have
\[-\Gamma\spr m.+\Gamma\spr m-n.\sim Q^{n,m}_j + p_{[m-n]}^*(D^{n,m}_{j+1})
\sim Q^{n,m}_{j+1} + p_{[m-n]}^*(D^{n,m}_{j}),\]
 hence
\[-\Gamma\spr m.\sim Q^{n,m}_j + e^{n,m}_{j+1}\sim Q^{n,m}_{j+1}+e^{n,m}_j.\]

 Therefore, it will
suffice to prove that
\eqspl{tfr-Q}{\tau_m(Q_{j+1}^{n,m-1}[\alpha]\beta\subpar{m})= \theta^*(\beta)
F^{n+1,m}_{j+1}[\alpha]+ Q_{j+1}^{n,m}[\alpha\beta] } and similarly for $j$,
which case is similar (see below). It will suffice to prove this without the
$\alpha,\beta$ twisting.

\par
To this end, note that, with $Q=Q^{n,m-1}_{j+1}$,
$q_{m-1}^*Q$
splits in two
parts, depending on whether the point $w$ added to a scheme $z\in Q$ is
in the off-node or nodebound portion of $z$. It is easy to see that the first
part gives rise to the 2nd  term in the RHS of \refp{tfr-Q}..
\par The analysis of the other part, which leads to the first summand in the RHS of \eqref{tfr-Q}
 is a bit more involved. In substance, what has to be proved in the case at hand is that
 $F^{n+1,m}_{j+1}$ appears with coefficient equal to 1.
To begin with, it is easy to see that we may assume $m=n+1$,
in which case $F$ is just a $\P^1$, namely $C^{m-1}_j$.
Now referring to \eqref{C-mm-1}, the nodebound portion of
$q_{m-1}\inv(C^{m-1}_j)$, as a set, is given by
$\td C^m_j\cup \td C^{m-1}_{j}\cup\td C^m_{j+1}$
 and that of $q_{m-1}\inv(Q)$ is $\td C^m_{j+1}$.
 It will now suffice to show that the 1-dimensional cycle $q_{m-1}\inv(Q)$
 contains $\td C^m_{j+1}$ with multiplicity 1.\par
The latter assertion will be an elementary consequence of the
equations on p. 440, l. 9-14 of \cite{Hilb}, describing the local model
$H_{\mm}$, as well as those on p. 433, describing the analogous local
model $H_m$, to which equations we will be referring constantly in the
remainder of the present proof. Note that $c_{m-i}$ (resp. $b'_{i-1}$) plays
the role of the affine coordinate $u_i/v_i$ (resp. $v'_{i-1}/u'_{i-1}$). Also
our $j+1$ is the $i$ there. We work on $q_{m-1}\inv(C^{m-1}_j)$. Now to
complete the proof, it will suffice to prove\par
\emph{Claim } :  In a
neighborhood of the point
\[(Q, Q^m_{j+1})=(Q^{m-1}_{j+1}, Q^m_{j+1})\in
 \td C^{m-1}_j\cap\td C^m_{j+1}\subset X\sbr m,m-1._B,\] $q_{m-1}\inv(Q)$ contains $\td C^m_{j+1}$ with
multiplicity 1.\par
 To see this note that the defining equations
of $C^{m-1}_j$ on $X\sbr m-1._B$ are given by setting all $a'_k$ and $d'_k$,
as well as $c'_{m-i-1}$ to zero . By loc. cit. p.433
l.9, this implies that we have $b'_1=...=b'_{i-2}=0$
on $q_{m-1}\inv(C^{m-1}_j)$ as well. At a general
point of $C^m_{j+1}$,  $c_{m-i}$ is nonzero. Therefore we may
consider $c_{m-i}$ as a unit. By loc. cit. p.440, eq. (15), we conclude $a_{m-i}=0$.
From this we see easily
that all $a_k=d_k=0$ except $d_{i-1}$, which is a local equation for
$\td C^m_{j+1}$, while $b'_{i-1}$ is a coordinate along $C^{m-1}_j$
having $Q^{m-1}_{j+1}$ as its unique zero. Now by p.440 l. 14, $b'_{i-1}$
and $d_{i-1}$ differ by the multiplicative unit $-c_{m-i}$,
therefore $b'_{i-1}$  cuts
out  $\td C^m_{j+1}$ with multiplicity 1,
which proves  our Claim.\end{proof}
For many purposes, it is possible and even more convenient
 to use \eqref{tfr-Q} directly, rather than through Theorem\ref{taut-tfr}, to compute the transfer
 on node sections. In particular, this form works better with
iterated node scroll/$Q$-sections. To state the result, it is
 convenient to introduce the following notation for the
 iterated scroll/$Q$-sections:
 \eqspl{}{
(FQ)^{n_1,...,n_r;n'_1,...,n'_s;m}_{j_1,...,j_r;j'_1,...,j'_s}(\theta_1,...,\theta_r;
\theta'_1,...,\theta'_s)[\alpha] :=
F^{n_1,m}_{j_1}(\theta_1)[...[Q^{n'_1,m-n_1-...-n_r}_{j'_1}(\theta'_1)...[\alpha]...]
}where $\theta_1,...,\theta'_s$ is a collection of distinct
nodes. Then we have
\begin{cor}
\eqspl{}{
&\tau_m((FQ)^{n_1,...,n_r;n'_1,...,n'_s;m-1}_{j_1,...,j_r;j'_1,...,j'_s}(\theta_1,...,\theta_r;
\theta'_1,...,\theta'_s)[\alpha]\beta_m)=\\
&(FQ)^{n_1,...,n_r;n'_1,...,n'_s;m}_{j_1,...,j_r;j'_1,...,j'_s}(\theta_1,...,\theta_r;
\theta'_1,...,\theta'_s)[\tau_{m-|n.|-|n'.|}(\alpha)\beta_m] \\
&+\sum\limits_{i=1}^s (\theta'_i)^*(\beta)(FQ)^{n_1,...,n_r,
n'_i+1;n'_1,...,\hat{n'_i},...,n'_s;m}_{j_1,...,j_r, j'_i;
j'_1,...,\hat{j'_i},...,j'_s}(\theta_1,...,\theta_r, \theta'_i;
\theta'_1,...,\hat{\theta'_i},...,\theta'_s)[\tau_{m-|n.|-|n'.|}(\alpha)\beta_m]
}
\end{cor}

\newsubsection{Transfer and Chern numbers}\label{transfer-chern-sec}
We are now ready to tackle the computation of Chern numbers,
and in fact all polynomials in the Chern classes
of the tautological bundle on the relative Hilbert scheme
$X\sbr m._B$. The computation is based on passing from
$X\sbr m._B$ to the corresponding full-flag Hilbert scheme
$W=W^m(X/B)$ studied in \cite{R} and a \emph{
diagonalization theorem} for the total Chern class of
(the pullback of) a tautological bundle on $W$, expressing it
either as a simple
(factorable) polynomial in diagonal classes
induced from the various $X\sbr n._B,\; n\leq m$, plus base classes, or, more conveniently,
as the product of the Chern class of a smaller tautological bundle and a diagonal class.
Given this, we can compute Chern numbers essentially by repeatedly applying
the transfer calculus of the last section.\par
We start by reviewing some results from \cite{R}.
Let
 \beq W^m=W^m(X/B)\stackrel{\pi\spr m.}{\longrightarrow} B \eeq denote the relative
 flag-Hilbert scheme of $X/B$, parametrizing flags
 of subschemes \beq z.=(z_1<...<z_m) \eeq
 where $z_i$ has length $i$ and $z_m$ is contained in some fibre
 of $X/B$. Let \beq w^m:W^m\to X\sbr m._B, w^{ m,i}.:W^m\to X\sbr i._B \eeq be the canonical
 (forgetful) maps. Let \beq a_i:W^m\to X \eeq be  the canonical map
 sending a flag $z.$ to the 1-point support of $z_i/z_{i-1}$
 and \beq a^m=\prod a_i:W^m\to X^m_B \eeq their (fibred) product,
 which might be called the 'ordered cycle map'.
 Let \beq \I_m<\O_{X\sbr m._B\times_BX}\eeq
 be the universal ideal of colength $m$.
 For any coherent sheaf on $X$, set
 \beq \Lambda_{m}(E)=p_{{X\sbr m._B*}}(p_{X}^*(E)\otimes
 (\O_{X\sbr m._B\times_BX}/\I_m))\eeq
  These are called the
 \emph{tautological sheaves} associated to $E$; they are locally free if $E$ is. Abusing notation,
 we will also denote by $\Lambda_m(E)$ the pullback of the
  tautological sheaf to appropriate flag Hilbert schemes
  mapping naturally to $X\sbr m._B$, such as $W^m$ or $X\sbr
  \mm._B$. With a similar convention, set
  \beql{}{\Delta\spr m.=\Gamma\spr m.-\Gamma\spr m-1.
  } (recall that $\Gamma\spr m.$ is half the 'physical' discriminant and becomes effective and reduced on $W^m$; thus $\Delta\spr m.$ is an effective (integral) divisor).
 The various tautological sheaves form
 a flag of quotients on $W^m$:
 \beql{fu--flag}{...\Lambda_{m.i}(E)\twoheadrightarrow \Lambda_{m,i-1}(E)
 \twoheadrightarrow...}
 This flag makes possible a simple formula for the total Chern class of the
 tautological bundles, namely the following
 \emph{diagonalization theorem} (\cite{R}, Cor. 3.2):
 \begin{thm}\label{tautbun-W}
 The total Chern class of the tautological bundle $\Lambda
 _m(E)$ pulled back to $W^3(X/B)$ is given by
 \beql{}{c(\Lambda_m(E))=\prod_{i=1}^m c(a_i^*(E)(-\Delta\spr i.))
 }
 \end{thm}
 An analogue of this, more useful for our purposed, holds already on the flaglet
 Hilbert scheme. It can be proved in the same way, or as an easy
 consequence of Thm \ref{tautbun-W}
 \begin{cor}\label{tautbun}
 We have an identity in $A.(X\sbr \mm._B)_\Q$:
 \beql{}{c(\Lambda_m(E))=c(\Lambda_{m-1}(E))c(a_m^*(E)
 (-\Delta\spr m.)).
 }

 \end{cor}
 \begin{proof}
 By Theorem \ref{tautbun-W}, the RHS and LHS pull back to the same class in
 $W^m$. As the projection $W^m\to X\sbr \mm._B$ is generically finite, they
 agree mod torsion.
 \end{proof}
    \begin{rem}
   If $E$ is a line bundle, then it is easy to see from Theorem
   \ref{tautbun-W} that
   $$c_1(\Lambda_m(E))=[m]_*c_1(E)-\Gamma\spr m.=
   \Gamma_1[c_1(E)]-\Gamma\spr m..\qed$$
   \end{rem}
\begin{example}\label{chern-lambda-3-example} On $W^3$:
\eqsp{&c(\Lambda_3(L))=1+L_1+L_2+L_3-\Gamma\spr 3.\\
&+L_1L_2+L_1L_3+L_2L_3-(\Gamma\spr 2.)^2
+\Gamma\spr 2.\Gamma\spr 3.+\Gamma\spr 2.L_2
-\Gamma\spr 2.L_3-\Gamma\spr 3. L_1-\Gamma\spr 3.L_2\\
&+L_1L_2L_3-\Gamma\spr 2.L_1L_3-\Gamma\spr 3.L_1L_2
+\Gamma\spr 2.L_1L_2+\Gamma\spr 2.\Gamma\spr 3.L_1-
(\Gamma\spr 2.)^2L_1\qed
}
\end{example}

More generally, again as a consequence of Theorem \ref{tautbun-W}, we obtain an explicit formula for the
Chern classes of $\Lambda_m(L)$ on $X\sbr m._B$, for a line bundle $L$. First some notation. For a distribution $\mu$, set
\[\chi(\mu)=\prod\limits_n ((n-1)!)^{\mu(n)}.\]
\begin{cor}
For a line bundle $L$, we have
\[c(\Lambda_m(L))=\sum\limits_{|\mu|\leq m}
(-1)^{|\mu|-\l(\mu)|}\frac{a(\mu)\chi(\mu)}{|\mu|!(m-|\mu|)!}
\Gamma_\mu[(1+[L])^{(\mu)}].\]
\end{cor}
\begin{proof}
Straightforward from Theorem \ref{tautbun-W}, see \cite{R2},
Theorem 4.2 (though
the formula there is slightly misstated).
\end{proof}
\begin{rem}
In the classical situation of a single smooth curve over a point,
multiplying diagonal classes is elementary. For example, in the case of
single-block classes, we have
\eqspl{}{
\Gamma\subp{n}.\Gamma\subp{\l}=2\binom{n}{2}\binom{\l}{2}\Gamma\subp{n+l-2}
[-\omega]+n\l\Gamma\subp{n+\l-1}+(1+\delta_{n,\l})\Gamma\subp{n|\l}
.}
Moreover in the notation of Remark \ref{Theta}, we have
\[\Gamma\subp{n+l-2}
[-\omega]=(2-2g)\theta^{n+\l-2}.\]
This allows computation, modulo the combinatorics, of all Chern polynomials of $\Lambda_m(L)$. However, Macdonald's formulation, which expresses everything as polynomials of $\Gamma$ and $\theta$, is probably more efficient in this case
(no singularities).

\end{rem}
\begin{example}
In the case of a good pencil, we can give an explicit formula for
the polyblock (singularity-free) portion of $(\Delta\spr m.)^k\Gamma_\mu$ in $X\sbr m,m-1._B$, as follows. Define for a partition $\mu=(n_1, n_2,...)$,
\[f_k(\mu)=(\sum n_i)^k-\sum\limits_j(\sum\limits_{i\neq j}
n_i)^k+\sum\limits_{j<j'}(\sum\limits_{i\neq j,j'}n_i)^k\pm...\]
and note that this can be identified with the sum of all terms
$\frac{k!}{\alpha_1!...\alpha_\l!}n_1^{\alpha_1}...n_\l^{\alpha_\l}$
where each $\alpha_i$ is $>0$. Then
\eqspl{}{
(\Delta\spr m.)^k.\Gamma_\mu\equiv\sum\limits_{\mu=\mu_1+\mu_2}
f_k(\mu_2)
\Gamma_{\mu_1}\star\Gamma_{|\mu_2|+1}[(-\omega)^{m-\l(\mu_2)+1}]
\mod {\mathrm{node\ classes}}} To see this expand $(\Delta\spr m.)^k$ as a multinomial in the
$D_{i,m}$ and for each monomial $M$
 break up $\mu$ as $\mu_1+\mu_2$ where no block (resp. every block) of $\mu_1$ (resp. $\mu_2$) has an element occurring
 in $M$. The product of $\Gamma_{\mu_2}$ with $M$ yields $\Gamma_{|\mu_2|+1}[(-\omega)^{m-\l(\mu_2)+1}]$, whence the result.
\end{example}
 Motivated by the Corollary we make the following definition.
\begin{defn} Let $R$ be a $\Q$-subalgebra of $A(X)$ containing
 1, the canonical class $\omega$ and the classes of all
 the marked points. The Chern tautological ring on $X\sbr m._B$,
denoted
$$TC^m_R=TC^m_R(X/B),$$ is the $R$-subalgebra of $A(X\sbr m._B)_\Q$ generated by
  the Chern classes of $\Lambda_m(E)$ and the discriminant class $\Gamma\spr m.$.
   \end{defn}

   The following is the main result of this paper.
  \begin{thm}
  There is a computable inclusion
  \beql{}{TC^m_R\to T^m_R.}
  More explicitly, any polynomial in the Chern classes of $\Lambda_m(E)$,
  in particular the Chern numbers, can be computably expressed as a linear
  combination of standard
  tautological classes: twisted diagonal classes, twisted node scrolls,
  and twisted node sections.
  \end{thm}
\begin{proof} For $m=1$ the statement is essentially vacuous.
For $m=2$ it is a consequence of the Module Theorem
\ref{taut-module}. For general $m$, we
assume inductively
the result is true for $m-1$. Given any polynomial $P$ in the
Chern classes of $\Lambda_m(E)$, Corollary \ref{tautbun} implies
that we can write its pullback on $X\sbr\mm._B$ as a sum of terms of
the form $p_{X\sbr m-1._B}^*(Q).(\Gamma\spr m.)^k.S$ where $Q\in
TC^{m-1}_R$. By induction, $Q\in T^{m-1}_R$, so by the Transfer Theorem
\ref{taut-tfr}, $\tau_m(Q)\in T^m_R$. By the projection formula and
the Module Theorem \ref{taut-module}, it follows that $P\in T^m_R$.
\end{proof}
\begin{rem}
This result suggests the natural question: is $T^m_R$ a ring? more
ambitiously, is the inclusion $TC^m_R\to T^m_R$ an equality?
\end{rem}
\newsubsection{Punctual transfer and Pl\"ucker formulas}\label{punctual-transfer-sec}
There is a useful variant of the transfer operation for punctual schemes, i.e. those supported at a single point, which are
parametrized by the small diagonal $\Gamma_{(m)}\subset X\sbr m._B$. This yields a quicker way to compute Chern classes and Chern numbers
of tautological bundles over the small diagonal (compared with
computing the analogous objects over the full Hilbert scheme and
restricting).
Working with the smaller-dimensional $\Gamma_{(m)}$ allows us to obtain geometrically meaningful numbers from the Chern and Segre classes themselves without considering higher-degree polynomials, resulting in some Plu\"cker-type formulas
which, unlike in the case of the full Hilbert scheme, we are able to give in closed form.
Note that unlike the Hilbert scheme $X\sbr m._B$,  the small diagonal $\Gamma\subp{m}$ is generally singular.\par
  This  transfer is based on the correspondence
\eqspl{}{
\begin{matrix}
&&\tilde\Gamma_{(m)}&&\\
&\swarrow&&\searrow&\\
\Gamma_{ (m)}&&&&\Gamma_{(m-1)}
\end{matrix}
}
where the punctual flaglet Hilbert scheme
$\tilde\Gamma\subp{m}$ is defined by the Cartesian diagram
\eqspl{}{
\begin{matrix}
\tilde\Gamma\subp{m}&\to &X\sbr m,m-1._B\\
\downarrow&&\downarrow\\
\Gamma\subp{m}&\to&X\sbr m._B
\end{matrix}
}
As in \S \ref{smalldiag} and Corollary \ref{flaglet-smalldiag-cor}, $\tilde\Gamma\subp{m}$ can be identified locally
over each node with the unique
dominant component of the fibre product over $X$
\[\prod\limits_{i=1}^{m-1}(X_{m,i}/X)\times_X\prod\limits_{i=1}^{m-2}
(X_{m-1,i}/X)\]
where each $X_{n,i}$ is the blowup of the ideal $(x^{n-i}, y^i)$.\
 We denote by
\[\tau_m^0:A^.(\Gamma_{(m-1)})\to A.(\Gamma_{(m)})\]
the induced map, i.e. $p_{[m]*}p^*_{[m-1]}$ from Chow cohomology. Also, let $C^m_i(\theta), C^{m-1}_i(\theta)$
denote the proper transforms of the appropriate node
scrolls ($\P^1$-bundles over $\theta$) from $\Gamma_{(m)},
\Gamma_{(m-1)}$, respectively (for simplicity of notation, we will
drop the tilde over the $C^\bullet_\bullet$ used in earlier sections), and
$C^m_i=\sum\limits_\theta C^m_i(\theta), C^{m-1}_i=\sum\limits_\theta C^{m-1}_i(\theta)$.
$C^m_i$ contains the distinguished sections $Q^m_i, Q^m_{i+1}$.
The following result is proved similarly to Proposition \ref{excdiv-smalldiag-prop}:
\begin{lem}\label{disc-on-small-flag-lem}
We have \eqspl{}{\Gamma\spr m..\tilde\Gamma_{(m)}=
-\binom{m}{2}\omega+\sum
\limits_{i=1}^{m-1} \nu_{m,i}C^m_i
+\sum\limits_{i=1}^{m-2}\nu_{m,i}^-C^{m-1}_i}
\eqspl{}{\Gamma\spr m-1..\tilde\Gamma_{(m)}=
-\binom{m-1}{2}\omega+\sum\limits_{i=1}^{m-2} \nu_{m-1,i}C^{m-1}_i
+\sum\limits_{i=1}^{m-1}\nu_{m-1,i}^+ C^{m}_i}
where
\eqspl{}{
\nu_{m,i}=\frac{i(m-i)m}{2}, 1\leq i\leq m-1,    \\
\nu_{m,i}^-=\frac{i(m-i-1)(m+1)}{2}, 1\leq i\leq m-2\\
\nu_{m-1,i}^+=\frac{i(m-i)(m-2)}{2}, 1\leq i\leq m-1.
}\end{lem}
\begin{proof}
The coefficient of $C^m_i$ in  $\Gamma\spr m..\tilde\Gamma_{(m)}$
is already computed in Propoosition \ref{excdiv-smalldiag-prop}.
As for the coefficient of $C^{m-1}_i$ in  $\Gamma\spr m..\tilde\Gamma_{(m)}$, it suffices to note that $C^{m-1}_i$ contracts
to $Q^m_{i+1}$ in $\Gamma\subp{m}$, where $\Gamma\spr m.$ is defined by $x^{\binom{m-i}{2}}y^{\binom{i+1}{2}}$, and $x,y$ have respective multiplicities equal to $i,(m-1-i)$ on $C^{m-1}_i$,
so the coefficient of $C^{m-1}_i$ equals
\[\binom{m-i}{2}i+\binom{i+1}{2}(m-i-1)= \frac{i(m-i-1)(m+1)}{2}.\]
The other
formula is proved similarly.
\end{proof}
\begin{cor}
\eqsp{(\Gamma\spr m.-\Gamma\spr m-1.).\Gamma\subp{m}=
-(m-1)\omega+\sum\limits_{i=1}^{m-1}i(m-1)C^m_i+\sum\limits_{i=1}^{m-2}
i(m-1-i)C^{m-1}_i.}
\end{cor}
Now set
\eqspl{}{
\nu_m=\sum\limits_{i=1}^{m-1}\nu_{m,i}=\frac{m^2(m^2-1)}{12}
}
\eqspl{}{\nu_m^-:=\sum\limits_{i=1}^{m-2}\nu_{m,i}^-
=\sum\limits_{i=1}^{m-1}\nu_{m-1,i}^+:=\nu_{m-1}^+\\
=\frac{m(m-2)(m^2-1)}{12}
}
We begin with a couple of simple examples that don't require
the full force of the transfer process.
\begin{example}\label{segre-example} Consider a family  $X/B$ of arbitrary base dimension,
with a map $f:X\to \P^m$. We wish to enumerate the locus of $m$-contact $\P^{m-2}$'s
in the family (e.g. cusps ($m=2$), inflexional tangent lines ($m=3$) etc.). More precisely, this is the locus of punctual
length-$m$ subschemes of fibres whose image under $f$ is contained in some $\P^{m-2}$; when (and only when) the subscheme is supported
at a point $x$ smooth on its fibre, the subscheme is determined by $x$ and denoted $mx$.
This locus, of expected codimension 2,
  is given by a degeneracy locus of a map over
$\Gamma\subp{ m}$:
\[(m+1)\O\to \Lambda_m(L)\]
By Porteous, the locus of these is given by
\[s_{2,m}=(c_1^2-c_2)(\Lambda_m(L))\cap[\Gamma\subp{m}]\]
On $\tilde\Gamma_{(m)}$, we can write
\[c(\Lambda_m(L))=c(\Lambda_{m-1}(L))(1+[L]+\Gamma\spr m-1.-\Gamma\spr m.)\]
Therefore,
\eqsp{& s_{2,m}-\tau_m^0(s_{2, m-1})=\\
&c_1(\Lambda_{m-1}(L))(L+\Gamma\spr m-1.-\Gamma\spr m.)+(L+\Gamma\spr m-1.-\Gamma\spr m.)^2\\
&=(L+\Gamma\spr m-1.-\Gamma\spr m.)(mL-\Gamma\spr m.)\\
&=mL^2-(m+1)L\Gamma\spr m.+mL\tau_m(\Gamma\spr m-1.)
-\Gamma\spr m.\tau_m(\Gamma\spr m-1.)+(\Gamma\spr m.)^2\\
&=mL^2+(m+1)\binom{m}{2}L\omega
-(m+1)\sum\limits_{i,\theta}\nu_{m,i}C^m_i(\theta)[\theta^*L]
-m\binom{m-1}{2}L\omega+m\sum\limits_{i,\theta}\nu_{m-1,i}^+C^m_i
(\theta)[\theta^*L]\\
&-\binom{m}{2}\binom{m-1}{2}\omega^2-\Gamma\spr m.
\sum\limits_{i,\theta}\nu_{m-1, i}^+C^m_i(\theta)\\&+\binom{m}{2}^2\omega^2+\Gamma\spr m.\sum
\limits_{i,\theta}\nu_{m,i}
C^m_i(\theta)}
Therefore
\eqspl{segre-smalldiag-eq}{
& s_{2,m}-\tau_m^0(s_{2, m-1})=\\
&=mL^2+ m(m-1)L\omega+(m-1)\binom{m}{2}\omega^2
-\sum\limits_{i,\theta}\frac{3i(m-i)m}{2}C^m_i(\theta)[\theta^*L]+\Gamma\spr m.\sum\limits_{i,\theta} i(m-i) C^m_i(\theta)
}
If $\dim(B)=1$, then $s_{2,m}$ can be viewed as a number, and \eqref{segre-smalldiag-eq} gives a recursion for
it, one which can be easily solved explicitly, as follows.
Note that in this case we have $\theta^*L=0$ and $\Gamma\spr m.C^m_i(\theta)=-1$, so \eqref{segre-smalldiag-eq} simplifies to
\eqsp{s_{2,m}-s_{2, m-1}&=mL^2+m(m-1)L\omega+(m-1)\binom{m}{2}\omega^2-
\sigma\sum\limits_{i,\theta} i(m-i)\\&
=mL^2+m(m-1)L\omega+(m-1)\binom{m}{2}\omega^2-\sigma\frac{m(m^2-1)}{6}.
}
This recursion is easily integrated, yielding the following
 Pl\"ucker-type formula in closed form:
\eqspl{segre-smallgiag-explicit-eq}{
s_{2,m}=\binom{m+1}{2}L^2+\frac{m(m^2-1)}{2}L\omega+\frac{m(m^2-1)(3m-2)}{24}\omega^2
-\frac{m(m^2-1)(m+2)}{24}\sigma
}
e.g.
\[s_{2,1}=L^2,  s_{2,2}=2L^2+3L\omega+\omega^2-\sigma,
 s_{2,3}=6L^2+12L\omega+7\omega^2-5\sigma,...\]

In case $\dim(B)>1$, \eqref{segre-smalldiag-eq} must be combined with the punctual transfer calculus (Proposition \ref{punctual-transfer-prop} below) to yield the recursion.\qed
\end{example}
Now we take up the punctual transfer proper.
To be precise, let $T^{0,m}_R(X/B)$ denote the group generated
by $R\subset A^.(X)$, the twisted node scrolls $C^m_i(\theta)[\beta],
\beta\in R$ on $\Gamma\subp{m}$, and their sections
$(-\Gamma\spr m.)C^m_i(\theta)[\beta]$.  Then we will define a 'pointwise transfer' map
\[\tau_m^0:T^{m-1,0}_R(X/B)\to T^{m,0}_R(X/B)\]
that fits in the diagram
\eqspl{}{
\begin{matrix}
T^{m-1,0}_R(X/B)&\to &T^{m,0}_R(X/B)\\
\downarrow&&\downarrow\\
A^\bullet(\Gamma\subp{m-1})&\to&A_\bullet(\Gamma\subp{m}).
\end{matrix}
} Recall that we are assuming $R$ contains $\Q$; this is essential here as $\tau^0_m$ is only defined over $\Q$.
Set
\eqspl{}{
\psi^m_i:=\psi_x^{\otimes\binom{m-i+1}{2}}\otimes\psi_y^{\otimes\binom{i}{2}}
} These obviously depend on $\theta$ and will be denoted
$\psi^m_i(\theta)$ when necessary.\par
The following result gives the main rules of the punctual
transfer calculus.
\begin{prop}\label{punctual-transfer-prop}
Notations as above, we have, for each node $\theta$,
\eqspl{punctual-transfer-on-C-eq}{
\tau_m^0(C^{m-1}_i(\theta))=\frac{m-i}{m-1}C^m_i(\theta)
+\frac{i+1}{m-1}C^m_{i+1}(\theta).
}
\eqspl{punctual-transfer-on-Gamma.C-eq}{
\tau^0_m&(-\Gamma\spr m-1..C^{m-1}_i(\theta))=
\frac{m-i-1}{m}Q_{i+1}^m(\theta)+\frac{i+1}{m}Q^m_{i+2}(\theta)
+ \frac{m-i}{m-1}C^m_i(\theta)[\psi^{m-1}_i]
+\frac{i+1}{m-1}C^m_{i+1}(\theta)[\psi^{m-1}_{i}]\\
&=-\Gamma\spr m..C^m_{i+1}(\theta)-C^m_{i+1}(\theta)[\frac{m-i-1}{m}
\psi^m_{i+2}+\frac{i+1}{m}\psi^m_{i+1}]+\frac{m-i}{m-1}C^m_i(\theta)[\psi^{m-1}_i]
+\frac{i+1}{m-1}C^m_{i+1}(\theta)[\psi^{m-1}_i]
}
\end{prop}
\begin{rem}\label{compare-q-gamma}
Because $Q^m_i=\P(\psi^m_i)$ and $Q^m_{i+1}$ are disjoint sections
on the $\P^1$-bundle \mbox{$C^m_i=\P(\psi^m_i\oplus\psi^m_{i+1})$,}
we have
\eqspl{qi-qi+1-comparison}{
Q^m_{i+1}\sim Q^m_i+C^m_i[\psi^m_{i+1}-\psi^m_i]= Q^m_i+C^m_i[-(m-i)\psi_x+i\psi_y],
}
\eqspl{gamma-q-comparison}{
-\Gamma\spr m..C^m_i\sim &Q^m_i+C^m_i[\psi^m_{i+1}]\\
\sim &Q^m_{i+1}+C^m_i[\psi^m_{i}]
}
\end{rem}\begin{rem}\label{punctual-transfer-pencil-rem}
If the base $B$ is 1-dimensional, both $-\Gamma\spr m-1.C^{m-1}_i(\theta)$
and  $-\Gamma\spr m.C^{m}_i(\theta)$ are points over the finite set $B(\theta)$
and in particular
\[\tau^0_m(-\Gamma\spr m-1.C^{m-1}_i(\theta))=-\Gamma\spr m.C^{m}_i(\theta)
\]

\end{rem}
\begin{proof}
Fixing (and suppressing) $\theta$,
we analyze $\Gamma\subp{m}$ locally over $B$ and near
the point
$m\theta$ in the small diagonal $ X\subset X\spr m._B$,
as in \S\ref{smalldiag}.
As we have seen, $\Gamma\subp{m}$ is given by the blowup of
\[ J_m=(g_i:=x^{\binom{m-i+1}{2}}y^{\binom{i}{2}}:i=1,...,m)
=\prod\limits_{i=1}^{m-1}(x^{m-i}, y^i)\]
and embeds in
\[X\times \prod\limits_{i=1}^{m-1}\P^1_{(u_i,v_i)}\]
where the $i$th factor can be identified with $C^m_i$.
The image map $\Gamma\subp{m}\to X\times C_i^m$
can be identified with the blowup $X_i$ of $X$ in the ideal
\[(g_i,g_{i+1})=x^{\binom{m-i}{2}}y^{\binom{i}{2}}(x^{m-i}, y^i)\]
or what is the same, the blowup of $(x^{m-i}, y^i)$.
The pullback of the exceptional divisor in the blowup of $(x^{m-i}, y^i)$,
locally defined by $x^{m-i}$ or $y^i$, yields a structure of Cartier divisor
of multiplicity (=generic length) $i(m-i)$ on $C^m_i(\theta)\subset \Gamma\subp{m}$. Along $C^m_i$, $\Gamma\subp{m}$ is defined
by the equation
\[x^{m-i}u_i=y^iv_i\]
Also, $C^m_i$ is endowed with the special sections $Q^m_i=(x^{m-i+1}, y_i)$, corresponding to $u_i=0$, and $Q^m_{i+1}=(x^{m-i}, y^{i+1})$, corresponding to $v_i=0$.\par
Globally over $B$, $C^m_i$ is a $\P^1$ bundle of the form
\eqsp{C^m_i=\P(\psi^m_i\oplus\psi^m_{i+1}),
}(using multiplicative notation for line bundles)
with $Q^m_i=\P(\psi^m_i)$,
$Q^m_{i+1}=\P( \psi^m_{i+1})$ and $-\Gamma\spr m.$ corresponding to $\O(1)$. This implies
\eqspl{Gamma-on-Q}{
\O(1).Q^m_j=-\Gamma\spr m..Q^m_j=\psi^m_j, j=i,i+1.
}
\par
Now on $\tilde\Gamma\subp{m}\subset X\times\prod \limits_{i=1}^{m-1}C^m_i\times\prod\limits_{i=1}^{m-2}C^{m-1}_i$, the
exceptional locus is a connected chain of the form
\[ \tilde C^m_1\cup \tilde C^{m-1}_1\cup...\cup
\tilde C^m_i\cup \tilde C^{m-1}_i\cup \tilde C^m_{i+1}\cup...\cup \tilde C^m_{m-1},\]
with each $\tilde C^n_i$ projecting isomorphically to $C^n_i$, and where
the intersection $\tilde C^m_i\cap\tilde  C^{m-1}_i$ (resp. $\tilde C^{m-1}_i\cap
\tilde C^m_{i+1}$) is set-theoretically
the section $(Q^{m-1}_i, Q^m_{i+1})$ (resp. $(Q^{m-1}_{i+1}, Q^m_{i+1})$) (the
multiplicities will be determined below).
From this it follows already that
\mbox{$\tau^0_m(C^{m-1}_{i-1})=a(i-1)C^m_{i-1}+b(i-1)C^m_{i}$,}
and it remains to identify $a$ and $b$. This is a consequence of
the following Lemma, which completes the proof of \eqref{punctual-transfer-on-C-eq}.
\begin{lem}\label{pullback-mult-lem}
We have:
\eqspl{transfer-from-m-1}{
p_{[m-1]}^*((m-1)C^{m-1}_i)=(m-i)\tilde C^m_i+(m-1)\tilde C^{m-1}_i+
(i+1)\tilde C^m_{i+1},
}
\eqspl{transfer-from-m}{
p_{[m]}^*(mC^m_i)=(i-1)\tilde C^{m-1}_{i-1}+m\tilde C^m_i+
(m-1-i)\tilde C^{m-1}_i.
}
\end{lem}
\begin{proof}[Proof of Lemma] We will prove \eqref{transfer-from-m-1}
as the case of \eqref{transfer-from-m} is similar.
We will analyze this near the flag $(Q^{m-1}_i\subset Q^m_{i+1})$, which
is the intersection $\tilde C^m_i\cap\tilde  C^{m-1}_{i}, i\leq m-2$. There, using
homogeneous coordinates $[u'_i, v'_i]$ on $C^{m-1}_i=\P^1$,
$\tilde\Gamma\subp{m}$ has the local equations
\eqspl{}{
x^{m-i}u_i=y^iv_i,\\
x^{m-1-i}u'_{i}=y^{i}v'_i
} which imply
\eqspl{v-times-u'=x}{(v_i/u_i)(u'_i/v'_i)=x.
}
This means that $x$ cuts out a divisor locally equal to the union
of the zero sets of $v_i$ ( i.e. $C^{m-1}_{i}$), and of  $u'_i$ (i.e. $C^m_i$).
Moreover the multiplicity, i.e. local length, of $v_i$, along $C^{m-1}_i$
is equal to that of $x$, which is the length of $\C[x, y]/(x^{m-i}-y^i, x)$, i.e.
$i$. Likewise, the multiplicity of $u'_i$ along $C^m_i$ is also equal to $i$. As we have seen in the proof of Proposition \ref{excdiv-smalldiag-prop},
$u'_i$, i.e. $u'_i/v'_i$, is a local defining equation on
$\Gamma\subp{m-1}$ for the Cartier divisor $(m-1)C^{m-1}_{i-1}$. It follows,
under the assumption $i\leq m-2$,
that in the pullback of $(m-1)C^{m-1}_{i-1}$ to $\tilde\Gamma\subp{m}$,
$C^m_i$ appears with multiplicity $i$. In the extreme case
of $C^{m-1}_{m-2}$, both it and $C^m_{m-1}$ clearly have multiplicity
1 on $\tilde\Gamma_{(m)}$, so $C^m_{m-1}$ appears with multiplicity
$m-1$ in the pullback of $(m-1)C^{m-1}_{m-2}$. 
\par
Similarly using the relation
\eqspl{u-times-v'=y}{ (u_{i+1}/v_{i+1})(v'_i/u'_i)=y,
}
the  pullback of $(m-1)C^{m-1}_{i+1}$,
contains $C^m_{i+1}$with multiplicity $m-i-1$ for all $i<m-2$.

Thus proves \eqref{transfer-from-m-1} by shifting the
index $i$.\par
For \eqref{transfer-from-m} briefly,
the same relations \eqref{v-times-u'=x}
and \eqref{u-times-v'=y} show that the respective multiplicities
of the pullback of $mC^m_{i+1}$ (resp. $mC^m_i$) along $C^{m-1}_i$
are equal to $i$ (resp. $m-1-i$).
\end{proof}
To obtain \eqref{punctual-transfer-on-Gamma.C-eq}, we continue
to develop intersection theory on the punctual flaglet Hilbert scheme,
using the same notations.
\begin{lem}\label{intersect-mult-lem}
We have
\eqspl{}{
\tilde C^{m-1}_i\tilde C^m_i=\frac{1}{i}(Q^{m-1}_i, Q^m_{i+1}),
}
\eqspl{}{
\tilde C^{m-1}_i\tilde C^m_{i+1}=\frac{1}{m-i-1}(Q^{m-1}_{i+1}, Q^m_{i+1}), i<m-1.
}
\end{lem}
\begin{proof}[Proof of Lemma]
For the first relation, the identity \eqref{v-times-u'=x}
shows that the Cartier divisors $i\tilde C^{m-1}_i, i\tilde C^m_i$
meet with multiplicity $i$ along $(Q^{m-1}_i, Q^m_i)$.
For the second relation, use similarly the identity \eqref{u-times-v'=y}.
\end{proof}
\begin{lem}
 We have:
\eqspl{}{
(\tilde C^m_i)^2=-\frac{m-1}{m}(\frac{1}{m-i}(Q^{m-1}_{i},Q^m_i)
+\frac{1}{i}(Q^{m-1}_i,Q^m_{i+1}))
}
\end{lem}
\begin{proof}[Proof of Lemma]
The projection formula, together with Lemma \ref{pullback-mult-lem} shows that
\[p_{[m]*}((\tilde C^m_i)^2)=(C^m_i)^2-p_{[m]*}(\tilde C^m_i.((i-1)\tilde C^{m-1}_{i-1}
+(m-1-i)\tilde C^{m-1}_i)).\]
Combining this with Corollary \ref{C-m-i-intersect-lem}, (ii) and
Lemma \ref{intersect-mult-lem}, we get the result by simple arithmetic.
\end{proof}
Similarly, we can show:
\begin{lem}
\eqspl{}{
(\tilde C^{m-1}_i)^2=-\frac{m}{m-1}(\frac{1}{i}(Q^{m-1}_i, Q^m_{i+1})
+\frac{1}{m-i-1}(Q^{m-1}_{i+1}, Q^m_{i+1})).
}
\end{lem}
Now we are in position to compute the punctual transfer of a section.
We have
\[-\Gamma\spr m-1..C^{m-1}_i=Q^{m-1}_{i+1}+C^{m-1}_i[\psi^{m-1}_i]=
(m-1)C^{m-1}_i.C^{m-1}_{i+1}+C^{m-1}_i[\psi^{m-1}_i].\]
The transfer of the second summand above is computed directly by
\eqref{punctual-transfer-on-C-eq} and yields the last two
terms in the first equality in \eqref{punctual-transfer-on-Gamma.C-eq}.
As for the first summand above, its
pullback on $\tilde\Gamma_{(m)}$ is given by
\[(m-1)(\tilde C^{m-1}_i+\frac{m-i}{m-1}\tilde C^m_i+\frac{i+1}{m-1}\tilde
C^m_{i+1})(\tilde C^{m-1}_{i+1}+\frac{m-i-1}{m-1}\tilde C^m_{i+1}
+\frac{i+2}{m-1}\tilde
C^m_{i+2}).\]
The only terms not clearly trivial come from
$\tilde C^{m-1}_i.\tilde C^m_{i+1}$,  $\tilde C^{m-1}_{i+1}.\tilde C^m_{i+1}$
and $(\tilde C^m_{i+1})^2$,
and those be computed using the above Lemmas,
yielding
\eqspl{transfer-as-q}{\frac{m-i-1}{m}(Q^{m-1}_{i+1}, Q^m_{i+1})+\frac{i+1}{m}(Q^{m-1}_{i+1},
Q^m_{i+2}).}
By projection, we obtain the first
equality in \eqref{punctual-transfer-on-Gamma.C-eq}. The second equality is
immediate from Remark \ref{compare-q-gamma}.

This completes the proof of Proposition \ref{punctual-transfer-prop}.
\end{proof}
\begin{rem}
Using \eqref{transfer-as-q} and
 Remark \ref{compare-q-gamma} again, one can formulate the
punctual transfer
entirely in terms of node scrolls $C^m_i$ and the canonical sections
$Q^m_i$, via
\eqspl{transfer-on-q}{
\tau_m^0(Q^{m-1}_i)&=Q^m_i+ C^m_i[\frac{-i(m-i)}{m}\psi_x+\frac{i^2}{m}\psi_y]\\
&=\frac{m-i}{m}Q^m_i+\frac{i}{m}Q^m_{i+1}
}
\end{rem}
\begin{rem}
By simple numerology, Proposition \ref{punctual-transfer-prop} can be used to reprove
Lemma \ref{disc-on-small-flag-lem}. We omit the details.\qed
\end{rem}
Proposition \ref{punctual-transfer-prop} leads, in a completely straightforward way, to a recursive formula for the total Chern class
of the tautological bundles $\Lambda_m(L)|_{\Gamma\subp{m}}$ for
a line bundles $L$. Write recursively
\eqspl{punctual-taut-generic-eq}
{c_m:=c(\Lambda_m(L)|_{\Gamma\subp{m}})&=\alpha_m+\sum\limits_\theta
(\sum\limits_{i=1}^{m-1} C^m_i(\theta)[\beta^i_m(\theta)]
+\sum\limits_{i=1}^{m-1}\gamma_m^i(\theta)Q^m_i(\theta)),\\
 &\alpha\in R, \beta(\theta), \gamma(\theta)\in A^.(B(\theta)),
 \alpha_1=1+L, \beta^._1=\gamma^._1=0.
 }
 where $B(\theta)\subset B$ is the normalization of the
 boundary divisor corresponding to $\theta$ and a cycle on $B(\theta)$ is viewed as a cycle on $B$
via the Gysin map. In view of Remark \ref{compare-q-gamma},
such expressions are not unique, but this does not matter.
In the ensuing computation we will suppress the
$\theta$ summation, which will be understood. \par
We have:
\[c_m=\tau^0_m(c_{m-1}(1+L+\Gamma\spr m-1.))+
(-\Gamma\spr m.)\tau^0_m(c_{m-1}).\]
The first summand yields
\eqsp{
&(1+L-\binom{m-1}{2}\omega+\frac{1}{m-1}\sum\limits_{i=1}^{m-2}
\nu_{m-1,i}((m-i)C^m_i+(i+1)C^m_{i+1}))\alpha_{m-1}\\
&+\frac{1}{m-1}\sum C^m_i[\beta^i_{m-1}(m-i)(1+L)]+C^m_{i+1}[\beta^i_{m-1}(i+1)(1+L)]\\
&-\frac{1}{m-1}\sum\beta_{m-1}^i\psi^{m-1}_{i+1}((m-i)C^m_i+(i+1)C^m_{i+1}\\
&+\frac{1}{m}\sum\limits_{i=1}^{m-1}((1+L)\gamma_{m-1}^i
-\beta^i_{m-1})((m-i)Q^m_i+iQ^m_{i+1})
}
*************
The second summand yields
\eqsp{&-\sum\nu_{m,i}C^m_i[\alpha_{m-1}]+\binom{m}{2}\omega\alpha_{m-1}\\
&+\frac{1}{m-1}\sum\beta^i_{m-1}((m-i)(Q^m_i+C^m_i[\psi^m_{i+1}])
+(i+1)(Q^m_{i+1}+C^m_{i+1}[\psi^m_{i+2}]))\\
&+\frac{1}{m}\sum\gamma_{m-1}^i((m-i)Q^m_i[\psi^m_i]+iQ^m_{i+1}[\psi^m_{i+1}])
}
 (where $\theta^*$ means pullback by
the node-section $\theta:B(\theta)\to X$    and products involving classes on $B(\theta)$ are performed on $B(\theta)$).\par
Thus in all
\eqspl{punctual-taut-recursion-eq}{
&\alpha_m=(1+L+(m-1)\omega)\alpha_{m-1}=\prod\limits_{i=1}^m
(1+L+(i-1)\omega);\\
&\beta_m^i=\frac{1}{m-1}(((m-i)(1+L-\psi^{m-1}_{i+1})+\psi^m_{i+1})\beta^i_{m-1}
+(i(1+L-\psi^{m-1}_i)+\psi^m_{i+1})\beta^{i-1}_{m-1})-\nu_{m,i}\alpha_{m-1}\\
&\gamma^i_m=\frac{1}{m}((m-i)(1+L+\psi^m_i)\gamma^i_{m-1})+
(i-1)(1+L+\psi^m_i)\gamma^{i-1}_{m-1})+\frac{1}{m(m-1)}((m-i)\beta^i_{m-1}
+(m+i-1)\beta^{i-1}_{m-1})
}
We have used the fact that
$\theta^*(\alpha_{m-1})=\theta^*((1+L)^{m-1})$, which follows from the above formula for the $\alpha_m$, plus the
fact that $\theta^*(\omega)=0.$
We have proven
\begin{cor}\label{punctual-taut-cor}
The Chern classes of the tautological bundle on the punctual
Hilbert scheme $\Gamma\subp{m}(X/B)$ are given by \eqref{punctual-taut-generic-eq}, where the coefficients
satisfy the recursion \eqref{punctual-taut-recursion-eq}
\end{cor}
\begin{example}
Given a family $X/B$ (of any base dimension), and a map \[f:X\to \P^n, n<m,\]
$c_{m-n}(\Lambda_m(f^*(\O(1)))|_{\Gamma\subp{m}}])$ represents
 the locus, finite if $\dim(B)=m-n-1$, of points in $X$
 where the fibre admits an $m$-contact hyperplane. If $n=1$, this is the locus of $(m-1)$st order
 ramification points. If $n=2$, it is the locus of $m$-th order hyperflexes, etc.\par
 The case $n=m-2$ can be worked out more explicitly, as in
 Example \ref{segre-example}.
 Let $L=f^*(\O(1))), c_{2, m}=c_{2}(\Lambda_m(L)|_{\Gamma\subp{m}})$.
 \eqsp{c_{2,m}-\tau_m^0(c_{2, m-1})&=((m-1)L-\Gamma\spr m-1.)(L-\Gamma\spr m.+\Gamma\spr m-1.)\\&
=(m-1)L^2+(m-1)\binom{m-1}{2}\omega^2+(\binom{m-1}{2}+(m-1)^2)L\omega
\\&-\sum\nu_{m-1, i}^+C^m_i[\theta^*L]+(m-1)\sum(\nu_{m-1, i}^+
-\nu_{m,i})C^m_i[\theta^*L]\\&+\Gamma\spr m.\sum\nu_{m-1,i}^+C^m_i-\sum\limits\nu_{m-1, i}\tau_m(\Gamma\spr m-1.C^{m-1}_i)}
Then simple computations and 
Remark \ref{punctual-transfer-pencil-rem} yield in the pencil case
\eqspl{c2-punctual-recursion-eq}{
&c_{2,m}-\tau_m^0(c_{2, m-1})=\\
&(m-1)L^2+(m-1)\binom{m-1}{2}\omega^2+\frac{(m-1)(3m-4)}{2}L\omega\\
&-\sum\limits_\theta\sum\limits_{i=1}^{m-1}
\frac{i(m-i)}{2}\left((3m-4)C^m_i(\theta)[\theta^*L]\right)
-\binom{m}{3}\sigma
}
Again in the pencil
case $\dim(B)=1$ these classes can be viewed as numbers, $\theta^*L$ and $\psi^{m-1}_i$ are all zero and the
above simplifies to
\eqsp{
&c_{2,m}-c_{2,m-1}=\\
&(m-1)L^2+(m-1)\binom{m-1}{2}\omega^2
+\frac{(m-1)(3m-4)}{2}L\omega-\binom{m}{3}\sigma
}
This recursion can be integrated easily using the elementary formula
\mbox{$\sum\limits_{n=k}^m \binom{n}{k}=\binom{m+1}{k+1}$,} yielding the following closed-form Pl\"ucker-type formula:
\eqspl{c2-punctual-explicit-eq}{
c_{2,m}=\binom{m}{2}L^2+(3\binom{m+1}{4}-\binom{m}{3})\omega^2
+(3\binom{m+1}{3}-2\binom{m}{2})L\omega-\binom{m+1}{4}\sigma
}
\end{example}

\newsection{Low-degree  examples}\label{examples}
\newsubsection{Trisecants to one space curve curve}
If $X$ is a smooth curve of degree $d$ and genus $g$ in $\P^3$, the virtual degree of its
trisecant scroll, i.e. the virtual number of trisecant lines to $X$ meeting
a generic line, is given by $c_3(\bigwedge\limits^2 (\Lambda_3(\O_X(1)))$,
which can be easily computed to be
\eqspl{}{
\frac{1}{6}(2d^3-12d^2+16d-3d(2g-2)+6(2g-2))
}
\newsubsection{Multisecants in a pencil }
Let $X/B$ be a family of nodal curves over a smooth curve,
and
 suppose \beq f:X\to\P^{2m-1}\eeq is a
 morphism. One, quite special, class of examples of this situation
 arises as what we call a {\it{generic rational pencil}}; that is,
 generally,
 the normalization of the family of rational curves of fixed degree $d$
  in $\P^r$ (so
 $r=2m-1$ here) that are incident to a generic collection $A_1,
 ...A_k$ of linear spaces, with
 \beq (r+1)d+r-4=\sum(\text{codim}(A_i)-1); \eeq see \cite{R3} and
 references therein, or \cite{RA}.
Our result applies to curves of
 arbitrary genus.
 \par Returning to the general situation,
 one expects a finite number $N_m$ of curves $f(X_b)$ to
 admit an m-secant $(m-2)$-plane.  Let $L=f^*\O(1), V=H^0(\P^{2m-1}), c(m,i)=c_i(\Lambda_m(L)).$ Then $N_m$ is the degree
 of the locus where the natural map $V\otimes\O_{X\sbr m._B}\to\Lambda_m(L)$ drops rank.
 By Porteous' formula
 \cite{ful}, this number can be evaluated as
\[\int\limits_{X\sbr m._B}\Delta\spr m+1._1(c(m,.)).\]
For convenience, we will denote $m!$ times this number simply by $\Delta\spr m+1._1$.
These numbers have been evaluated by completely different means
 (namely, 'test pencils' as in the Harris-Mumford paper)
 by Ethan Cotteril \cite{cot}, \cite{cotteril1}, \cite{cotteril2}. For $m=3$,
our intersection theory yields the same number after prolonged manual calculations:
\eqspl{trisec}{
1/6((3d^2-27d+60)L^2+(-12d+72)L.\omega+(-3d+28)\omega^2
-3b(2g-2)+(3d-20)\sigma)}
where $d$ is the degree of a fibre in $\P^{2m-1}$ and $\sigma$ is the number of singular fibres.\par
More generally, for any 1-parameter family as above, natural numbers $r,s$ with $rs=m+1$ and a map \[f:X\to\P^{(r+1)(s-1)-1},\] we have a 'Porteous number'
\[\Delta\spr s._r=m!\int\limits_{X\sbr m._B}\Delta\spr s._r(c(m,.))=\int\limits_{W^m(X/B)}\Delta\spr s._r(c(m,.))\]
which counts $m!$ times the virtual number of $m$-secant $\P^{m-r-1}$-planes in the family.

Beyond $m=3$, manual calculation via our intersection calculus seems impractical;
fortunately, we have macnodal, discussed next.

\subsection{Gwoho Liu's Macnodal program}\label{macnodal}
Gwoho Liu has written a Java program \cite{macnodal} which implements the above intersection theory ,
both on a fixed Hilbert scheme
$X\sbr m._B$ and, together with the transfer calculus of
\S \ref{transfer}, on a flag scheme $W^m(X/B)$.
This is sufficient to evaluate all Chern numbers of
tautologocal bundles and in particular all the aforementiond multisecant numbers.
We proceed with some examples of Macnodal computations.
Details about Macnodal, and further examples, will be given elsewhere.

\begin{example}[pencils]
As above, let $c(m,i)=c_i(\Lambda_m(L)), L=f^*(\O(1), f:X\to\P^{2m-1}$,
for a good pencil $X/B$. When $m=3,r=1$, Macnodal yields
the formula above. When $m=4, r=1$, it computes the number of
$m$-secant $(m-2)$-planes in the family (times $m!$) as
\eqsp{\Delta^{(5)}_1&=
c(4,1)^5-4c(4,2)c(4,1)^3+3c(4,3)c(4,1)^2+3c(4,2)^2c(4,1)
-2c(4,1)c(4,4)-2c(4,3)c(4,2)\\
&=          (-1008+452d-72d^2+4d^3-24dg+168g)L^2+(432-98d+6d^2-12g)\sigma\\
&\ \ \ \ +(-1440+360d-24d^2+4
8g)\omega.L+(-720+130d-6d^2+12g)\omega^2
}

Similarly, we have\small
\eqspl{}{
&c(5,1)^6=
(1483200-1022400d+280800d^2-10800d^2g-36000d^3+1800d^4+122400dg-363600g+5400g^2)
L.L\\&+(-400800+162000d-23400d^2+1200d^3-3600dg+25800g)\sigma
\\&+(1843200-820800d+12960
0d^2-7200d^3+21600dg-144000g)w.L\\&+(775200-262800d+30600d^2-1200d^3+3600dg-33000g)
w.w}
\normalsize
 This is one of 10 terms in $\Delta\spr 6._1$ (5-secant planes of
 a pencil in $\P^9$):
\eqspl{}
{\Delta_1\spr 6.=\\
&(19560-9270d+1735d^2-60d^2g-150d^3+5d^4+1020dg-4500g+60g^2)L.L\\
&+(-10720+2960d-290
d^2+10d^3-60dg+640g)\sigma\\
&+(33600-10160d+1080d^2-40d^3+240dg-2400g)w.L\\
&+(20000-464
0d+370d^2-10d^3+60dg-800g)w.w
.}
\par For 5-secant planes of
a pencil in $\P^5$:
\eqspl{}{
\Delta_2\spr 3.=\\
&(17400-11070d+2805d^2-120d^2g-330d^3+15d^4+1500dg-4860g+60g^2)L.L\\
&+(-4640+1630d-2
10d^2+10d^3-60dg+440g)\sigma\\
&+(21600-8520d+1200d^2-60d^3+300dg-2160g)w.L\\
&+(9280-267
0d+270d^2-10d^3+60dg-520g)w.w
}
\par 5-secant lines for a pencil in $\P^3$:
\eqspl{}{\Delta^2_3=\\
&(5160-4020d+1250d^2-60d^2g-180d^3+10d^4+600dg-1620g+60g^2)L.L\\
&+(-640+210d-20d^2+40g)\sigma\\
&+(4800-2200d+360d^2-20d^3+60dg-480g)w.L+(1520-450d+40d^2-80g)w.w
}
\par
7-secant $\P^4$-s for a pencil in $\P^9$:\tiny
\eqspl{}{
\Delta_2\spr 4.=\\
&(23420880-15283296d+4337172d^2-592200d^2g+5040d^2g^2-679140d^3+44100d^3g+61425d^
4-1260d^4g-3024d^5+63d^6\\&+3633840dg-95760dg^2-8646120g+468720g^2-2520g^3)L.L
\\&+(-10
425240+4489296d-819630d^2+39060d^2g+78960d^3-1260d^3g-3990d^4+84d^5-418740dg+252
0dg^2+1563240g-28560g^2)\sigma\\
&+(40007520-18793404d+3700620d^2-187740d^2g-379575d^
3+6300d^3g+20160d^4-441d^5+1922760dg-12600dg^2-6804000g+136080g^2)w.L\\
&+(23373000-
8857296d+1396710d^2-46620d^2g-114240d^3+1260d^3g+4830d^4-84d^5+586740dg-2520dg^2
-2536800g+34440g^2)w.w
}\normalsize 
\par
7-secant planes for a pencil in $\P^4$:
\eqspl{}{\Delta_4\spr 2.=\\
&
(2242800-1966440d+750680d^2-109200d^2g+1260d^2g^2-157920d^3+10920d^3g
\\&+19145d^4-4
20d^4g-1260d^5+35d^6+500640dg-17640dg^2-896280g+65520g^2-840g^3)L.L
\\&+(-400680+200
340d-40180d^2+840d^2g+3780d^3-140d^4-12180dg+48720g-840g^2)\sigma\\
&+(2646000-1619
100d+415800d^2-20160d^2g-55335d^3+840d^3g+3780d^4-105d^5
\\&+164220dg-1260dg^2-46620
0g+12600g^2)\omega.L\\
&+(1040760-483420d+90440d^2-1680d^2g-7980d^3+280d^4+25620dg-110040
g+1680g^2)\omega^2
}
\par
8-secant $\P^4$-s of a pencil in $\P^7$:\tiny
\eqspl{}{
\Delta\spr 3._3=\\
&(-524966400+399120960d-135158352d^2+25242000d^2g-564480d^2g^2+26227992d^3\\
&-280728
0d^3g+20160d^3g^2-3129840d^4+159600d^4g+228480d^5-3696d^5g-9408d^6+168d^7\\
&-116556
384dg+5382720dg^2-20160dg^3+222243840g-17539200g^2+201600g^3)L.L\\
&+(193334400-9897
4848d+22152760d^2-1627080d^2g+10080d^2g^2-2761220d^3+102480d^3g+201180d^4\\
&-2520d^
4g-8092d^5+140d^6+11973360dg-215040dg^2-34487040g+1186080g^2-3360g^3)\sigma\\
&+(-8
53493760+482291712d-118440672d^2+9760800d^2g-60480d^2g^2+16089360d^3-651840d^3g\\
&-
1268400d^4+16800d^4g+54768d^5-1008d^6-67085760dg+\\
&1249920dg^2+179141760g-6632640g
^2+20160g^3)\omega.L\\
&+(-469929600+216819456d-43094520d^2+2299080d^2g-10080d^2g^2+46921
00d^3\\
&-122640d^3g-293020d^4+2520d^4g+9884d^5-140d^6-19634160dg+248640dg^2\\
&+6443136
0g-1575840g^2+3360g^3)\omega^2
}\normalsize

In all the cases where $r=1$ (i.e. computing $\Delta\spr s._r$)
or the pencil is in $\P^r$, these formulas
match those previously obtained by Cotteril \cite{cotteril2}
by different methods (which don't seem to generalize readily
to the case of higher base dimension).
\end{example}
\begin{example}[beyond pencils]
Macnodal is applicable in any base dimension, though
in the case $\dim(B)>1$, there are apparently no examples in the
literature to compare to. Here is a 2-dimensional example.
Let $X/B$ let a 2-parameter family, $E$ a rank-6 bundle over $B$
with Chern classes $d_i$, and
$\phi: E\to\Lambda_5(L)$ a map. The locus
where $\phi$ has rank at most 3,
 is enumerated by
\small
{\vskip4ex\noindent\font\macnodalfont=cmtt10\macnodalfont\rightskip=0ptplus1fil \chardef\macnodalchar=77\macnodalchar\chardef\macnodalchar=97\macnodalchar\chardef\macnodalchar=99\macnodalchar\chardef\macnodalchar=110\macnodalchar\chardef\macnodalchar=111\macnodalchar\chardef\macnodalchar=100\macnodalchar\chardef\macnodalchar=97\macnodalchar\chardef\macnodalchar=108\macnodalchar\chardef\macnodalchar=62\macnodalchar~\discretionary{}{}{}\chardef\macnodalchar=99\macnodalchar\chardef\macnodalchar=99\macnodalchar\chardef\macnodalchar=40\macnodalchar\chardef\macnodalchar=68\macnodalchar\chardef\macnodalchar=101\macnodalchar\chardef\macnodalchar=108\macnodalchar\chardef\macnodalchar=116\macnodalchar\chardef\macnodalchar=97\macnodalchar\chardef\macnodalchar=40\macnodalchar\chardef\macnodalchar=100\macnodalchar\chardef\macnodalchar=105\macnodalchar\chardef\macnodalchar=109\macnodalchar\chardef\macnodalchar=95\macnodalchar\chardef\macnodalchar=98\macnodalchar\chardef\macnodalchar=61\macnodalchar\chardef\macnodalchar=50\macnodalchar\chardef\macnodalchar=44\macnodalchar\chardef\macnodalchar=50\macnodalchar\chardef\macnodalchar=44\macnodalchar\chardef\macnodalchar=50\macnodalchar\chardef\macnodalchar=44\macnodalchar\chardef\macnodalchar=50\macnodalchar\chardef\macnodalchar=41\macnodalchar\chardef\macnodalchar=41\macnodalchar\vskip0ex}\vskip2explus2ex\noindent{\leftskip=0ptplus1fil\rightskip=\leftskip\lineskiplimit=1000pc\lineskip=2explus2ex\hbox{$\displaystyle (144-204d+106d^2-12d^2g-24d^3+2d^4+84dg-156g+12g^2)\int\limits_{B_0}d_1^2$}\hskip1em\hbox{$\displaystyle +(1152-708d+156d^2-12d^3+48dg-252g)\int\limits_{B_0}d_1\pi_*(L^2)$}\hskip1em\hbox{$\displaystyle +(1152-384d+36d^2-60g)\int\limits_{B_0}d_1\pi_*(L\omega )$}\hskip1em\hbox{$\displaystyle +(402-90d+6d^2-12g)\int\limits_{B_0}d_1\pi_*(\omega ^2)$}\hskip1em\hbox{$\displaystyle +(216-258d+107d^2-12d^2g-18d^3+d^4+108dg-216g)\int\limits_{B_0}d_2$}\hskip1em\hbox{$\displaystyle +(387-132d+12d^2-12g)\int\limits_{B_0}\pi_*(L^2)^2$}\hskip1em\hbox{$\displaystyle +(378-60d)\int\limits_{B_0}\pi_*(L^2)\pi_*(L\omega )$}\hskip1em\hbox{$\displaystyle +(92-12d)\int\limits_{B_0}\pi_*(L^2)\pi_*(\omega ^2)$}\hskip1em\hbox{$\displaystyle +(-2778+582d-30d^2+48g)\int\limits_{B_0}\pi_*(L^2\omega )$}\hskip1em\hbox{$\displaystyle +(-1668+624d-84d^2+4d^3-24dg+168g)\int\limits_{B_0}\pi_*(L^3)$}\hskip1em\hbox{$\displaystyle +51\int\limits_{B_0}\pi_*(L\omega )^2$}\hskip1em\hbox{$\displaystyle +24\int\limits_{B_0}\pi_*(L\omega )\pi_*(\omega ^2)$}\hskip1em\hbox{$\displaystyle +(-1962+266d-6d^2+12g)\int\limits_{B_0}\pi_*(L\omega ^2)$}\hskip1em\hbox{$\displaystyle +3\int\limits_{B_0}\pi_*(\omega ^2)^2$}\hskip1em\hbox{$\displaystyle +(-540+48d)\int\limits_{B_0}\pi_*(\omega ^3)$}\hskip1em\hbox{$\displaystyle +(-222+66d-6d^2+12g)\int\limits_{B_1}d_1$}\hskip1em\hbox{$\displaystyle +(-76+12d)\int\limits_{B_1}\pi_*(L^2)$}\hskip1em\hbox{$\displaystyle -24\int\limits_{B_1}\pi_*(L\omega )$}\hskip1em\hbox{$\displaystyle -6\int\limits_{B_1}\pi_*(\omega ^2)$}\hskip1em\hbox{$\displaystyle +(207-24d)\int\limits_{B_1}\psi^{(1)}_1$}\hskip1em\hbox{$\displaystyle +(1014-178d+6d^2-12g)\int\limits_{B_1}\theta^*_1(L)$}\hskip1em\hbox{$\displaystyle +3\int\limits_{B_2}1$}\vskip6explus2ex}\noindent
\normalsize
Here $B_i\to B$ is the normalization of the locus of $i$-nodal curves,
$\theta_1\to X_{B_1}$ is the node of the restricted family,
and $\psi_1\spr 1.$ is the sum of the cotangent classes at the
branches of $X$ at $\theta_1$. Also, in integrals over $B_i$,
$\omega$ refers to the dualizing sheaf of the \emph{normalization}
of $X_{B_i}$.
\par
More examples and details about using Macnodal will be given elsewhere.

\end{example}
\newsubsection{Double points}
Let $X/B$ be an arbitrary nodal family and $f:X\to\mathbf P^n$ a morphism. Consider the relative double points of $f$, i.e. double
points on fibres. This locus  on $X\sbr 2._B$ is the
degeneracy locus of a bundle map
\[\phi:(n+1)\O\to\Lambda_2(L), L:=f^*\O(1).\]
By Porteous, the virtual fundamental class of this locus
is given by the Segre class $s_n(\Lambda_2(L)^*)$,
whose pullback on the ordered (= flag) Hilbert scheme $X\scl 2._B$ equals
\eqspl{}{
\sum\limits_{i=0}^n(L_1)^{n-i}(L_2-\Gamma)^i, \Gamma=\Gamma\scl 2..
}
The powers of $\Gamma$ can be evaluated using Corollary \ref{polpowers,m=2}.
Pushing the result down to $X^2_B$ for simplicity yields
\eqsp{
\sum\limits_{i=0}^nL_1^{n-i}L_2^i-
\sum\limits_{i=0}^nL_1^{n-i}(\sum\limits_{j=1}^i
(\Gamma[\omega^{j-1}]+\sum\limits_{s,k}
\delta_{s*}(\psi_x^{j-2-k}\psi_y^k))
L^{i-j})
}
To describe the direct image of this on $B$, we need some notation. Recall that $\kappa_j=\pi_*(\omega^{j+1})$. Extending this, we may set
\eqspl{}{
\kappa_j(L)=\pi_*(L^{j+1}), \kappa_{i,j}(L,M)=
\pi_*(L^{i+1}M^{j+1}).
} Note that in our case $\kappa_j(L)$ may be interpreted
as the class of the locus of curves meeting a generic $\P^{n-j}$. Also,
for each boundary datum $(T_s, \delta_s, \theta_s)$, $T_s$
admits a map to $\P^n$ via either the $x$ or $y$-section
(the two maps are the same), via which we can pull back $L^j$,
which corresponds to the locus of boundary curves whose node
$\theta_s$ meets $\P^{n-j}$.
 Then pushing the above down to $B$ yields
 the following
 \begin{prop}
Notations as above, the virtual class of the locus
of fibres not embedded by $f$ is
\eqsp{
m_{2,B}=\half(-1)^n(\sum\limits_{i=1}^{n-1}\kappa_{i-1}(L)
\kappa_{n-i-1}(L)-\kappa_{n-j-1, j-2}(L,\omega)
+\sum\limits_{s,k}
\delta_{s*}(L^{n-j}\psi_x^{j-2-k}\psi_y^k))\qed
}
\end{prop}
More generally, for any smooth variety $Y$ of dimension $n$
 and map $f:X\to Y$, one can use the double-point formula of \cite{R}, Th. 3.3ter, p. 1208, to evaluate the
class of the double-point locus in $X^2_B$ in terms of the diagonal
class $\Delta_Y$ on $Y\times Y$ as
\eqspl{m2}{
2m(f)_{X^2_B}&=(f^2)^*(\Delta_Y)-\sum\limits_{i\geq 1}(-\Gamma^i)c_{n-i}(T_Y)\\
&=(f^2)^*(\Delta_Y)-\sum\limits_{i\geq 1}(
-\Gamma[\omega^{i-1}]+\half\sum\limits_{s,j}\delta_{s*}(\psi_x^{i-j-3}\psi_y^j))
c_{n-i}(T_Y)
} (here each boundary term corresponding to a node $\theta$
is embedded in the diagonal in $X^2_B$ via $\theta$).
Applying this set-up to the case $L=\omega, Y=\P(\bE)$,
$\bE$ being the Hodge bundle, and replacing $T_Y$ and $\Delta_Y$
be their relative analogues over $B$, we note that
\eqspl{}{
[\Delta_{Y/B}]=\{\frac{c(\bE^*)}{
(1-[L_1])(1-[L_2])}\}_{g-1}, L_i=p_i^*\O(1),\\
c_i(T_{Y/B})=\sum(-1)^k\binom{g-k}{i-k}\lambda_kL^{i-k},
\lambda_k=c_k(\bE).
}
Applying this to \eqref{m2}, we note that $L$ pulls back to
$\omega$, which meets each $\theta$ trivially, so we obtain
\eqspl{}{
2m_{X^2_B}(f)=\sum (-1)^k\lambda_k[\omega^i,\omega^{g-1-i-k}]
+\sum\limits_{i\geq\max(1,k)}
(-1)^k\binom{g-k}{n-i-k}\lambda_k\Gamma[\omega^{g-2-k}]\\
-
\half\sum\limits_{i\geq 1, s,j}(-1)^{g-i-1}
\delta_{s*}(\psi_x^{i-j-3}\psi_y^j))
\lambda_{g-i-1}\\
=\sum (-1)^k\lambda_k[\omega^i,\omega^{g-1-i-k}]
+(2^{g-1}-1)\Gamma[\omega^{g-2}]-\sum\limits_{k=1}^{g-1}
(-1)^k2^{g-1-k}\lambda_k\Gamma[\omega^{g-2-k}]\\
-
\half\sum\limits_{i\geq 1, s,j}(-1)^{g-i-1}
\delta_{s*}(\psi_x^{i-j-3}\psi_y^j))
\lambda_{g-i-1}
}
Multiplying by $\omega_1$ and projecting to $B$ we obtain (compare
\cite[\S 7]{Mu}):
\eqspl{}{
2(2g-2)m_2=\sum (-1)^k\lambda_k\kappa_{i-1}\kappa_{g-2-i-k}
+(2^{g-1}-1)\kappa_{g-2}+\sum\limits_{k=1}^{g-1}
(-1)^k2^{g-1-k}\lambda_k\kappa_{g-2-k}
} where we set $\kappa_0=2g-2$ for simplicity.
This formula is correct over the locus $\frak M_g\cup\Delta_0^0$
of curves with at most 1 nonseparating node, but
breaks down over the curves with a separating node or a separating pair of nodes,
because the naive notion of canonical curve
in $\P^{g-1}$ is ill-behaved and requires substantial
modification. See \cite{grd} for some work in this direction.
\subsection{ Hyperelliptic locus in genus 3}
Our purpose here is to compute the fundamental class of the
closure $\he$ of the locus of smooth hyperelliptic curves of
genus 3 in the stable moduli $\overline{\mathcal M}_3$. Consider a big family $X/B$ parametrized by the locally closed subscheme of a suitable Hilbert scheme of curves in
some $\P^N$ corresponding to stable curves. $B$ has the two boundary divisors $\delta_0, \delta_1$, and we have
\[X_{\delta_1}=X_1\bigcup\limits_{p_1\leftrightarrow p_2} X_2\]
with each $X_i$ of genus $i$.
Let $L=\omega_{X/B}(-X_1)$ where $X_1/\delta_1$ is the obvious divisor on $X$ (of
relative genus 1),
and let $\be_1=\pi_*(L)$.
The fibre of $\be_1$ over a curve in $\delta_1$ corresponds
to the complete linear system that is the hyperelliptic system
on the genus 2 part and the pencil $|2p_1|$ on the genus-1 part.
Then we have an exact sequence
\[\exseq{\be_1}{\be}{\pi_*(\omega_{X_1/\delta_1}\otimes\O_{p_1})},\]
hence $c(\be_1)=c(\be)(1-\delta_1)$ and in particular
\[c_1(\be_1)=\lambda_1-\delta_1.\]
Also, it is easy to see that
\[\pi_*(L^2)=\pi_*(\omega^2)-3\delta_1=\kappa_1-3\delta_1,
\pi_*(L.\omega)=\kappa_1-\delta_1.\]
Now consider the rank-1 locus of the natural map over
$W^2(X/B)$:
\[ \be_1\to\Lambda_2(L).\]
By Porteous, the fundamental class of this locus
 in $W^2(X/B)$ is computed by
 \eqsp{
 \Delta&=\{\frac{c(\be_1)}{c(\Lambda_2(L))}\}_2\\
 &=-(\lambda_1-\delta_1)(L_1+L_2-\Gamma\spr 2.)+
 L_1^2+(L_2-\Gamma\spr 2.)^2+L_1(L_2-\Gamma\spr 2.)
 }
Because this locus maps to $B$ with 1-dimensional fibres,
we multiply this by $L_1$ and then project to $B$.
Using Mumford's relation $\kappa=12\lambda_1-\delta$,
we compute easily that
\[\pi_*(L_1\Delta)=36 \lambda_1-4\delta_0-6\delta_1.\]
Now the locus $L_1\Delta$ clearly covers the hyperelliptic locus
4 times. Additionally, there are three loci over $\delta_1$:\par
 (A) the hyperelliptic pencil on $X_2$;\par
  (B) the pencil $|2p_1|$
on $X_1$;\par
(C) $X_1\times p_2'$ where $p_2'$ is the image of $p_2$
 under the hyperelliptic involution.\par
Each of these, intersected with $L_1$, maps to $\delta_1$
with degree 2, for a total of 6.
 Therefore
\[[\he]=\frac{1}{4}(36\lambda_1-4\delta_0-12\delta_1)=
9\lambda_1-\delta_0-3\delta_1,\]
a formula first obtained by Harris and Mumford \cite{har-mum}.
In higher genus, extra components of excessive dimension
appear, requiring substantial modifications (see \cite{grd}).
\bibliographystyle{amsplain}
\bibliography{mybib}
\end{document}